\documentclass{article}
\usepackage[cm]{fullpage}
\usepackage{graphicx} 

\usepackage{amsmath, amssymb, amsfonts, amsthm}
\usepackage{mathdots}
\usepackage{subcaption}
\usepackage{xcolor}
\usepackage{enumitem}
\usepackage[numbers]{natbib}

\usepackage{algorithm}
\usepackage{algpseudocode}

\usepackage{hyperref}
\usepackage{colortbl}
\definecolor{niceblue}{rgb}{0.0,0.19,0.56}
\hypersetup{colorlinks,linkcolor={blue},citecolor={niceblue},urlcolor={blue}}

\usepackage[capitalise]{cleveref}

\usepackage{todonotes}
\makeatletter
\renewcommand*\env@matrix[1][*\c@MaxMatrixCols c]{%
  \hskip -\arraycolsep
  \let\@ifnextchar\new@ifnextchar
  \array{#1}}
\makeatother

\usepackage{array,booktabs,multirow,makecell,ragged2e,tabularx,rotating,pdflscape}
\newcolumntype{L}[1]{>{\RaggedRight\arraybackslash}p{#1}}
\newcolumntype{Y}{>{\RaggedRight\arraybackslash}X}

\newif\ifAutoRotateSEGComparisonTable
\AutoRotateSEGComparisonTablefalse
\newcommand{\SEGComparisonTableWidth}{%
  \ifAutoRotateSEGComparisonTable
    \linewidth
  \else
    \textheight
  \fi
}

\usepackage{nicefrac}

\allowdisplaybreaks

\title{On Same-Sample and Independent-Sample Stochastic Extragradient \\ for Monotone Variational Inequalities}

\author{TaeHo Yoon\textsuperscript{*}
\and Nicolas Loizou\textsuperscript{*}
}

\date{}

\newcommand{\imu}{\mathbf{i}}

\newcommand{\cut}[1]{{}}

\newcommand{\cC}{{\mathcal{C}}}
\newcommand{\cD}{{\mathcal{D}}}
\newcommand{\cE}{{\mathcal{E}}}
\newcommand{\cF}{{\mathcal{F}}}

\newcommand{\cO}{{\mathcal{O}}}

\newcommand{\cX}{{\mathcal{X}}}
\newcommand{\cY}{{\mathcal{Y}}}
\newcommand{\cZ}{{\mathcal{Z}}}

\newcommand{\order}[1]{{\mathcal{O}\left(#1\right)}}

\usepackage[bb=boondox]{mathalfa}

\usepackage{xparse}
\DeclareFontFamily{U}{ntxmia}{}
\DeclareFontShape{U}{ntxmia}{m}{it}{<-> ntxmia }{}
\DeclareFontShape{U}{ntxmia}{b}{it}{<-> ntxbmia }{}
\DeclareSymbolFont{lettersA}{U}{ntxmia}{m}{it}
\SetSymbolFont{lettersA}{bold}{U}{ntxmia}{b}{it}

\ExplSyntaxOn
\NewDocumentCommand{\varmathbb}{m}
 {
  \tl_map_inline:nn { #1 }
   {
    \use:c { varbb##1 }
   }
 }
\tl_map_inline:nn { ABCDEFGHIJKLMNOPQRSTUVWXYZ }
 {
  \exp_args:Nc \DeclareMathSymbol{varbb#1}{\mathord}{lettersA}{\int_eval:n { `#1+67 }}
 }
\exp_args:Nc \DeclareMathSymbol{varbbk}{\mathord}{lettersA}{169}
\ExplSyntaxOff

\renewcommand{\Re}{\operatorname{Re}}

\newcommand{\proj}{\Pi}

\newcommand{\reals}{\mathbb{R}}

\newcommand{\expec}[2]{\mathbb{E}_{#2}\left[#1\right]}
\newcommand{\condexp}[2]{\mathbb{E}\left[#1\,\middle|\,#2\right]}
\newcommand{\prob}[1]{\mathrm{Prob}\left[#1\right]}
\newcommand{\ones}{\mathbf{1}}

\newcommand{\inprod}[2]{\left\langle #1,#2 \right\rangle}
\newcommand{\sqnorm}[1]{\left\| #1 \right\|^2}

\newcommand{\norm}[1]{\left\|#1\right\|}

\newcommand{\pr}[1]{\left( #1 \right)}

\newcommand{\halfitr}{\hat{x}_k}
\newcommand{\halfxi}{\xi_{\hat{i}_k}}

\newcommand{\extstep}{\hat{\gamma}_k}
\newcommand{\err}{\mathrm{Err}}

\usepackage{thmtools}
\definecolor{shadecolor}{gray}{0.9}
\declaretheoremstyle[
headfont=\normalfont\bfseries,
notefont=\mdseries, notebraces={(}{)},
bodyfont=\normalfont,
postheadspace=0.5em,
spaceabove=\topsep,
mdframed={
  skipabove=8pt,
  skipbelow=8pt,
  hidealllines=true,
  backgroundcolor={shadecolor},
  innerleftmargin=4pt,
  innerrightmargin=4pt}
]{shaded}

\declaretheorem[style=shaded,within=section]{definition}
\declaretheorem[style=shaded,sibling=definition]{theorem}
\declaretheorem[style=shaded,sibling=definition]{proposition}
\declaretheorem[style=shaded,sibling=definition]{assumption}
\declaretheorem[style=shaded,sibling=definition]{corollary}

\declaretheorem[style=shaded,sibling=definition]{lemma}

\newcommand{\attn}[1]{{\color{red}#1}}

\begin{document}

\renewcommand\thefootnote{\fnsymbol{footnote}}
\footnotetext[1]{Department of Applied Mathematics and Statistics \& Mathematical Institute for Data Science, Johns Hopkins University.
\texttt{\{tyoon7, nloizou\}@jhu.edu}}

\renewcommand\thefootnote{\arabic{footnote}}
\setcounter{footnote}{0}

\maketitle

\begin{abstract}

We study stochastic extragradient (SEG) methods for solving monotone variational inequality problems (VIPs) over a feasible set $\cX$.
Although extragradient is a foundational algorithm for VIPs and its deterministic convergence theory is well developed, its stochastic counterpart remains less understood.
Most existing analyses focus on \emph{independent-sample SEG (I-SEG)} and assume either that the domain $\cX$ is compact or that the variance of the stochastic operator is uniformly bounded.
The behavior of \emph{same-sample SEG (S-SEG)}, a natural variant with materially different properties, has received far less attention.
In this work, we address these gaps in the literature.
We first show that S-SEG is sensitive to samplewise Lipschitz parameters: mean Lipschitzness and bounded variance alone do not ensure convergence, even on a compact set.
Then, for possibly unbounded domains, we establish a high-probability restricted-gap convergence for each SEG variant under a relaxed set of assumptions,
and show that certain fundamental improvements to these results are impossible in general.
Finally, we show that a known asymmetric double step-size selection that guarantees almost sure last-iterate convergence for I-SEG can fail for S-SEG: there exists a stochastic monotone VIP for which S-SEG diverges almost surely even under the modified step-sizes.

\vspace{.1cm}
\noindent
\textbf{Keywords\phantom{.} } Variational inequality $\cdot$ Stochastic optimization $\cdot$ Extragradient $\cdot$ Gap function $\cdot$ Almost sure convergence

\vspace{.1cm}
\noindent
\textbf{Mathematics Subject Classification\phantom{.} } 65K15 $\cdot$ 90C33 $\cdot$ 62L20 $\cdot$ 90C15 $\cdot$ 49J40 $\cdot$ 47H05
\end{abstract}

\tableofcontents

\section{Introduction}

A variational inequality problem (VIP)\footnote{
    We state the Minty (or ``weak'') form of VIP \cite{mintyGeneralizationDirectMethod1967} here, while formulation by Stampacchia (``strong'') VIP \cite{stampacchiaFormesBilineairesCoercitives1964}:  
    $\underset{x\in\cX}{\text{find}} \,\, \inprod{F(x)}{u - x} \ge 0 , \,\, \forall u \in \cX$, appeared earlier.
    The two formulations are equivalent if $F$ is monotone and continuous.
} \cite{stampacchiaFormesBilineairesCoercitives1964, mintyGeneralizationDirectMethod1967}, stated as
\begin{align}
\label{eqn:VIP}
\begin{array}{cc}
    \underset{x\in\cX}{\text{find}} & \inprod{F(u)}{x - u} \le 0 , \,\, \forall u \in \cX
\end{array}
\end{align}
where $\cX \subseteq \reals^d$ is a closed convex domain, is a general problem class that arises naturally in 
optimization and game theory \cite{facchineiFiniteDimensionalVariationalInequalities2003}.
For a VIP, the objective operator $F\colon \reals^d \to \reals^d$ is typically assumed to be monotone.
For example, $F = \nabla f$ may be the gradient of a differentiable and convex function $f \colon \reals^d \to \reals$, 
or $F(x^{(1)}, x^{(2)}) = \left( \nabla_{x^{(1)}} f(x^{(1)}, x^{(2)}) , -\nabla_{x^{(2)}} f(x^{(1)}, x^{(2)}) \right)$ may be the saddle operator associated with a 
differentiable and convex-concave function $f\colon \reals^{d_1} \times \reals^{d_2} \to \reals$ \cite{rockafellarConvexAnalysis1970}.
In this sense, a monotone VIP encodes the first-order optimality conditions for constrained convex minimization and convex-concave minimax optimization, 
and can be viewed as their natural generalization.

In modern optimization, objectives are often constructed from large datasets.
For VIPs, a corresponding theoretical model is the objective operator $F$ with the finite sum structure
\begin{align}
\label{eqn:objective-F}
    F(x) = \frac{1}{n} \sum_{i=1}^n F(x; \xi_i) .
\end{align}
When $d$ and $n$ are large, the standard paradigm is to use stochastic first-order algorithms, 
which update the iterates using sample operator evaluations $F(\cdot; \xi_i)$ rather than evaluations of the full-batch operator $F$.
This will be the setting of interest for our work.

A variety of algorithms for solving VIP~\eqref{eqn:VIP} or its unconstrained case (which finds $x \in \reals^d$ such that $F(x) = 0$) exist in the literature.
Korpelevich's celebrated Extragradient (EG) algorithm \cite{Korpelevich1976_extragradient} is one of the algorithmic foundations for solving VIP.
It has been generalized to the non-Euclidean Mirror-Prox algorithm \cite{nemirovskiProxmethodRateConvergence2004}, whose stochastic extension \cite{juditskySolvingVariationalInequalities2011} has become a standard reference for stochastic VIP.
More recent works have actively studied (stochastic) EG and its extensions in the context of machine learning \citep{GidelBerardVignoudVincentLacoste-Julien2019_variational, mertikopoulosOptimisticMirrorDescent2019, golowichTightLastiterateConvergence2020,golowichLastIterateSlower2020, HsiehIutzelerMalickMertikopoulos2020_explore,MishchenkoKovalevShulginRichtarikMalitsky2020_revisiting,MokhtariOzdaglarPattathil2020_unified,MokhtariOzdaglarPattathil2020_convergence,antonakopoulosSiftingNoiseUniversal2021,GorbunovLoizouGidel2022_extragradient,GorbunovBerardGidelLoizou2022_stochastic,chaeStochasticExtragradientFlipflop2024,choudhury2025extragradient}.
Complementary lines of work seek single-call algorithms or different algorithmic structures.
Popov's method is an early alternative to EG requiring only one new operator evaluation per iteration \citep{popovModificationArrowHurwiczMethod1980}, whose modern optimistic mirror descent interpretation was later proposed and popularized in the context of zero-sum games \citep{rakhlinOptimizationLearningGames2013, mertikopoulosOptimisticMirrorDescent2019} and machine learning \citep{ChavdarovaGidelFleuretLacoste-Julien2019_reducing,hsiehConvergenceSinglecallStochastic2019,golowichTightLastiterateConvergence2020,gorbunovLastiterateConvergenceOptimistic2022,ChoudhuryGorbunovLoizou2023_singlecall}.
The golden ratio algorithm considered in \citep{malitskyGoldenRatioAlgorithms2020, alacaogluGoldenRatioVariational2023} also uses one operator evaluation per iteration and admits an adaptive line-search-free variant.
Lightweight second-order algorithms such as Hamiltonian gradient descent or consensus optimization achieve linear convergence under favorable conditions
\citep{MeschederNowozinGeiger2017_numerics, liangInteractionMattersNote2019,loizouStochasticHamiltonianGradient2020,AbernethyLaiWibisono2021_lastiterate,LoizouBerardGidelMitliagkasLacoste-Julien2021_stochastic}.
Recently, a number of works adapted Halpern iteration for fixed-point problems \citep{Halpern1967_fixed, Lieder2021_convergence} to minimax,  monotone inclusion problems and VIPs, and achieved accelerated residual norm convergence \citep{Diakonikolas2020_halpern, YoonRyu2021_accelerated,LeeKim2021_fast,Tran-DinhLuo2021_halperntype, CaiZheng2023_accelerated,yoonAcceleratedMinimaxAlgorithms2025,botExtragradientMethodFlexible2026}.
This inspired the development of distinct interpretations and new mechanisms for acceleration in these problem classes \citep{tran-dinhHalpernsFixedpointIterations2024,botFastOptimisticGradient2025,sedlmayerFastOptimisticMethod2023,yoonOptimalAccelerationMinimax2024,tran-dinhExtragradienttypeMethods2024,yoonHinvarianceTheoryComplete2025,yoonTheoryCompositionDuality}.
Some of these accelerated methods have been applied to stochastic VIPs or monotone inclusion problems via variance reduction \citep{cai2022stochastic, alacaogluComplexitySimplePrimaldual2025,tran-dinhVFOSAVariancereducedFast2025}.

Among all the aforementioned algorithms, in this paper, we focus on stochastic Extragradient (SEG), which is not yet fully studied in spite of being a fundamental algorithm for stochastic VIPs.

\subsection{Extragradient and its stochastic variants}

Although gradient descent is foundational in first-order optimization, its analogue for a general monotone VIP, defined by the update rule 
\begin{align}
\label{eqn:operator-descent}
    x_{k+1} = \proj_\cX \left( x_k - \gamma_k F(x_k) \right) ,
\end{align}
where $\proj_\cX$ denotes projection onto the closed convex set $\cX$ and $\gamma_k$ are the step-sizes, may not converge.  
A standard counterexample is the rotation operator $F \colon \reals^2 \to \reals^2$ given by $F(x^{(1)}, x^{(2)}) = (x^{(2)}, -x^{(1)})$ with $\cX = \reals^2$, 
arising from the two-player minimax problem with bilinear objective $f(x^{(1)}, x^{(2)}) = x^{(1)} x^{(2)}$.
For this operator, every iteration of \eqref{eqn:operator-descent} causes the iterates to further deviate from the unique solution $x_\star = (0,0)$.
While certain convergence results for \eqref{eqn:operator-descent} exist, they require additional conditions such as boundedness of the operator \cite{nedicSubgradientMethodsSaddlePoint2009} or 
stability of the optimum \cite{liangInteractionMattersNote2019, BeznosikovGorbunovBerardLoizou2023_stochastic}; even then, convergence can be slow \cite{MeschederNowozinGeiger2017_numerics}.
Consequently, there is broad consensus that additional mechanisms are required to solve \eqref{eqn:VIP} efficiently \cite{yadavStabilizingAdversarialNets2018,GidelBerardVignoudVincentLacoste-Julien2019_variational,gidelNegativeMomentumImproved2019,mertikopoulosOptimisticMirrorDescent2019,YoonRyu2021_accelerated,yoonOptimalAccelerationMinimax2024,zhengDissipativeGradientDescent2024,LeeChoYun2024_fundamental,shugart2025negative}.

Extragradient (EG) \cite{Korpelevich1976_extragradient} is a classical algorithm proposed for solving saddle point problems (minimax optimization) using two operator evaluations per iteration.
In the language of VIPs, it has the update rule
\begin{align*}
    \halfitr & = \proj_\cX \left( x_k - \extstep F(x_k) \right) \\
    x_{k+1} & =  \proj_\cX \left( x_k - \gamma_k F(\halfitr) \right) .
\end{align*}
The first line is commonly called the extrapolation step, and the second the update step; the latter uses the operator evaluated at the extrapolated iterate $\halfitr$.
Despite its simplicity, EG is remarkably effective. For a VIP with a monotone and Lipschitz continuous operator $F$, it achieves an $\cO(1/k)$ convergence rate with respect to the gap function, one of the standard measures for this problem class \cite{nemirovskiProxmethodRateConvergence2004}. This rate is optimal up to a constant \cite{Nemirovsky1992_informationbased,nemirovskiProxmethodRateConvergence2004}.
A wide range of algorithms now used for minimax optimization and VIPs can be viewed broadly as generalizations or extensions of EG \cite{popovModificationArrowHurwiczMethod1980,tsengModifiedForwardbackwardSplitting2000,SolodovSvaiter1999_hybrid,rakhlinOptimizationLearningGames2013,mertikopoulosOptimisticMirrorDescent2019}.
Consequently, the convergence theory of EG has been studied extensively and is now well established \citep{MonteiroSvaiter2010_complexity, AzizianMitliagkasLacoste-JulienGidel2020_tight,MokhtariOzdaglarPattathil2020_unified,MokhtariOzdaglarPattathil2020_convergence, GorbunovLoizouGidel2022_extragradient}.

\begin{figure}[H]
    \centering
    \begin{subfigure}[t]{0.35\textwidth}
        \centering
        \includegraphics[width=\linewidth]{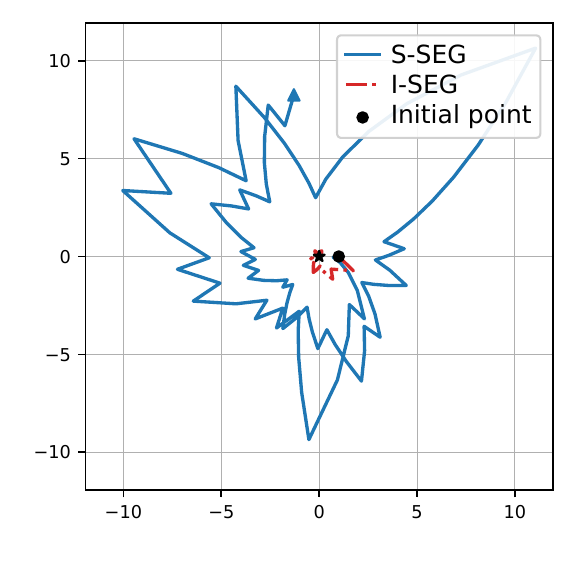}
        \caption{$(\rho,\nu)=(2,0)$}
        \label{fig:sseg-counterexample-plane-trajectory}
    \end{subfigure}
    \hspace{1.5cm}
    \begin{subfigure}[t]{0.35\textwidth}
        \centering
        \includegraphics[width=\linewidth]{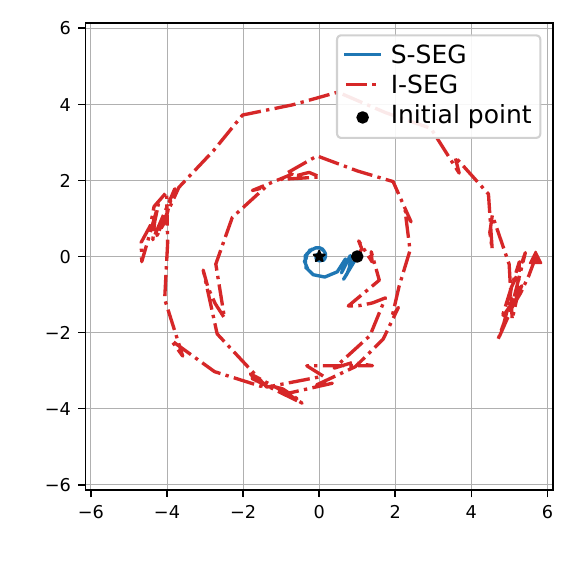}
        \caption{$(\rho,\nu)=(0,2)$.}
        \label{fig:iseg-counterexample-plane-trajectory}
    \end{subfigure}
    \caption{Trajectories of S-SEG and I-SEG on two different two-dimensional stochastic monotone VIPs defined in \cref{lemma:stochastic-vip-counterexample}: $J = \begin{bmatrix} 0 & -1 \\ 1 & 0 \end{bmatrix}$, and $F(x;\xi_1) = (\rho I + (1+\nu)J)x$ and $F(x;\xi_2) = (-\rho I + (1-\nu)J)x$ for $x \in \cX = \reals^2$.
    The two problems are constructed via different choices of $\rho$ and $\nu$, and both versions of SEG are run with step-size $\extstep = \gamma_k = \frac{1}{\sqrt{k+25}}$ for $N=100$ iterations in both cases. In (a), S-SEG dynamics is unstable while I-SEG converges to the solution $(0,0)$, indicated with the ``$\star$'' marker, while in (b), the roles of S-SEG and I-SEG are reversed. \emph{Thus, neither sampling rule uniformly dominates the other.}
    }
    \label{fig:seg-counterexample-plane-trajectories}
\end{figure}

Given the success of EG for deterministic VIPs, a natural stochastic counterpart for objective operators of the form \eqref{eqn:objective-F} is stochastic extragradient (SEG), for which we have two variants:
\begin{center}
\vspace{-0.3cm}
\begin{minipage}{0.45\textwidth}
\begin{align}
\label{eqn:I-SEG}
\tag{I-SEG}
\begin{aligned}
    \halfitr & = \proj_\cX \left( x_k - \extstep F(x_k; \xi_{i_k}) \right) \\
    x_{k+1} & =  \proj_\cX \left( x_k - \gamma_k F(\halfitr; \halfxi) \right) 
\end{aligned}    
\end{align}
\end{minipage}
\hfill
\begin{minipage}{0.45\textwidth}
\begin{align}
\label{eqn:S-SEG}
\tag{S-SEG}
\begin{aligned}
    \halfitr & = \proj_\cX \left( x_k - \extstep F(x_k; \xi_{i_k}) \right) \\
    x_{k+1} & =  \proj_\cX \left( x_k - \gamma_k F(\halfitr; \xi_{i_k}) \right) .
\end{aligned}    
\end{align}
\end{minipage}
\end{center}
Here \textit{independent-sample SEG} (\ref{eqn:I-SEG}) selects independent random indices $i_k, \hat{i}_k$ from $[n] = \{1,\dots,n\}$ in the extrapolation and update steps,
and there is a large body of convergence theory for I-SEG and its extensions \citep{juditskySolvingVariationalInequalities2011, yousefianStochasticMirrorproxAlgorithms2018, kannanOptimalStochasticExtragradient2019, mertikopoulosOptimisticMirrorDescent2019, GidelBerardVignoudVincentLacoste-Julien2019_variational, ChavdarovaGidelFleuretLacoste-Julien2019_reducing, HsiehIutzelerMalickMertikopoulos2020_explore, antonakopoulosSiftingNoiseUniversal2021, beznosikovDecentralizedLocalStochastic2022}.
However, the dynamics of SEG depend heavily on how the stochastic indices are chosen.
In particular, \emph{same-sample SEG} (\ref{eqn:S-SEG}) which uses a shared sample $\xi_{i_k} = \halfxi$ across the two steps, can follow trajectories entirely different from those of I-SEG on different stochastic monotone VIPs (Figure~\ref{fig:seg-counterexample-plane-trajectories}).
In our view, S-SEG is at least as natural as I-SEG and may even be preferable in practice because it does not require two independent operator samples at each iteration.
However, only a small body of prior work has focused on S-SEG \citep{MishchenkoKovalevShulginRichtarikMalitsky2020_revisiting, GorbunovBerardGidelLoizou2022_stochastic}, and a relatively limited range of theoretical guarantees have been established. 
Thus, the theoretical understanding of S-SEG remains less developed than that of I-SEG, and the precise distinctions between the two variants remain underexplored.

\subsection{Notion of convergence for VIP}

For VIP, the following (restricted) gap function \cite{larssonClassGapFunctions1994,nesterovDualExtrapolationIts2007} is a standard measure of convergence used for deriving quantitative rates:
\begin{align}
\label{eqn:gap-function}
    \err_\cZ (x) = \max_{u\in \cX \cap \cZ} \inprod{F(u)}{x-u} \ge 0 .
\end{align}
Here $\cZ$ is a compact set used to enforce the finiteness of \eqref{eqn:gap-function}.
In cases where $\cX$ is already compact, we take $\cZ = \cX$ and simply write $\err(x) = \err_\cX(x)$, and otherwise, when $\cX$ is unbounded, we take the ball $\cZ = B(y,D)$ for some $y\in \cX$ and radius $D>0$.
Assuming that $B(y,D)$ contains at least one solution to \eqref{eqn:VIP}, for $x$ in the interior of $B(y,D)$, we have $\err_{B(y,D)}(x) = 0$ if and only if $x$ is a solution to \eqref{eqn:VIP} \cite[Lemma~1]{nesterovDualExtrapolationIts2007}.
The seminal result of \cite{nemirovskiProxmethodRateConvergence2004} showed $\err(\overline{\hat{x}}_N) = \cO(1/N)$ along the ergodic average of extrapolated iterates $\overline{\hat{x}}_N$ for EG, 
and \cite{juditskySolvingVariationalInequalities2011} proved $\expec{\err(\overline{\hat{x}}_N)}{} = \cO(1/\sqrt{N})$ for I-SEG.

Another type of convergence of interest is whether the iterate sequence $x_k$ converges to a VIP solution $x_\star$ almost surely with respect to the randomness in the algorithm, i.e., whether $x_k \to x_\star$ with probability $1$.
Here, we typically focus on the last iterate rather than the averaged iterate.
For deterministic EG, the original work \cite{Korpelevich1976_extragradient} already proved the last-iterate convergence, 
and for I-SEG, prior works including \cite{kannanOptimalStochasticExtragradient2019,mertikopoulosOptimisticMirrorDescent2019,antonakopoulosSiftingNoiseUniversal2021,HsiehIutzelerMalickMertikopoulos2020_explore} have provided almost-sure convergence guarantees.
On the other hand, whether the analogous results hold for S-SEG is mostly unknown, to the best of our knowledge.

\ifAutoRotateSEGComparisonTable
\begin{landscape}
\begin{table}[p]
\else
\begin{sidewaystable}[p]
\fi
\centering
\scriptsize
\setlength{\tabcolsep}{4pt}
\renewcommand{\arraystretch}{1.18}
\caption{Comparison of convergence results for \ref{eqn:I-SEG} and \ref{eqn:S-SEG} considered in this work. New results are highlighted in red. Here, we globally assume that $F$ is monotone and Lipschitz, and $F(\cdot;\xi_i)$ is an unbiased sampling of $F$.
Please see \cref{section:technical-summary} for a detailed explanation of novel results, and \cref{section:prior-work} for a review of existing results.
}
\label{tab:sseg-iseg-comparison}
\begin{tabularx}{\SEGComparisonTableWidth}{@{}L{1.75cm}L{2.35cm}L{2.35cm}L{3.55cm}YY@{}}
\toprule
Result type & Domain & Step-size & Assumptions & I-SEG & S-SEG \\
\midrule
\multirow{5}{1.75cm}{\RaggedRight \makecell[c]{Gap function \\ (Averaged \\ iterate)} }
& \multirow{2}{2.35cm}{\RaggedRight \makecell[l]{$\mathcal{X}$ compact \\ (\cref{section:prior-work}, \\ \cref{section:S-SEG-compact-analysis})} }
& \multirow{2}{2.35cm}{\RaggedRight\makecell[l]{$\gamma_k = \extstep$\\$\sum_{k=0}^N \gamma_k = o(N)$\\$\gamma_N = \omega(1/N)$}}
& \parbox[t]{\linewidth}{%
\cref{assumption:bounded-variance} (I-SEG)\\
\cref{assumption:uniform-lipschitzness} (S-SEG)%
}
& Converges with rate $\cO\left(N^{-1/2}\right)$ with appropriate step-size choice \citep{juditskySolvingVariationalInequalities2011, yousefianStochasticMirrorproxAlgorithms2018}.
& Converges with rate $\cO\left(N^{-1/2}\right)$ with appropriate step-size choice \citep{GidelBerardVignoudVincentLacoste-Julien2019_variational}. \\
\cmidrule(lr){4-6}
& & & \makecell[l]{  \cref{assumption:bounded-variance} but not \\ \cref{assumption:uniform-lipschitzness} }
& 
Converges with rate $\cO\left(N^{-1/2}\right)$ with appropriate step-size choice \citep{juditskySolvingVariationalInequalities2011, yousefianStochasticMirrorproxAlgorithms2018}.
& \makecell[l]{ May converge to $x_\infty \notin \cX_\star$ almost surely \\ \attn{(New: \cref{proposition:same-sample-nonlipschitz-realizations-positive-measure}).} } \\
\cmidrule(lr){2-6}
& \multirow{3}{2.35cm}{\RaggedRight \makecell[l]{$\mathcal{X}$ unbounded \\ (\cref{section:gap-function-general-domain})} }
& \makecell[l]{$\gamma_k = \extstep$\\$\sum_k \gamma_k^2 < \infty$\\
$k\gamma_k \to \infty$}
& \parbox[t]{\linewidth}{%
\cref{assumption:generalized-variance-bound} (I-SEG)\\
\cref{assumption:uniform-lipschitzness} (S-SEG)%
}
& \makecell[l]{ High-probability conditional expectation bound \\ holds \attn{(New: \cref{theorem:independent-sample-SEG-unbounded-domain-high-probability-bound}). } \\ Similar unconditional bound is generally not \\
possible \attn{(New: \cref{theorem:unbounded-gap-high-probability-tightness}). } } 
& \makecell[l]{ High-probability conditional expectation bound \\ holds \attn{(New: \cref{theorem:same-sample-SEG-unbounded-domain-high-probability-bound}). } \\
Similar unconditional bound is generally not \\ 
possible \attn{(New: \cref{theorem:unbounded-gap-high-probability-tightness}). } } 
\\
\cmidrule(lr){3-6}
&
& \makecell[l]{$\gamma_k = \extstep$\\$\sum_k \gamma_k^2 = \infty$}
& \parbox[t]{\linewidth}{%
\cref{assumption:generalized-variance-bound} (I-SEG)\\
\cref{assumption:uniform-lipschitzness} (S-SEG)%
}
& \makecell[l]{ $x_N$, $\overline{\hat{x}}_N$ may both diverge almost surely \attn{(New:} \\  \attn{\cref{theorem:independent-sample-unbounded-divergence-generalized-variance}). }
However, expected gap function \\ converges with stronger \cref{assumption:bounded-variance} \citep{GidelBerardVignoudVincentLacoste-Julien2019_variational}. } 
& \makecell[l]{ $x_N$, $\overline{\hat{x}}_N$ may both diverge almost surely \\ \attn{(New: \cref{theorem:same-sample-unbounded-divergence-nonzero-mean}).} } \\
\cmidrule(lr){3-6}
&
& \makecell[l]{$\sum_k \gamma_k \hat{\gamma}_k = \infty$\\$\sum_k \gamma_k^2 < \infty$\\$\sum_k \gamma_k \hat{\gamma}_k^2 < \infty$}
& \parbox[t]{\linewidth}{%
\cref{assumption:generalized-variance-bound} (I-SEG)\\
\cref{assumption:uniform-lipschitzness} (S-SEG)%
}
& \makecell[l]{ Converges with rate $\cO\left( \left( \sum_{k=0}^N \gamma_k \right)^{-1} \right)$ \\ \attn{(New: \cref{proposition:independent-sample-DSEG-restricted-gap}). } } 
& \makecell[l]{ $x_N, \overline{\hat{x}}_N^\gamma$ may both diverge almost surely \\ \attn{(New: \cref{theorem:same-sample-DSEG-gamma-average-divergence}). } } \\
\midrule
\multirow{2}{1.75cm}{\RaggedRight\makecell[c]{Almost sure\\convergence\\(Last iterate)}}
& \multirow{2}{2.35cm}{\RaggedRight\makecell[l]{ $\mathcal{X}$ unbounded \\ (\cref{section:almost-sure-convergence}) }}
& \makecell[l]{$\sum_k \extstep^2 < \infty$\\$\sum_k \gamma_k^2 < \infty$}
& \parbox[t]{\linewidth}{%
\cref{assumption:quasi-strict-monotone} \\
\cref{assumption:generalized-variance-bound} (I-SEG)\\
\cref{assumption:uniform-lipschitzness} (S-SEG)
}
& $x_k \to x_\star$ almost surely \citep{kannanOptimalStochasticExtragradient2019,mertikopoulosOptimisticMirrorDescent2019}.
& \makecell[l]{ $x_k \to x_\star$ almost surely \attn{(New: \cref{theorem:S-SEG-almost-sure-strict-monotone}).} } \\
\cmidrule(lr){3-6}
&
& \makecell[l]{$\sum_k \gamma_k \hat{\gamma}_k = \infty$\\$\sum_k \gamma_k^2 < \infty$\\$\sum_k \gamma_k \hat{\gamma}_k^2 < \infty$}
& \parbox[t]{\linewidth}{%
\cref{assumption:generalized-variance-bound} (I-SEG)\\
\cref{assumption:uniform-lipschitzness} (S-SEG)
}
& $x_k \to x_\star$ almost surely \citep{HsiehIutzelerMalickMertikopoulos2020_explore}.
& \makecell[l]{ $x_k$ may diverge almost surely \\ \attn{(New: \cref{theorem:same-sample-DSEG-gamma-average-divergence}). } } \\
\bottomrule
\end{tabularx}
\ifAutoRotateSEGComparisonTable
\end{table}
\end{landscape}
\else
\end{sidewaystable}
\fi

\subsection{Summary of main contributions}
\label{section:technical-summary}

Here, we provide a detailed summary of our main technical contributions. 
See Table~\ref{tab:sseg-iseg-comparison}, which serves two complementary purposes:
it provides a compact review of the existing convergence theory for I-SEG and S-SEG,
and it identifies the open questions resolved in this paper.
The comparison is organized along the three axes that govern the behavior of both versions of SEG: (i)
the convergence criterion (restricted gap function convergence of averaged iterates or almost-sure last-iterate convergence), (ii) the boundedness of the domain, and (iii) the step-size selection.
In Table~\ref{tab:sseg-iseg-comparison}, the previously known results are displayed in black text with citations, while the contributions of this paper are
highlighted in red.

The main message is that the choice between independent and same samples in the two versions of SEG is not simply an innocuous implementation detail.
Under symmetric and square-summable step-sizes, the two variants admit a similar high-probability gap function convergence theory on possibly unbounded domains.
Outside this regime, however, their behavior differs sharply: first, S-SEG guarantees rely on uniform Lipschitzness of sample operators, a condition beyond what is needed for I-SEG, and second, an asymmetric step-size choice from \cite{HsiehIutzelerMalickMertikopoulos2020_explore} that restores convergence for I-SEG can still lead S-SEG to diverge.
Furthermore, our positive results---convergence guarantees---are accompanied by negative results---counterexamples---showing that,
in the directions considered here, the assumptions and the probabilistic form of the
results in our theorems cannot generally be improved.
Altogether, our results delineate the boundary between convergence and failure for the two principal variants of SEG.

We summarize the main contributions below, beginning with the results for possibly
unbounded domains, which constitute the technical core of our analysis.

\begin{itemize}[leftmargin=*]
    \item[$\diamond$]  \textbf{A sharp restricted-gap theory on general closed, convex, and possibly
    unbounded domains.}

    Prior work establishes gap-function convergence for I-SEG on unbounded domains under uniformly bounded variance (\cref{assumption:bounded-variance})
    \citep{GidelBerardVignoudVincentLacoste-Julien2019_variational}.
    In contrast, the corresponding unbounded-domain behavior of S-SEG, as well as the
    behavior of I-SEG under the weaker generalized variance
    model (\cref{assumption:generalized-variance-bound}), is not clearly understood.

    For nonincreasing and symmetric step-sizes
    $\extstep=\gamma_k$ satisfying $\sum_{k=0}^{\infty}\gamma_k^2<\infty$ and $k\gamma_k\to\infty$, we prove conditional expected gap function bounds on a high-probability event for both I-SEG and S-SEG, under their respective assumptions
    (Theorems~\ref{theorem:independent-sample-SEG-unbounded-domain-high-probability-bound}
    and~\ref{theorem:same-sample-SEG-unbounded-domain-high-probability-bound}).
    These bounds immediately yield high-probability convergence guarantees via
    \cref{corollary:conditional-expectation-to-high-probability}.
    In particular, for $\gamma_k=\Theta((k+1)^{-\alpha})$ with
    $\frac12<\alpha<1$, the conditional expected gap function decreases at rate
    $\cO(N^{-(1-\alpha)})$.

    We complement these guarantees with two impossibility results via counterexamples.
    First, square summability cannot be removed in general: with the non-square-summable
    choice $\extstep=\gamma_k=1/\sqrt{k+1}$, for both I-SEG and S-SEG, the last-iterate and average-iterate norms may both be unbounded (Theorems~\ref{theorem:same-sample-unbounded-divergence-nonzero-mean} and \ref{theorem:independent-sample-unbounded-divergence-generalized-variance}).
    Second, even with square-summable symmetric step-sizes, the high-probability guarantee cannot generally be replaced with an unconditional guarantee on the expected gap function over a prespecified bounded region, as the limiting radius of the iterates can have unbounded support (\cref{theorem:unbounded-gap-high-probability-tightness}).
    Together, these positive and negative results explain both why conditioning is needed,
    and why the high-probability form of our gap function theory is essentially
    unavoidable.

    \item [$\diamond$] \textbf{Almost-sure last-iterate convergence for \ref{eqn:S-SEG}, and a fundamental
    separation from \ref{eqn:I-SEG}.}

    Prior work establishes almost-sure last-iterate convergence of I-SEG under
    quasi-strict monotonicity; see
    \cref{proposition:I-SEG-almost-sure-strictly-monotone} and
    \citep{kannanOptimalStochasticExtragradient2019,mertikopoulosOptimisticMirrorDescent2019}.
    We prove the corresponding result for S-SEG: if $F$ is quasi-strictly monotone,
    $\sum_k\gamma_k=\infty$, and $\gamma_k, \extstep$ are both square-summable, then the last iterate converges almost surely to a solution
    (\cref{theorem:S-SEG-almost-sure-strict-monotone}).

    Without quasi-strict monotonicity, the two variants behave fundamentally differently.
    For I-SEG on $\cX=\reals^d$, the asymmetric step-size rule~\eqref{eqn:DSEG}
    is known to yield almost-sure last-iterate convergence \citep{HsiehIutzelerMalickMertikopoulos2020_explore}; we additionally show that it
    gives an expected restricted-gap bound of order $\cO((\sum_{k=0}^N\gamma_k)^{-1})$ for the weighted average of the iterates
    (\cref{proposition:independent-sample-DSEG-restricted-gap}).
    In contrast, the same step-size principle does not rescue S-SEG.
    We construct a uniformly Lipschitz stochastic monotone VIP with a unique solution, for which both the last iterate and the weighted average of the iterates diverge almost surely
    under DSEG step-sizes
    (\cref{theorem:same-sample-DSEG-gamma-average-divergence}).
    This separation shows that independence of the two samples $\xi_{i_k}$ and $\xi_{\hat{i}_k}$ is a central ingredient for the convergence mechanism of DSEG, rather than a proof convenience.

    \item [$\diamond$] \textbf{Bounded variance is insufficient for \ref{eqn:S-SEG}, even on compact domains.}

    On compact domains, the basic gap function convergence picture is largely known:
    with an appropriate choice of step-size, I-SEG and S-SEG achieve the standard
    $\cO(N^{-1/2})$ expected gap convergence rate under their respective assumptions; see \cref{proposition:I-SEG-gap-compact}, \cite[Proposition~4]{yousefianStochasticMirrorproxAlgorithms2018} and
    \cref{proposition:S-SEG-gap-compact}.
    Our contribution in this regime is to show that the stronger uniform Lipschitzness used in the S-SEG analysis is substantive.
    We construct a counterexample of monotone and Lipschitz stochastic VIP defined on a compact domain, whose stochastic oracle is unbiased and has bounded variance, and whose solution $x_\star$ is unique, but whose sample operators are not Lipschitz continuous.
    For a positive-measure set of initial points, the S-SEG iterates remain separated from $x_\star$, and the averaged iterates also fail to converge to $x_\star$
    (\cref{proposition:same-sample-nonlipschitz-realizations-positive-measure}).
    This shows that bounded variance and Lipschitzness of $F$ cannot replace
    uniform Lipschitzness in the S-SEG theory, so \cref{assumption:uniform-lipschitzness} used in all of our S-SEG theorems is not merely an artifact of the proofs.
\end{itemize}

\section{Background and preliminaries}

In this section, we establish a list of assumptions required for stating different convergence theories of SEG, 
and explore closely related prior results on SEG.

\subsection{Theoretical assumptions and their hierarchy}
\label{section:assumptions}

An operator $G\colon \mathrm{dom}(G) \to \reals^d$ is \textit{monotone} if its domain $\mathrm{dom}(G)$ is convex and
\begin{align}
\nonumber
    \inprod{G(x) - G(y)}{x - y} \ge 0 , \quad \forall x, y \in \mathrm{dom}(G) .
\end{align}
For $L>0$, we say $G$ is $L$-\textit{Lipschitz} if 
\begin{align}
\nonumber
    \norm{G(x) - G(y)} \le L\norm{x - y} , \quad \forall x, y \in \mathrm{dom}(G) .
\end{align}
Throughout the paper, we globally assume that \eqref{eqn:VIP} has a nonempty solution set $\cX_\star \subseteq \cX$. 
Furthermore, we commonly use the following assumptions on the problem class.

\begin{assumption}[Mean monotonicity]
\label{assumption:monotonicity}
The mean operator $F = \frac{1}{n}\sum_{i=1}^n F(\cdot;\xi_i)$ is monotone.
\end{assumption}

\begin{assumption}[Mean Lipschitzness]
\label{assumption:lipschitzness}
The mean operator $F = \frac{1}{n}\sum_{i=1}^n F(\cdot;\xi_i)$ is $L$-Lipschitz for some $L>0$. 
\end{assumption}

\begin{assumption}[Uniform Lipschitzness]
\label{assumption:uniform-lipschitzness}
In \eqref{eqn:objective-F}, $F(\cdot;\xi_i)\colon \cX \to \reals^d$ are \emph{uniformly $L$-Lipschitz} for some $L>0$, i.e., each $F(\cdot;\xi_i)$ is $L$-Lipschitz.
\end{assumption}

Clearly, \cref{assumption:uniform-lipschitzness} implies \cref{assumption:lipschitzness}; while S-SEG requires \cref{assumption:uniform-lipschitzness} for convergence (as we show in \cref{proposition:same-sample-nonlipschitz-realizations-positive-measure}),
I-SEG only requires the weaker \cref{assumption:lipschitzness} and does not seem to directly benefit from assuming \cref{assumption:uniform-lipschitzness}.

Next, we state the assumptions on the stochastic gradients. 
Define by $\emptyset = \cF_0 \subseteq \hat{\cF}_0 \subseteq \cF_1 \subseteq \hat{\cF}_1 \subseteq \dots$ the natural filtration generated by the algorithm;
more precisely, $\cF_k = \sigma\left( \xi_{i_0}, \xi_{\hat{i}_0}, \dots, \xi_{i_{k-1}}, \xi_{\hat{i}_{k-1}} \right)$, $\hat{\cF}_k = \sigma\left(\cF_k, \xi_{i_k}\right)$ respectively encompass all randomness involved in determining $x_k$ and $\hat{x}_k$.
For S-SEG, we have $\hat{\cF}_k = \cF_{k+1}$.

\begin{assumption}[Unbiasedness]
\label{assumption:unbiasedness}
At each iteration $k \ge 0$, the stochastic operator is an unbiased estimator of $F(x_k)$ conditioned on $\cF_k$, i.e., $\expec{F(x_k; \xi_{i_k})\,\middle|\,\cF_k}{} = F(x_k)$.
For I-SEG, we also assume $\expec{F(\halfitr; \halfxi)\,\middle|\,\hat{\cF}_k}{} = F(\halfitr)$.
\end{assumption}

\begin{assumption}[Bounded variance]
\label{assumption:bounded-variance}
There is $\sigma > 0$ such that $\expec{\sqnorm{F(x_k) - F(x_k;\xi_{i_k})}\,\middle|\,\cF_k}{} \le \sigma^2$ for $k\ge 0$. 
For I-SEG, we also assume $\expec{\sqnorm{F(\halfitr) - F(\halfitr; \halfxi)}\,\middle|\,\hat{\cF}_k}{} \le \sigma^2$.
\end{assumption}

While \cref{assumption:bounded-variance} is common in the stochastic optimization literature, in this work, we more often consider the following more general bound.

\begin{assumption}[Generalized variance bound]
\label{assumption:generalized-variance-bound}
There exist a solution $x_\star \in \cX_\star$ and $\sigma, A \ge 0$ such that
\begin{align}
\label{eqn:generalized-variance-bound}
\begin{aligned}
    \expec{\sqnorm{F(x_k) - F(x_k;\xi_{i_k})}\,\middle|\,\cF_k}{} & \le A\sqnorm{x_k - x_\star} + \sigma^2 \\
    \expec{\sqnorm{F(\halfitr) - F(\halfitr; \halfxi)}\,\middle|\,\hat{\cF}_k}{} & \le A\sqnorm{\halfitr - x_\star} + \sigma^2 
\end{aligned}
\end{align}
for $k=0,1,\dots$, where the second line is considered only for I-SEG.
\end{assumption}

Under \cref{assumption:lipschitzness}---which is a minimal and global assumption made throughout this paper---\cref{assumption:generalized-variance-bound} is equivalent to the condition $\condexp{\sqnorm{F(x; \xi_i)}}{x} \le B^2 \sqnorm{x - x_\star} + G^2$ (for all $x$ where new stochastic sampling occurs, where $B,G \ge 0$ are fixed), 
which dates back to \citep{blumApproximationMethodsWhich1954,gladyshevStochasticApproximation1965} and was also considered in more recent works including \cite{wangStochasticFirstorderMethods2016,HsiehIutzelerMalickMertikopoulos2020_explore, neuDealingUnboundedGradients2024,alacaogluWeakerVarianceAssumptions2025,alacaogluSolvingStochasticVariational2026}.
Notably, this type of variance bound is suitable for modeling unconstrained stochastic VIP with linear sample operators featuring random matrices, unlike the stricter \cref{assumption:bounded-variance}.
In fact, \cref{assumption:generalized-variance-bound} is implied by uniform Lipschitzness (\cref{assumption:uniform-lipschitzness}) and unbiasedness of stochastic operators (\cref{assumption:unbiasedness}) in the finite sum setting that we consider in this paper:

\begin{proposition}
\label{proposition:uniform-lipschitzness-implies-variance-bound}
Suppose that Assumptions~\ref{assumption:uniform-lipschitzness} and \ref{assumption:unbiasedness} hold.
Then, for any $x_\star \in \cX_\star$, \eqref{eqn:generalized-variance-bound} holds with $A=2L^2$ and 
\[
    \sigma^2 = 2\underset{i=1,\dots,n}{\max} \sqnorm{F(x_\star; \xi_i)} < \infty .
\]
Furthermore, we have $\sqnorm{F(x_k; \xi_{i_k})} \le 2L^2 \sqnorm{x_k - x_\star} + \sigma^2$ almost surely.
\end{proposition}

\begin{proof}
By unbiasedness and Lipschitzness of $F(\cdot; \xi_{i_k})$, we have
\begin{align*}
    \condexp{\sqnorm{F(x_k) - F(x_k;\xi_{i_k})}}{\cF_k} & = \condexp{\sqnorm{F(x_k; \xi_{i_k})}}{\cF_k} - \sqnorm{F(x_k)} \\
    & \le \condexp{2\sqnorm{F(x_k; \xi_{i_k}) - F(x_\star; \xi_{i_k})} + 2\sqnorm{F(x_\star; \xi_{i_k})}}{\cF_k} \\
    & \le 2L^2 \sqnorm{x_k - x_\star} + 2 \max_{i=1,\dots,n} \sqnorm{F(x_\star; \xi_i)} .
\end{align*}
For I-SEG, $\condexp{\sqnorm{F(\halfitr) - F(\halfitr;\halfxi)}}{\hat{\cF}_k}$ can be bounded similarly.
Finally, $\sqnorm{F(x_k;\xi_{i_k})}$ can be bounded as
\begin{align}
    \sqnorm{F(x_k; \xi_{i_k})} & \le 2\sqnorm{F(x_k; \xi_{i_k}) - F(x_\star; \xi_{i_k})} + 2\sqnorm{F(x_\star; \xi_{i_k})} \le 2L^2 \sqnorm{x_k - x_\star} + 2 \max_{i=1,\dots,n} \sqnorm{F(x_\star; \xi_i)} .
    \label{eqn:second-moment-bound}
\end{align}
\end{proof}

Note that in particular, when $\cX$ is compact, Assumptions~\ref{assumption:uniform-lipschitzness} and \ref{assumption:unbiasedness} imply \cref{assumption:bounded-variance}.
In theorems on S-SEG, we will not directly assume Assumptions~\ref{assumption:bounded-variance} or \ref{assumption:generalized-variance-bound}, because they will be tacitly implied through other assumptions.

Finally, we introduce the following notion, which will be used in the discussion of almost-sure iterate convergence.

\begin{assumption}[Quasi-strict monotonicity]
\label{assumption:quasi-strict-monotone}
For any $x_\star \in \cX_\star$ and $x \in \cX \setminus \cX_\star$, we have $\inprod{F(x)}{x-x_\star} > 0$.
\end{assumption}

This assumption, also referred to as the acute angle condition \citep{kannanOptimalStochasticExtragradient2019} or strict coherence property \citep{mertikopoulosOptimisticMirrorDescent2019}, 
was used for deriving the almost-sure convergence of I-SEG variants.
It states that any solution $x_\star$ is stable and attracts the negative operator directions, but excludes some monotone operators such as purely rotational ones.
On the other hand, \citet{HsiehIutzelerMalickMertikopoulos2020_explore} states that I-SEG may not converge for merely monotone stochastic VIP (without \cref{assumption:quasi-strict-monotone}) with $\extstep = \gamma_k$ 
satisfying the standard Robbins--Monro type conditions $\sum_{k=0}^\infty \gamma_k = \infty$ and $\sum_{k=0}^\infty \gamma_k^2 < \infty$.
Our \cref{theorem:unbounded-gap-high-probability-tightness} demonstrates that quasi-strict monotonicity is also essential for the last-iterate convergence of S-SEG with symmetric step-size $\extstep = \gamma_k$.

\subsection{Prior work and closely related work}
\label{section:prior-work}

In the following, we introduce important prior work of high relevance, broadly dividing them into results on gap function and almost-sure last-iterate convergence.
At the end of the section, we briefly explain the challenges of analyzing \ref{eqn:S-SEG} compared to the case of \ref{eqn:I-SEG}.

\subsubsection{Results on gap function convergence}

In this section, we review some existing results on the gap function convergence.
After each statement, we highlight its dependence on restrictive assumptions such as the compactness of the domain or bounded variance condition, or other limitations.
Unlike the results presented here, the main theorems of this paper are not subject to those limitations.

We start with the analyses of \ref{eqn:I-SEG}.

\begin{proposition}{\cite[Corollary~1]{juditskySolvingVariationalInequalities2011}}
\label{proposition:I-SEG-gap-compact}
Let $\cX$ be compact and let $D = \mathrm{diam}(\cX)$.
Suppose that Assumptions~\ref{assumption:monotonicity}, \ref{assumption:unbiasedness} and \ref{assumption:bounded-variance} hold, and that $F$ satisfies 
$\norm{F(x) - F(y)} \le L\norm{x - y} + M$ for all $x, y \in \cX$ (which is weaker than \cref{assumption:lipschitzness}).
Then I-SEG with constant step-size $\gamma_k = \extstep \equiv \min\left\{ \frac{1}{\sqrt{3}L} , D\sqrt{\frac{1}{7(N+1)(M^2 + 2\sigma^2)}} \right\}$ exhibits the rate
\begin{align*}
    \expec{\err(\overline{\hat{x}}_N)}{} = \expec{\max_{u\in \cX} \inprod{F(u)}{\overline{\hat{x}}_N - u}}{} \le \max\left\{ \frac{7LD^2}{8(N+1)} , 7D \sqrt{\frac{M^2 + 2\sigma^2}{6(N+1)}} \right\} = \cO\left( \frac{LD^2}{N} + \frac{D(M+\sigma)}{\sqrt{N}}\right)
\end{align*}
where $\overline{\hat{x}}_N = \frac{1}{N+1} \sum_{k=0}^N \halfitr$.
\end{proposition}

\paragraph{Remark.}
This is a classical result for \textit{compact domains, assuming uniformly bounded variance.}
It shows that the dependence of I-SEG on $L$ is only in the non-dominant term, while the dominant term $\cO\left(\frac{D(M+\sigma)}{\sqrt{N}}\right)$ depends on $\sigma$ and $M$, the parameters quantifying stochasticity and nonsmoothness of $F$.
That is, the stochastic components in the I-SEG analysis are not sensitive to the Lipschitzness parameter of $F$ or $F(\cdot; \xi_i)$.
Similar $\cO\left(\frac{1}{\sqrt{N}}\right)$ rate on gap function at averaged iterates can also be proved for I-SEG with step-sizes decaying at the rate $\extstep = \gamma_k = \cO\left(\frac{1}{\sqrt{k}}\right)$; see, e.g., \cite[Proposition~4]{yousefianStochasticMirrorproxAlgorithms2018}.

\begin{proposition}{\cite[Theorem'~3]{GidelBerardVignoudVincentLacoste-Julien2019_variational}}
\label{proposition:I-SEG-gap-convex-closed}
Let $\cX$ be convex and closed.
Under Assumptions~\ref{assumption:monotonicity}, \ref{assumption:lipschitzness}, \ref{assumption:unbiasedness} and \ref{assumption:bounded-variance}, I-SEG with constant step-size $\gamma_k \equiv \gamma \le \frac{1}{2\sqrt{3}L}$ exhibits the rate
\begin{align*}
    \expec{\err_{B(y,D)}(\overline{\hat{x}}_N)}{} = \expec{\max_{u\in \cX \cap B(y,D)} \inprod{F(u)}{\overline{\hat{x}}_N - u}}{} \le \frac{D^2}{\gamma(N+1)} + \frac{7\gamma\sigma^2}{2}
\end{align*}
for any $D>0$, where $\overline{\hat{x}}_N = \frac{1}{N+1} \sum_{k=0}^N \halfitr$.
\end{proposition}

\paragraph{Remark.}
The above result shows that using the gap function, the $\cO\left(\frac{1}{\sqrt{N}}\right)$ convergence of I-SEG can be extended to possibly unbounded domain $\cX$, but \textit{it still requires uniformly bounded variance}.
Note that Lipschitzness of $F$ is assumed but uniform Lipschitzness of $F(\cdot; \xi_i)$ is not.
Furthermore, while not formally stated as a theorem, \cite[Theorem'~3]{GidelBerardVignoudVincentLacoste-Julien2019_variational} in fact proved the more general rate
\[
    \expec{\err_{B(y,D)}(\overline{\hat{x}}_N^\gamma)}{} \le \frac{D^2}{\sum_{k=0}^N \gamma_k} + \frac{7\sigma^2}{2\sum_{k=0}^N \gamma_k} \sum_{k=0}^N \gamma_k^2
\]
for $\overline{\hat{x}}_N^\gamma = \frac{1}{\sum_{k=0}^N \gamma_k} \sum_{k=0}^N \gamma_k \hat{x}_k$, which can be used to derive 
a similar $\cO\left(\frac{1}{\sqrt{N}}\right)$ rate on gap function with decreasing step-sizes $\extstep = \gamma_k = \cO\left(\frac{1}{\sqrt{k}}\right)$.

Next, we turn to \ref{eqn:S-SEG}.

\begin{proposition}
\label{proposition:S-SEG-gap-compact}
Let $\cX$ be compact, $x_0 \in \cX$ and let $D = \mathrm{diam}(\cX)$.
Under Assumptions \ref{assumption:monotonicity}, \ref{assumption:uniform-lipschitzness} and \ref{assumption:unbiasedness}, 
S-SEG with deterministic, nonincreasing step-sizes $\extstep = \gamma_k \in \left(0 , \frac{1}{\sqrt{2}L}\right]$ (i.e, $0 < \gamma_k \le \dots \le \gamma_0 \le \frac{1}{\sqrt{2}L}$) exhibits the rate 
\begin{align*}
    \expec{\err(\overline{\hat{x}}_N)}{} \le \frac{D^2}{2(N+1)\gamma_{N+1}} + (8L^2 D^2 + 2\sigma^2) \frac{1}{N+1}\sum_{k=0}^N \gamma_k + \frac{\sigma D}{\sqrt{N+1}}
\end{align*}
where $\overline{\hat{x}}_N = \frac{1}{N+1} \sum_{k=0}^N \halfitr$ and $\sigma^2 = 2L^2 D^2 + 2\max_{i=1,\dots,n} \sqnorm{F(x_\star, \xi_i)}$ is a uniform upper bound on operator variance, implied by \cref{proposition:uniform-lipschitzness-implies-variance-bound}.
\end{proposition}

\paragraph{Remark.}
\cref{proposition:S-SEG-gap-compact} provides a qualitatively similar conclusion as \citep[Theorem'~5]{GidelBerardVignoudVincentLacoste-Julien2019_variational}, i.e., convergence of S-SEG on \textit{compact domains}, but we present it in a slightly different form with a self-contained, distinct proof detailed in Section~\ref{section:S-SEG-compact-analysis}. 
One immediate difference observed in \cref{proposition:S-SEG-gap-compact} from the previous results on I-SEG is that it uses \cref{assumption:uniform-lipschitzness}, which implies both Assumptions~\ref{assumption:lipschitzness} and \ref{assumption:bounded-variance} used in \cref{proposition:I-SEG-gap-compact}, when $\cX$ is compact.
Hence, it appears that S-SEG requires a stronger assumption compared to I-SEG.
In fact, this is not merely an artifact of the analysis but is a substantive requirement (see \cref{proposition:same-sample-nonlipschitz-realizations-positive-measure}).
With concrete step-size selections, \cref{proposition:S-SEG-gap-compact} directly yields the following $\cO\left( \frac{LD^2 + \nicefrac{\sigma^2}{L} + \sigma D}{\sqrt{N}} \right)$ convergence results, where $L$ is the uniform Lipschitzness parameter, which also indicate the sensitivity of S-SEG to the geometry of sample operators.

\begin{corollary}
\label{corollary:S-SEG-compact-constant-step}
Under the conditions of \cref{proposition:S-SEG-gap-compact}, for $N\ge 1$, S-SEG using the constant step-size $\gamma_k \equiv \gamma = \frac{1}{L\sqrt{N+1}}$ for $k=0,\dots,N+1$ exhibits the rate
\begin{align*}
    \expec{\err(\overline{\hat{x}}_N)}{} \le \frac{1}{\sqrt{N+1}} \left(9LD^2 + \frac{2\sigma^2}{L} + \sigma D \right) = \cO\left( \frac{LD^2 + \nicefrac{\sigma^2}{L} + \sigma D}{\sqrt{N}} \right)
\end{align*}
\end{corollary}

\begin{corollary}
\label{corollary:S-SEG-compact-diminishing-step}
Under the conditions of \cref{proposition:S-SEG-gap-compact}, S-SEG with diminishing step-size $\gamma_k = \extstep = \frac{1}{L\sqrt{k+2}}$ exhibits the rate 
\begin{align*}
    \expec{\err(\overline{\hat{x}}_N)}{} \le \frac{1}{N+1} \left( \frac{LD^2}{2} \sqrt{N+3} + \left( 16LD^2 + \frac{4\sigma^2}{L} \right) \sqrt{N+2} + \sqrt{N+1} \sigma D \right) = \cO\left( \frac{LD^2 + \nicefrac{\sigma^2}{L} + \sigma D}{\sqrt{N}} \right) .
\end{align*}
\end{corollary}

Finally, we state a prior result on S-SEG with possibly unbounded domains.

\begin{proposition}{\cite[Theorem~3]{MishchenkoKovalevShulginRichtarikMalitsky2020_revisiting}}
\label{proposition:S-SEG-flipped-measure}
Let $\cX$ be convex and closed.
Under Assumptions \ref{assumption:monotonicity}, \ref{assumption:uniform-lipschitzness} and \ref{assumption:unbiasedness},
S-SEG with $\gamma_k = \extstep = \cO\left(\frac{1}{\sqrt{k}L}\right)$ satisfies
\begin{align*}
    \max_{u\in \cX \cap B(y,D)} \, \expec{\inprod{F(u)}{\overline{\hat{x}}_N - u}}{} \le \frac{1}{\sqrt{N}L} \left(\frac{L^2 D^2}{2} + \sigma^2 \right) 
\end{align*}
for any $D>0$, where $\overline{\hat{x}}_N = \frac{1}{N+1} \sum_{k=0}^N \halfitr$.
\end{proposition}

\paragraph{Interpretation.}
\cref{proposition:S-SEG-flipped-measure}, on a first look, seems to offer a cleaner analysis under a more general setting, compared to \cref{proposition:S-SEG-gap-compact} and its corollaries.
However, the caveat is that \textit{it bounds} $\err_{B(y,D)}(\expec{\overline{\hat{x}}_k}{})$, \textit{rather than the standard error measure} $\expec{\err_{B(y,D)}(\overline{\hat{x}}_k)}{}$.
In general,
\begin{align*}
    \err_{B(y,D)}(\expec{\overline{\hat{x}}_k}{}) = \max_{x\in \cX \cap B(y,D)} \, \expec{\inprod{F(x)}{\overline{\hat{x}}_k - x}}{} \le \expec{\max_{x\in \cX \cap B(y,D)} \inprod{F(x)}{\overline{\hat{x}}_k - x}}{} = \expec{\err_{B(y,D)}(\overline{\hat{x}}_k)}{}
\end{align*}
so a bound on $\err_{B(y,D)}(\expec{\overline{\hat{x}}_k}{})$ is weaker.
Interpreting this result, it states that $\expec{\overline{\hat{x}}_k}{}$ will be an approximate solution; however, to obtain this vector, one should run the stochastic algorithm multiple times and take the average of $\overline{\hat{x}}_k$ generated from each run.
This does not adequately represent the type of guarantee that we typically want, and it does not prevent each individual run from failing loudly, as the following example demonstrates.

\begin{proposition}{\cite[Example~1.1]{alacaogluComplexitySimplePrimaldual2025}}
\label{proposition:wrong-measure-convergence}
Let $F\colon \reals^2 \to \reals^2$ be the monotone operator defined as $F(x_1, x_2) = (x_2, -x_1)$.
For this $F$, the unique solution to \eqref{eqn:VIP} is $x = (0, 0)$, but an algorithm that outputs
\begin{align*}
    x_k = \begin{cases}
        (k, k) & \text{ with probability } \frac{1}{2} \\
        (-k, -k) & \text{ with probability } \frac{1}{2}   
    \end{cases}
\end{align*}
at each iteration $k=1,2,\dots$ satisfies $\err_\cZ \left(\expec{x_k}{} \right) = 0$ for all $k \ge 1$,
for any compact set $\cZ \subset \reals^2$.
\end{proposition}

\paragraph{Interpretation.}
Here, we see that the algorithm does not use any information of the operator $F$ and 
outputs a random sequence whose norm diverges to infinity,
but it appears to \textit{``achieve a perfect convergence''} if $\err_\cZ(\expec{x_k}{})$ is taken as a measure of convergence.
This indicates the need to provide a guarantee in terms of $\expec{\err_{B(y,D)}(\overline{\hat{x}}_k)}{}$, rather than the one used in \cref{proposition:S-SEG-flipped-measure}.

\subsubsection{Existing results on almost-sure last-iterate convergence}

While it is common to use same extrapolation and update step sizes $\extstep = \gamma_k$ for SEG, it is well known (for I-SEG) that such step-size choice fails to provide last iterate convergence for general monotone problems \citep{HsiehIutzelerMalickMertikopoulos2020_explore}, but requires an additional condition such as \cref{assumption:quasi-strict-monotone}.

\begin{proposition}[\cite{mertikopoulosOptimisticMirrorDescent2019}, Theorem~4.3]
\label{proposition:I-SEG-almost-sure-strictly-monotone}
Let $\cX$ be compact.
Under Assumptions~\ref{assumption:monotonicity}, \ref{assumption:lipschitzness}, \ref{assumption:unbiasedness}, \ref{assumption:bounded-variance} and \ref{assumption:quasi-strict-monotone},
I-SEG with $\extstep = \gamma_k > 0$, satisfying $\sum_{k=0}^\infty \gamma_k = \infty$ and $\sum_{k=0}^\infty \gamma_k^2 < \infty$, converges almost surely to a solution to \eqref{eqn:VIP}\footnote[1]{The original result was stated in a more general form, but we simplify it for coherent exposition.}.
\end{proposition}

On the other hand, with step-sizes $\extstep, \gamma_k$ decaying at different rates, I-SEG can converge almost surely to a solution for any monotone Lipschitz stochastic VIP, without the need for the stronger assumption of quasi-strict monotonicity.

\begin{proposition}{\cite[Theorem~1]{HsiehIutzelerMalickMertikopoulos2020_explore}}
\label{proposition:I-SEG-almost-sure-DSEG}
Let $\cX=\reals^d$.
Under Assumptions~\ref{assumption:monotonicity},
\ref{assumption:lipschitzness},
\ref{assumption:unbiasedness}, and
\ref{assumption:generalized-variance-bound},
I-SEG with \emph{double step-size extragradient (DSEG)}
step-sizes satisfying $0<\gamma_k\le\hat{\gamma}_k <
\frac{1}{3\max\{L,\sqrt{A}\}}$ for all $k=0,1,\dots$, and
\begin{align}
\label{eqn:DSEG}
\tag{DSEG}
    \sum_{k=0}^\infty \gamma_k\hat{\gamma}_k=\infty,
    \qquad
    \sum_{k=0}^\infty \gamma_k^2<\infty,
    \qquad
    \sum_{k=0}^\infty \hat{\gamma}_k^2\gamma_k<\infty
\end{align}
converges almost surely to some $x_\infty\in\cX_\star$, where $F(x_\infty) = 0$ since $\cX=\reals^d$.
\end{proposition}

Interestingly, in \cref{subsection:DSEG-fails-to-fix-S-SEG}, we show that even the DSEG step-sizes cannot resolve the non-convergence of S-SEG on certain monotone and Lipschitz stochastic VIP (see \cref{theorem:same-sample-DSEG-gamma-average-divergence}).

\subsubsection{A challenge of analyzing S-SEG}

Overall, we observe that the literature has pursued the study of I-SEG more extensively, and this is likely because analyzing S-SEG is trickier due to its bias term
\begin{align}
\label{eqn:SEG-bias-term-for-half-iterate}
    \inprod{F(\halfitr; \halfxi) - F(\halfitr)}{x_k} .
\end{align}
This quantity becomes zero under total expectation for independent-sample SEG because
\begin{align*}
    \expec{\inprod{F(\halfitr; \halfxi) - F(\halfitr)}{x_k}}{} = \expec{\expec{\inprod{F(\halfitr; \halfxi) - F(\halfitr)}{x_k} \,\middle|\, \hat{\cF}_k}{}}{} = \expec{\inprod{\expec{F(\halfitr; \halfxi) - F(\halfitr) \,\middle|\, \hat{\cF}_k}{}}{x_k}}{} = 0 .
\end{align*}
However, for same-sample SEG where $\halfitr$ relies on $\xi_{i_k}$, the stochastic operator $F(\halfitr;\xi_{i_k})$ is no longer an unbiased estimator of $F(\halfitr)$ conditioned on $\halfitr$, so \eqref{eqn:SEG-bias-term-for-half-iterate} does not vanish.
This is not merely the limitation in the analyses,
but can be viewed as the fundamental limitation preventing some positive results established for I-SEG from carrying over to S-SEG (e.g., see \cref{theorem:same-sample-DSEG-gamma-average-divergence}).

\section{Gap function analysis: New results with minimal assumptions}

In this section, we provide several missing components within the existing gap function analyses of SEG in the literature.
We first consider the compact-domain analysis, and present the proof of \cref{proposition:S-SEG-gap-compact}, a convergence result for S-SEG.
We then show that the uniform Lipschitzness requirement (\cref{assumption:uniform-lipschitzness}) cannot be removed without replacement for S-SEG, by demonstrating an example of monotone VIP violating this assumption for which S-SEG fails to converge.

Next, we explore the extension of existing convergence results for I-SEG and S-SEG to general (possibly unbounded) domain $\cX$.
For I-SEG, convergence has been already established in \cite{GidelBerardVignoudVincentLacoste-Julien2019_variational} for uniformly bounded variance (\cref{assumption:bounded-variance}), but not under the more general \cref{assumption:generalized-variance-bound}.
We show in that case, we can guarantee a quantitative gap function convergence with high-probability; then we complement this with the impossibility results for some fundamentally stronger or more general guarantees.
On the other hand, we show that with the \ref{eqn:DSEG} step-size of \cite{HsiehIutzelerMalickMertikopoulos2020_explore}, I-SEG converges with respect to the gap function.

For S-SEG, we analyze it under the minimal \cref{assumption:uniform-lipschitzness}, which is stronger than but is of comparable generality with \cref{assumption:generalized-variance-bound}, 
and prove a similar high-probability bound on the gap function with symmetric (i.e., $\extstep = \gamma_k$) and square-summable step-sizes, also complemented by negative results showing that this cannot be fundamentally improved.
Finally, unlike in the case of I-SEG, \ref{eqn:DSEG} step-size does not benefit S-SEG (see \cref{section:almost-sure-convergence} for details).

\subsection{Analysis in compact domain}
\label{section:gap-function-compact-domain}

This section covers gap function convergence results in the case where the domain $\cX$ is compact.
We start with a self-contained proof of \cref{proposition:S-SEG-gap-compact}.

\subsubsection{Warm-up: Convergence of S-SEG}
\label{section:S-SEG-compact-analysis}

We first present a lemma that will be useful for subsequent proofs throughout the rest of the paper.

\begin{lemma}
\label{lemma:first-key-bound}
The $k$-th step of SEG satisfies, for any $u\in \cX$,
\begin{align*}
    0 & \le \frac{1}{2} \left(\sqnorm{u - x_k} - \sqnorm{u - x_{k+1}} \right) - \frac{1}{2}\sqnorm{x_{k+1} - \halfitr} - \frac{1}{2}\sqnorm{\halfitr - x_k} \\
    & \quad + \inprod{\extstep F(x_k; \xi_{i_k}) - \gamma_k F(\halfitr; \xi_{\hat{i}_k})}{x_{k+1} - \halfitr} - \gamma_k \inprod{F(\halfitr; \xi_{\hat{i}_k})}{\halfitr - u} .
\end{align*}
\end{lemma}

\begin{proof}
When $y \in \reals^d$, $z = \proj_\cX (y)$ and $u \in \cX$, we have
\begin{align}
\label{eqn:projection-inequality-general}
    \inprod{u - z}{z - y} \ge 0 .
\end{align}
Applying \eqref{eqn:projection-inequality-general} with $y = x_k - \extstep F(x_k;\xi_{i_k})$, $z = \halfitr$ and $u = x_{k+1}$ gives
\begin{align}
    0 & \le \inprod{x_{k+1} - \halfitr}{\halfitr - x_k + \extstep F(x_k; \xi_{i_k})} \nonumber \\
    & = \inprod{x_{k+1} - \halfitr}{\halfitr - x_k} + \inprod{x_{k+1} - \halfitr}{\extstep F(x_k; \xi_{i_k})} \nonumber \\
    & = \frac{1}{2} \left( \sqnorm{x_{k+1} - x_k} - \sqnorm{x_{k+1} - \halfitr} - \sqnorm{\halfitr - x_k} \right) + \inprod{\extstep F(x_k; \xi_{i_k})}{x_{k+1} - \halfitr} \label{eqn:projection-inequality-extrapolation}
\end{align}
Similarly, taking $y = x_k - \gamma_k F(\halfitr; \xi_{\hat{i}_k})$, $z = x_{k+1}$ in \eqref{eqn:projection-inequality-general} gives
\begin{align}
    0 & \le \inprod{u - x_{k+1}}{x_{k+1} - x_k + \gamma_k F(\halfitr; \xi_{\hat{i}_k})} \nonumber \\
    & = \inprod{u - x_{k+1}}{x_{k+1} - x_k} + \inprod{\gamma_k F(\halfitr; \xi_{\hat{i}_k})}{u - x_{k+1}} \nonumber \\
    & = \frac{1}{2} \left( \sqnorm{u - x_k} - \sqnorm{u - x_{k+1}} - \sqnorm{x_{k+1} - x_k} \right) + \inprod{\gamma_k F(\halfitr; \xi_{\hat{i}_k})}{u - x_{k+1}} .
    \label{eqn:projection-inequality-update}
\end{align}
Adding \eqref{eqn:projection-inequality-extrapolation} and \eqref{eqn:projection-inequality-update} together, we obtain
\begin{align*}
    0 & \le \frac{1}{2} \left(\sqnorm{u - x_k} - \sqnorm{u - x_{k+1}} \right) - \frac{1}{2}\sqnorm{x_{k+1} - \halfitr} - \frac{1}{2}\sqnorm{\halfitr - x_k} \\
    & \quad + \inprod{\extstep F(x_k; \xi_{i_k})}{x_{k+1} - \halfitr} + \inprod{\gamma_k F(\halfitr; \xi_{\hat{i}_k})}{u - x_{k+1}}  \\
    & = \frac{1}{2} \left(\sqnorm{u - x_k} - \sqnorm{u - x_{k+1}} \right) - \frac{1}{2}\sqnorm{x_{k+1} - \halfitr} - \frac{1}{2}\sqnorm{\halfitr - x_k} \\
    & \quad + \inprod{\extstep F(x_k; \xi_{i_k}) - \gamma_k F(\halfitr; \xi_{\hat{i}_k})}{x_{k+1} - \halfitr} - \gamma_k \inprod{F(\halfitr; \xi_{\hat{i}_k})}{\halfitr - u} .
\end{align*}
\end{proof}

\begin{proof}[\textbf{Proof of \cref{proposition:S-SEG-gap-compact}}]
Adding $\inprod{\gamma_k F(\halfitr; \xi_{i_k})}{\halfitr - u}$ to the both sides of \cref{lemma:first-key-bound} and using $i_k = \hat{i}_k$ and $\extstep = \gamma_k$, we obtain
\begin{align*}
    \inprod{\gamma_k F(\halfitr; \xi_{i_k})}{\halfitr - u} & \le \frac{1}{2} \left(\sqnorm{u - x_k} - \sqnorm{u - x_{k+1}} \right) - \frac{1}{2}\sqnorm{x_{k+1} - \halfitr} - \frac{1}{2}\sqnorm{\halfitr - x_k} \\
    & \quad + \inprod{\extstep F(x_k; \xi_{i_k})}{x_{k+1} - \halfitr} + \inprod{\gamma_k F(\halfitr; \xi_{i_k})}{\halfitr - x_{k+1}} \\
    & = \frac{1}{2} \left(\sqnorm{u - x_k} - \sqnorm{u - x_{k+1}} \right) - \frac{1}{2}\sqnorm{x_{k+1} - \halfitr} - \frac{1}{2}\sqnorm{\halfitr - x_k} \\
    & \quad + \gamma_k \inprod{F(x_k; \xi_{i_k}) - F(\halfitr; \xi_{i_k})}{x_{k+1} - \halfitr} .
\end{align*}
Next, dividing throughout by $\gamma_k$ and adding $\inprod{F(\halfitr) - F(\halfitr;\xi_{i_k})}{\halfitr - u}$ to the both sides gives
\begin{align}
\label{eqn:one-iteration-key-bound-first}
\begin{aligned}
    \inprod{F(\halfitr)}{\halfitr - u} & \le \frac{1}{2\gamma_k} \left(\sqnorm{u - x_k} - \sqnorm{u - x_{k+1}} \right) - \frac{1}{2\gamma_k}\sqnorm{x_{k+1} - \halfitr} - \frac{1}{2\gamma_k}\sqnorm{\halfitr - x_k} \\
    & \quad + \inprod{F(x_k; \xi_{i_k}) - F(\halfitr; \xi_{i_k})}{x_{k+1} - \halfitr} + \inprod{F(\halfitr) - F(\halfitr;\xi_{i_k})}{\halfitr - u} .
\end{aligned}
\end{align}
We use Young's inequality and Lipschitz continuity of $F(\cdot; \xi_{i_k})$ to bound:
\begin{align*}
    \inprod{F(x_k; \xi_{i_k}) - F(\halfitr; \xi_{i_k})}{x_{k+1} - \halfitr} & \le \frac{\gamma_k}{2} \sqnorm{F(x_k; \xi_{i_k}) - F(\halfitr; \xi_{i_k})} + \frac{1}{2\gamma_k} \sqnorm{x_{k+1} - \halfitr} \\
    & \le \frac{\gamma_k L^2}{2} \sqnorm{x_k - \halfitr} + \frac{1}{2\gamma_k} \sqnorm{x_{k+1} - \halfitr} 
\end{align*}
and plug this into \eqref{eqn:one-iteration-key-bound-first} to obtain
\begin{align}
    \inprod{F(\halfitr)}{\halfitr - u} & \le \frac{1}{2\gamma_k} \left(\sqnorm{u - x_k} - \sqnorm{u - x_{k+1}} \right) - \frac{1}{2\gamma_k} (1-\gamma_k^2 L^2) \sqnorm{\halfitr - x_k} + \inprod{F(\halfitr) - F(\halfitr;\xi_{i_k})}{\halfitr - u} \nonumber \\
    & = \left( \frac{1}{2\gamma_k}\sqnorm{u - x_k} - \frac{1}{2\gamma_{k+1}}\sqnorm{u - x_{k+1}} \right) + \frac{1}{2}\left( \frac{1}{\gamma_{k+1}} - \frac{1}{\gamma_k} \right) \sqnorm{u - x_{k+1}} \nonumber \\
    & \quad - \frac{1}{2\gamma_k} (1-\gamma_k^2 L^2) \sqnorm{\halfitr - x_k} + \inprod{F(\halfitr) - F(\halfitr;\xi_{i_k})}{\halfitr - u} \nonumber \\
    & \begin{aligned}
        & \le \left( \frac{1}{2\gamma_k}\sqnorm{u - x_k} - \frac{1}{2\gamma_{k+1}}\sqnorm{u - x_{k+1}} \right) + \frac{1}{2}\left( \frac{1}{\gamma_{k+1}} - \frac{1}{\gamma_k} \right) D^2 \\
        & \quad - \frac{1}{2\gamma_k} (1-\gamma_k^2 L^2) \sqnorm{\halfitr - x_k} + \inprod{F(\halfitr) - F(\halfitr;\xi_{i_k})}{\halfitr - u} 
    \end{aligned}
    \label{eqn:one-iteration-key-bound-second}
\end{align}
where the last inequality holds because $\gamma_{k+1} \le \gamma_k$ and $\norm{u - x_{k+1}} \le D$.
Next, we bound the last term in \eqref{eqn:one-iteration-key-bound-second}:
\begin{align*}
    & \inprod{F(\halfitr) - F(\halfitr;\xi_{i_k})}{\halfitr - u} \\
    & = \inprod{F(\halfitr) - F(x_k)}{\halfitr - u} + \inprod{F(x_k) - F(x_k; \xi_{i_k})}{\halfitr - u} + \inprod{F(x_k; \xi_{i_k}) - F(\halfitr; \xi_{i_k})}{\halfitr - u} \\
    & \le \frac{1}{16\gamma_k L^2} \sqnorm{F(\halfitr) - F(x_k)} + 4\gamma_k L^2 \sqnorm{\halfitr - u} + \inprod{F(x_k) - F(x_k; \xi_{i_k})}{\halfitr - u} \\
    & \quad + \frac{1}{16\gamma_k L^2} \sqnorm{F(x_k; \xi_{i_k}) - F(\halfitr; \xi_{i_k})} + 4\gamma_k L^2 \sqnorm{\halfitr - u} \\
    & \le \frac{1}{8 \gamma_k} \sqnorm{\halfitr - x_k} + 8\gamma_k L^2 D^2 + \inprod{F(x_k) - F(x_k; \xi_{i_k})}{\halfitr - x_k} + \inprod{F(x_k) - F(x_k; \xi_{i_k})}{x_k - u} \\
    & \le \frac{1}{8 \gamma_k} \sqnorm{\halfitr - x_k} + 8\gamma_k L^2 D^2 + 2\gamma_k \sqnorm{F(x_k) - F(x_k;\xi_{i_k})} + \frac{1}{8\gamma_k} \sqnorm{\halfitr - x_k} + \inprod{F(x_k) - F(x_k; \xi_{i_k})}{x_k - u} \\
    & = \frac{1}{4\gamma_k} \sqnorm{\halfitr - x_k} + 8\gamma_k L^2 D^2 + 2\gamma_k \sqnorm{F(x_k) - F(x_k;\xi_{i_k})} + \inprod{F(x_k) - F(x_k; \xi_{i_k})}{x_k - u} .
\end{align*}
Plugging this back into \eqref{eqn:one-iteration-key-bound-second} and using monotonicity of $F$ gives
\begin{align*}
    \inprod{F(u)}{\halfitr - u} & \le \left( \frac{1}{2\gamma_k}\sqnorm{u - x_k} - \frac{1}{2\gamma_{k+1}}\sqnorm{u - x_{k+1}} \right) + \frac{1}{2}\left( \frac{1}{\gamma_{k+1}} - \frac{1}{\gamma_k} \right) D^2 \underbrace{- \frac{1}{2\gamma_k} \left( \frac{1}{2} - \gamma_k^2 L^2 \right) \sqnorm{\halfitr - x_k}}_{\le 0} \\
    & \quad + 8\gamma_k L^2 D^2 + 2\gamma_k \sqnorm{F(x_k) - F(x_k;\xi_{i_k})} + \inprod{F(x_k) - F(x_k; \xi_{i_k})}{x_k - u} .
\end{align*}
We now sum this up for $k=0,\dots,N$, take maximum over $u\in \cX$ to get
\begin{align}
\label{eqn:decreasing-step-inner-product-sum-bound}
    & \max_{u\in \cX} \sum_{k=0}^N \inprod{F(u)}{\halfitr - u} \nonumber \\
    & \begin{aligned}
        & \le \max_{u\in\cX} \left\{ \frac{1}{2\gamma_0} \sqnorm{u - x_0} - \frac{1}{2\gamma_{N+1}} \sqnorm{u - x_{N+1}} \right\} + \frac{1}{2} \left(\frac{1}{\gamma_{N+1}} - \frac{1}{\gamma_0}\right) D^2 + 8L^2 D^2 \sum_{k=0}^N \gamma_k \\
        & \quad + \sum_{k=0}^N 2\gamma_k \sqnorm{F(x_k) - F(x_k;\xi_{i_k})} + \sum_{k=0}^N \inprod{F(x_k) - F(x_k; \xi_{i_k})}{x_k} + \max_{u\in\cX} \sum_{k=0}^N \inprod{F(x_k) - F(x_k; \xi_{i_k})}{-u}    
    \end{aligned}
     \\
    & \le \frac{D^2}{2\gamma_{N+1}} + 8L^2 D^2 \sum_{k=0}^N \gamma_k + \sum_{k=0}^N 2\gamma_k \sqnorm{F(x_k) - F(x_k;\xi_{i_k})} + \sum_{k=0}^N \inprod{F(x_k) - F(x_k; \xi_{i_k})}{x_k - x_0} \nonumber \\
    & \quad + \max_{u\in\cX} \sum_{k=0}^N \inprod{F(x_k) - F(x_k; \xi_{i_k})}{x_0 - u} . \nonumber
\end{align}
Now take the total expectation---under which the term $\sum_{k=0}^N \inprod{F(x_k) - F(x_k; \xi_{i_k})}{x_k - x_0}$ vanishes---to obtain
\begin{align}
\label{eqn:deterministic-bound-final}
    & (N+1) \expec{\err(\overline{\hat{x}}_N)}{} \le \frac{D^2}{2} \expec{\frac{1}{\gamma_{N+1}}}{} + (8L^2 D^2 + 2\sigma^2) \expec{\sum_{k=0}^N \gamma_k}{} + \expec{\max_{u\in\cX} \sum_{k=0}^N \inprod{F(x_k) - F(x_k; \xi_{i_k})}{x_0 - u}}{} .
\end{align}
Finally, we apply the standard trick (\citep[Lemma~3.1]{NemirovskiJuditskyLanShapiro2009_robust}, \citep[Proposition~11]{bachUniversalAlgorithmVariational2019}) to bound the final term.
Denoting $\Delta_k = F(x_k) - F(x_k; \xi_{i_k})$, we have $\condexp{\Delta_k}{\cF_k} = 0$.
Therefore, $\expec{\inprod{\Delta_j}{\Delta_k}}{} = 0$ for any $j < k$, which implies
\begin{align*}
    \expec{\max_{u\in\cX} \sum_{k=0}^N \inprod{\Delta_k}{x_0 - u}}{} \le D\expec{\norm{\sum_{k=0}^N \Delta_k}}{} \le D \left( \expec{\sqnorm{\sum_{k=0}^N \Delta_k}}{} \right)^{1/2} = D \left( \sum_{k=0}^N \expec{\sqnorm{\Delta_k}}{} \right)^{1/2} \le D\sigma \sqrt{N+1} .
\end{align*}
Applying this bound for \eqref{eqn:deterministic-bound-final} and dividing throughout by $N+1$ proves the desired result.
\end{proof}

\subsubsection{Importance of uniform Lipschitzness for S-SEG}
\label{section:S-SEG-uniform-Lipschitzness-necessity}

We will now show that the uniform (samplewise) Lipschitzness assumption is a basic requirement for convergence of S-SEG.
Without such assumption, S-SEG may fail to display any meaningful convergence even if $\cX$ is compact, as the following result shows.

\begin{proposition}
\label{proposition:same-sample-nonlipschitz-realizations-positive-measure}
Let $\gamma_k , \extstep$ be any predetermined positive sequence satisfying $\gamma_k \le \extstep$ for all $k\ge 0$.
There exists a stochastic VIP~\eqref{eqn:VIP} with compact domain $\cX$ satisfying Assumptions~\ref{assumption:monotonicity}, \ref{assumption:lipschitzness}, \ref{assumption:unbiasedness} and \ref{assumption:bounded-variance} with unique solution $x_\star = 0$, 
but each stochastic realization $F(\cdot;\xi_i)$ ($i=1,2$) is not Lipschitz continuous, 
and there exists a closed subset of $\cY \subset \cX$ of positive measure such that $x_\star \notin \cY$ and for every initialization $x_0 \in \cY$, S-SEG iterates satisfy $x_k \in \cY$ for all $k=0,1,\dots$, and $\liminf_{N\to\infty} \overline{\hat{x}}_N \ge \frac{1}{2}$ almost surely where $\overline{\hat{x}}_N = \frac{1}{N+1} \sum_{k=0}^N \halfitr$.
\end{proposition}

\begin{proof}
Let $\cX = [0,a] \subset \reals$, where $a \ge 2$. 
\[
    F(x;\xi_1)=
    \begin{cases}
    2x & \text{if } x<1 \\
    -1 & \text{if } x\ge 1
    \end{cases} ,
    \qquad
    F(x;\xi_2)=
    \begin{cases}
    0 & \text{if } x<1 \\
    3 & \text{if } x\ge 1
    \end{cases} .
\]
At each iteration, $i_k$ will be sampled uniformly from $\{1,2\}$. Then the mean operator is
\[
F(x) = \frac{1}{2} \bigl(F(x;\xi_1)+F(x;\xi_2)\bigr) = 
\begin{cases}
    x & \text{if } x<1 \\
    1 & \text{if } x\ge 1 
    \end{cases}
= \min\{x,1\} .
\]
This is an increasing function on $\reals$, so $F$ is monotone, and additionally, $1$-Lipschitz.
However, both $F(\cdot;\xi_1)$ and $F(\cdot;\xi_2)$ have jump discontinuities, so they are not Lipschitz.
Clearly, Assumption~\ref{assumption:unbiasedness} holds if $i_k$ is sampled uniformly random in $\{1,2\}$.
Also, because $\cX = [0,a]$, we have $|F(x;\xi_i) - F(x)| \le a$ for all $x\in \cX$ and $i=1,2$, so \cref{assumption:bounded-variance} holds with $\sigma = a$.
We see that $x_\star=0$ is a solution to \eqref{eqn:VIP} because
\[
    \inprod{F(u)}{0-u} = -uF(u)\le 0
\]
for any $u \in \cX$.
Positive values of $x$ cannot be a solution and because taking $u=\frac{x}{2}$ we have $\inprod{F(u)}{x-u} = \frac{x}{2} \min\left\{ \frac{x}{2}, 1 \right\} > 0$.
This shows that $\cX_\star=\{0\}$.

Now define $\cY = [1,a]$. Consider S-SEG initialized at any $x_0 \in \cY$.
We show by induction that $x_k \in \cY$ for all $k\ge 0$.
Assuming $x_k \in \cY$ for some $k\ge 0$, if $i_k = 1$, then $F(x_k;\xi_1) = -1$, so we have
\begin{align}
\label{eqn:halfitr-at-least-1-when-ik-is-1}
    \halfitr = \proj_\cX (x_k + \extstep) \ge 1 \implies \halfitr \in \cY .
\end{align}
Hence we again have $F(\halfitr; \xi_1) = -1$, which implies $x_{k+1} = \proj_\cX(x_k + \gamma_k)\ge 1 \implies x_{k+1} \in \cY$.

If $i_k=2$, so $F(x_k;\xi_2)=3$ and we have $\halfitr = \proj_\cX (x_k - 3\extstep)$.
If $\halfitr \ge 1$, then we must have $x_k - 3\extstep \ge 1$.
Because in this case we have $\halfitr \in \cY \implies F(\halfitr; \xi_2) = 3$, we obtain
\[
    x_k - 3\gamma_k \ge x_k - 3\extstep \ge 1 \implies x_{k+1} = \proj_\cX(x_k - 3\gamma_k) \in \cY .
\]
Otherwise, $\halfitr < 1$ and we have $F(\halfitr; \xi_2) = 0$, which implies $x_{k+1} = \proj_\cX(x_k - 0) = x_k \in \cY$.
Thus, in either case, $x_{k+1} \in \cY$. 
This completes the induction and shows that $x_k \in \cY$ for all $k \ge 0$.

Now for the extrapolated iterates, because we always have $x_k \in \cY$, by \eqref{eqn:halfitr-at-least-1-when-ik-is-1} we have $\halfitr \ge 1$ if $i_k = 1$.
When $i_k = 2$, it may be the case that $\halfitr < 1$, but still we have $\halfitr \ge 0$ because $\halfitr \in \cX = [0, a]$.
This shows that $\overline{\hat{x}}_N \ge \frac{1}{N+1} \sum_{k=0}^N \mathbf{1}_{\{i_k = 1\}}$, and by the law of large numbers, $\frac{1}{N+1} \sum_{k=0}^N \mathbf{1}_{\{i_k = 1\}} \to \frac{1}{2}$ almost surely.
Thus $\liminf_{N\to\infty} \overline{\hat{x}}_N \ge \frac{1}{2}$ almost surely.
\end{proof}

\subsection{Analysis in general closed convex domain}
\label{section:gap-function-general-domain}

We now turn to the gap function convergence results in the case where the domain $\cX$ is closed and convex but possibly unbounded.
We first provide the high-probability convergence results for both I-SEG and S-SEG.

\subsubsection{High-probability gap function convergence under general variance bound}
\label{section:S-SEG-gap-unbounded-high-probability}

The proof of \cref{proposition:S-SEG-gap-compact} in Section~\ref{section:S-SEG-compact-analysis} relies on the fact that the iterates stay within a fixed radius, enforced via projection steps onto the compact domain $\cX$.
Other prior results for I-SEG often use a similar assumption that the second moment of the stochastic operator is bounded.
However, for more general problems with possibly unbounded domain $\cX$ and weaker variance condition in \cref{assumption:generalized-variance-bound},
we cannot provide a uniform, path-independent bound on $\sup_{k\ge 0} \sqnorm{x_k}$ or $\sup_{k\ge 0} \sqnorm{F(x_k; \xi_{i_k})}$.
Even if we use square-summable step-sizes $\gamma_k=\extstep$, for which $\sqnorm{x_k-x_\star}$ converges almost surely for $x_\star\in\cX_\star$ and hence ${x_k}$ is almost surely bounded, 
$\sup_{k\ge 0}\sqnorm{x_k - x_\star}$ may still be arbitrarily large along rare stochastic paths; see Section~\ref{section:impossibility-via-counterexamples} for details.

As a workaround, in this section, we provide a uniform quantitative bound on $\sqnorm{x_k - x_\star}$ that holds with high probability, with step-sizes satisfying $\gamma_k = \extstep$ and $\sum_{k=0}^\infty \gamma_k^2 < \infty$.
This allows us to prove a quantitative bound on the expected gap function $\err_{B(y,D)} (x)$ conditioned on that high-probability event, as follows.

\begin{theorem}[I-SEG, conditional expected gap bound]
\label{theorem:independent-sample-SEG-unbounded-domain-high-probability-bound}
Let $\cX$ be a closed, convex, not necessarily bounded domain.
Suppose Assumptions~\ref{assumption:monotonicity}, \ref{assumption:lipschitzness}, \ref{assumption:unbiasedness} and \ref{assumption:generalized-variance-bound} hold.
Let $x_\star\in\cX_\star$ and let $\sigma_\star^2 = \sigma^2+2\sqnorm{F(x_\star)}$.
Then, for any $\delta \in (0,1)$, there is an event $\cE_\delta$ with $\prob{\cE_\delta}\ge 1-\delta$ such that I-SEG with nonincreasing predetermined step-size $\extstep=\gamma_k \in \left(0,\frac{1}{\sqrt{A + 2L^2}}\right]$ with $\Gamma:=\sum_{k=0}^\infty \gamma_k^2<\infty$ satisfies
\[
    \sup_{k\ge 0}\sqnorm{x_k-x_\star}
    \le
    \frac{D_\delta^2}{4} := \frac{\left(\sqnorm{x_0-x_\star}+6\sigma_\star^2\Gamma\right)e^{6(A + 2L^2)\Gamma}}{\delta}.
\]
On this event $\cE_\delta$, there exists a constant $C>0$ such that $\overline{\hat{x}}_N=\frac{1}{N+1}\sum_{k=0}^N \halfitr$ 
satisfies
\[
    \condexp{
        \err_{x_\star,D_\delta/2}(\overline{\hat{x}}_N)
    }{\cE_\delta}
    \le
    \frac{D_\delta^2}{2(N+1)\gamma_{N+1}} +
    \frac{1}{1-\delta} \left[
        \frac{8\left( (A + 2L^2) D_\delta^2 + \sigma_\star^2 \right)}{N+1} \sum_{k=0}^N \gamma_k + 2D_\delta \sqrt{\frac{ (A + 2L^2) D_\delta^2 + \sigma_\star^2}{N+1}} 
    \right] .
\]
In particular, if $\gamma_k=\Theta((k+1)^{-\alpha})$ with $\frac12<\alpha<1$, then $\expec{\err_{x_\star,D_\delta/2}(\overline{\hat{x}}_N)\,\middle|\,\cE_\delta}{} = \cO\left(N^{-(1-\alpha)}\right)$.
\end{theorem}

\begin{theorem}[S-SEG, conditional expected gap bound]
\label{theorem:same-sample-SEG-unbounded-domain-high-probability-bound}
Let $\cX$ be a closed, convex, not necessarily bounded domain.
Suppose Assumptions \ref{assumption:monotonicity}, \ref{assumption:uniform-lipschitzness} and \ref{assumption:unbiasedness} hold.
Let $x_\star \in \cX_\star$.
Then, for any $\delta \in (0,1)$, 
there is an event $\cE_\delta$ with $\prob{\cE_\delta}\ge 1-\delta$ such that S-SEG with nonincreasing predetermined step-size $\extstep = \gamma_k \in \left(0 , \frac{1}{9L}\right]$ with $\Gamma := \sum_{k=0}^\infty \gamma_k^2 < \infty$ satisfies
\begin{align*}
    \sup_{k\ge 0} \sqnorm{x_k - x_\star} \le \frac{D_\delta^2}{4} := \frac{\left( \sqnorm{x_0 - x_\star} + 4\sigma^2 \Gamma \right) e^{24L^2 \Gamma}}{\delta}
\end{align*}
where $\sigma^2 = 2 \max_{i=1,\dots,n} \sqnorm{F(x_\star; \xi_i)}$.
On this event $\cE_\delta$, $\overline{\hat{x}}_N=\frac{1}{N+1}\sum_{k=0}^N \halfitr$ satisfies
\begin{align*}
    \expec{\err_{x_\star, D_\delta/2} (\overline{\hat{x}}_N)\,\middle|\, \cE_\delta}{} 
    \le \frac{D_\delta^2}{2(N+1)\gamma_{N+1}} + 
    \frac{1}{1-\delta} \left[ \frac{17L^2 D_\delta^2 + 2\sigma^2}{N+1} \sum_{k=0}^N \gamma_k + D_\delta \sqrt{ \frac{ \nicefrac{L^2 \delta D_\delta^2}{2} + 2\sigma^2 }{N+1} }  \right]
\end{align*}
In particular, if $\gamma_k=\Theta((k+1)^{-\alpha})$ with $\frac12<\alpha<1$, then $\expec{\err_{x_\star,D_\delta/2}(\overline{\hat{x}}_N)\,\middle|\,\cE_\delta}{} = \cO\left(N^{-(1-\alpha)}\right)$.
\end{theorem}

\begin{corollary}[SEG, high-probability convergence]
\label{corollary:conditional-expectation-to-high-probability}
Under the conditions of either \cref{theorem:independent-sample-SEG-unbounded-domain-high-probability-bound} or \cref{theorem:same-sample-SEG-unbounded-domain-high-probability-bound}, for any $\rho \in (0,1)$, the considered version of SEG respectively satisfies
\begin{align*}
    \prob{\err_{x_\star, D_\delta/2} (\overline{\hat{x}}_N) \le \frac{r_N(\delta)}{\rho}} \ge 1 - \delta - \rho ,
\end{align*}    
where $r_N(\delta)$ is the upper bound on $\condexp{\err_{x_\star, D_\delta/2} (\overline{\hat{x}}_N)}{\cE_\delta}$ given in each theorem.
\end{corollary}

\begin{proof}
Given the bound $\condexp{\err_{x_\star, D_\delta/2} (\overline{\hat{x}}_N)}{\cE_\delta} \le r_N(\delta)$ with $\prob{\cE_\delta} \ge 1-\delta$, the conditional Markov inequality implies
\begin{align*}
    \prob{\err_{x_\star, D_\delta/2} (\overline{\hat{x}}_N) \le \frac{r_N(\delta)}{\rho}} & \ge \prob{ \left\{\err_{x_\star, D_\delta/2} (\overline{\hat{x}}_N) \le \frac{r_N(\delta)}{\rho} \right\} \cap \cE_\delta } \\
    & \ge \prob{\cE_\delta} (1-\rho) \ge (1-\delta) (1-\rho) \ge 1 - \delta - \rho .
\end{align*}
\end{proof}

To prove Theorems~\ref{theorem:independent-sample-SEG-unbounded-domain-high-probability-bound} and \ref{theorem:same-sample-SEG-unbounded-domain-high-probability-bound}, we first present a handy lemma on controlling the growth of stochastic sequences.

\begin{lemma}
\label{lemma:almost-sure-sequence-lemma}
Let $(\Omega,\cF,\mathbb{P},(\cF_k)_{k\ge 0})$ be a filtered probability space.
Let $a_k, b_k, c_k, d_k$ be nonnegative, integrable, $\cF_k$-measurable random variables such that 
$\condexp{a_{k+1}}{\cF_k} \le (1 + d_k) a_k - b_k + c_k$, and $\sum_{k=0}^\infty c_k < \infty$, $\sum_{k=0}^\infty d_k < \infty$ almost surely.
Then
\begin{enumerate}
    \item[(i)] $a_k$ converges almost surely and $\sum_{k=0}^\infty b_k < \infty$ almost surely.
    \item[(ii)] If $a_0$ is fixed and $c_k, d_k$ are deterministic sequences with $C = \sum_{k=0}^\infty c_k < \infty$ and $D = \sum_{k=0}^\infty d_k < \infty$, then $\expec{a_k}{} \le (a_0 + C) e^D$ for all $k\ge 0$ and $\sup_k a_k < R_\delta := \frac{(a_0 + C) e^D}{\delta}$ with probability at least $1-\delta$.
    \end{enumerate}
\end{lemma}

\begin{proof}
(i) is standard; see, e.g., \cite[Theorem~1]{robbinsConvergenceTheoremNon1971}.

For (ii), let $p_0 = 1$ and $p_k = \prod_{j=0}^{k-1} (1+d_j)$ for $k\ge 1$. Then $p_k \le \exp\left( \sum_{j=0}^{k-1} \log (1+d_j) \right) \le \exp\left( \sum_{j=0}^{k-1} d_j \right) \le e^D$.
Hence
\begin{align*}
    \expec{\frac{a_{k+1}}{p_{k+1}} \,\middle|\, \cF_k }{} \le \frac{a_k (1+d_k)}{p_{k+1}} + \frac{c_k}{p_{k+1}} = \frac{a_k}{p_k} + \frac{c_k}{p_{k+1}} ,
\end{align*}
which shows that $M_k = \frac{a_k}{p_k} + \sum_{j=k}^\infty \frac{c_j}{p_{j+1}}$ is a (nonnegative) supermartingale, satisfying $\expec{M_{k+1} \,\middle|\, \cF_k}{} \le M_k$.
In particular, $\expec{a_k}{} \le p_k \expec{M_k}{} \le p_k \expec{M_0}{} \le e^D (a_0 + C)$ where we use $M_0 \le a_0 + \sum_{k=0}^\infty c_k = a_0 + C$.
Also, by Ville's inequality, we have
\begin{align*}
    \mathrm{Prob}\left[ \sup_{k\ge 0} M_k \ge \lambda \right] \le \frac{\expec{M_0}{}}{\lambda} \le \frac{a_0 + C}{\lambda}
\end{align*}
for any $\lambda > 0$. 
Therefore, for any $\delta > 0$, we have $\sup_{k\ge 0} M_k < \frac{a_0 + C}{\delta}$ with probability $\ge 1-\delta$, and on this event,
\begin{align*}
    \sup_{k\ge 0} a_k \le \left( \sup_{k\ge 0} p_k \right) \left(\frac{a_0 + C}{\delta} \right) \le \frac{(a_0 + C) e^D}{\delta} = R_\delta .
\end{align*}
\end{proof}

To apply \cref{lemma:almost-sure-sequence-lemma}, we provide the following key bounds on $\sqnorm{x_k - x_\star}$ for each version of SEG.

\begin{lemma}
\label{lemma:I-SEG-dist-to-sol-growth-control}
Given the setting and notation of \cref{theorem:independent-sample-SEG-unbounded-domain-high-probability-bound}, I-SEG satisfies
\begin{align*}
    \condexp{\sqnorm{\halfitr-x_\star}}{\cF_k} & \le \left( 1+\gamma_k^2 (A+2L^2) \right) \sqnorm{x_k - x_\star} + \gamma_k^2 \sigma_\star^2 , \\
    \condexp{\sqnorm{x_{k+1} - x_\star}}{\cF_k} & \le \left( 1+6\gamma_k^2(A+2L^2) \right) \sqnorm{x_k - x_\star} + 6\gamma_k^2\sigma_\star^2 .
\end{align*}
\end{lemma}

\begin{proof}
Using the variance bound \eqref{eqn:generalized-variance-bound}, for $z=x_k$ or $z=\halfitr$ and the corresponding $\xi=\xi_{i_k}$ or $\xi=\xi_{\hat{i}_k}$,
\begin{align}
    \condexp{\sqnorm{F(z;\xi)}}{z} &=
    \sqnorm{F(z)} + \condexp{\sqnorm{F(z;\xi)-F(z)}}{z} \nonumber \\
    & \le
    2\left( \sqnorm{F(x_\star)} + \sqnorm{F(z) - F(x_\star)} \right) +
    A\sqnorm{z - x_\star} + \sigma^2 \nonumber \\
    & \le 2\left( \sqnorm{F(x_\star)} + L^2 \sqnorm{z - x_\star} \right) +
    A\sqnorm{z - x_\star} + \sigma^2 \nonumber \\
    & = \nu\sqnorm{z-x_\star}+\sigma_\star^2 
    \label{eqn:I-SEG-second-moment-growth}
\end{align}
where $\nu = A+2L^2$. 
Then, we have
\begin{align*}
    \sqnorm{\hat{x}_k - x_\star} & = \sqnorm{\proj_\cX \left( x_k - \gamma_k F(x_k; \xi_{i_k}) \right) - \proj_\cX \left( x_\star \right)} \\
    & \le \sqnorm{x_k - \gamma_k F(x_k; \xi_{i_k}) - x_\star} \\
    & = \sqnorm{x_k - x_\star} - 2\inprod{\gamma_k F(x_k; \xi_{i_k})}{x_k - x_\star} + \gamma_k^2 \sqnorm{F(x_k; \xi_{i_k})} .
\end{align*}
Taking expectation conditioned on $\cF_k$, the second term becomes $-2\gamma_k \inprod{F(x_k)}{x_k - x_\star} \le 0$, and the last term can be bounded using \eqref{eqn:I-SEG-second-moment-growth}.
As a result, we obtain
\begin{align}
\label{eqn:I-SEG-high-probability-lemma-xkhat-bound}
    \expec{\sqnorm{\halfitr-x_\star}\,\middle|\,\cF_k}{}
    \le
    \sqnorm{x_k-x_\star}
    +
    \gamma_k^2\left(\nu\sqnorm{x_k-x_\star}+\sigma_\star^2\right)
    =
    \left( 1 + \gamma_k^2 \nu \right) \sqnorm{x_k - x_\star} + \gamma_k^2 \sigma_\star^2 
\end{align}
where we used $\inprod{F(x_k)}{x_k-x_\star}\ge 0$.
Using \cref{lemma:first-key-bound} with $u=x_\star$, with $\extstep = \gamma_k$ we have
\begin{align}
\label{eqn:I-SEG-high-probability-lemma-xkp1-bound}
\begin{aligned}
    \sqnorm{x_{k+1} - x_\star} & \le \sqnorm{x_k - x_\star} - \sqnorm{x_{k+1} - \halfitr} - \sqnorm{\halfitr - x_k} \\
    & \quad + \underbrace{2 \gamma_k \inprod{F(x_k; \xi_{i_k}) - F(\halfitr; \halfxi)}{x_{k+1} - \halfitr}}_{(\mathrm{I})} - 2\gamma_k \inprod{F(\halfitr; \halfxi)}{\halfitr - x_\star} .
\end{aligned}
\end{align}
The last term, conditioned on $\hat{\cF}_k$, is nonpositive in expectation. 
Next, we can bound $(\mathrm{I})$ as
\begin{align*}
    (\mathrm{I}) \le \gamma_k^2 \sqnorm{F(x_k; \xi_{i_k}) - F(\halfitr; \halfxi)} + \sqnorm{x_{k+1} - \halfitr} 
    & \le 2\gamma_k^2 \pr{ \sqnorm{F(x_k; \xi_{i_k})} + \sqnorm{F(\halfitr; \halfxi)} } + \sqnorm{x_{k+1} - \halfitr} .
\end{align*}
Plugging this into \eqref{eqn:I-SEG-high-probability-lemma-xkp1-bound}, taking the conditional expectation on $\cF_k$ and using \eqref{eqn:I-SEG-second-moment-growth} and \eqref{eqn:I-SEG-high-probability-lemma-xkhat-bound} with the tower property, we obtain
\begin{align*}
    \condexp{\sqnorm{x_{k+1} - x_\star}}{\cF_k} & \le \sqnorm{x_k - x_\star} - \underbrace{\condexp{\sqnorm{\halfitr - x_k}}{\cF_k}}_{ \ge 0 } + 2\gamma_k^2 \condexp{ \sqnorm{F(x_k; \xi_{i_k})}}{\cF_k} + 2\gamma_k^2 \condexp{ \sqnorm{F(\halfitr; \halfxi)}}{\cF_k} \\
    & \le \sqnorm{x_k - x_\star} + 2\gamma_k^2 \nu\sqnorm{x_k - x_\star} + 2\gamma_k^2 \sigma_\star^2 + 2\gamma_k^2 \left( \nu \condexp{\sqnorm{\halfitr - x_\star}}{\cF_k} + \sigma_\star^2 \right) \\
    & \le (1 + 2\gamma_k^2 \nu) \sqnorm{x_k - x_\star} + 4 \gamma_k^2 \sigma_\star^2 + 2\gamma_k^2 \nu  \left( \left( 1 + \gamma_k^2 \nu \right) \sqnorm{x_k - x_\star} + \gamma_k^2 \sigma_\star^2  \right) \\
    & \le (1 + 6\gamma_k^2 \nu) \sqnorm{x_k - x_\star} + 6\gamma_k^2 \sigma_\star^2    
\end{align*}
where the last line uses $\gamma_k^2 \nu \le 1$, which follows by our step-size choice.

\end{proof}

\begin{lemma}
\label{lemma:S-SEG-dist-to-sol-growth-control}
Given the setting and notation of  \cref{theorem:same-sample-SEG-unbounded-domain-high-probability-bound}, 
S-SEG satisfies
\begin{align*}
    \condexp{\sqnorm{x_{k+1} - x_\star}}{\cF_k} & \le (1+d_k) \sqnorm{x_k - x_\star} - b_k + c_k ,
\end{align*}
with
\begin{align*}
    b_k & = \mathbb{E}\left[\left( \frac{1}{4} - 17\gamma_k^2 L^2 \right)\sqnorm{\halfitr - x_k} + 2\gamma_k \inprod{F(\halfitr)}{\halfitr - x_\star} \,\middle|\, \cF_k \right] \ge 0 , 
    \qquad 
    c_k = 4\gamma_k^2 \sigma^2 , 
    \qquad d_k = 24\gamma_k^2 L^2 . 
\end{align*}
\end{lemma}

\begin{proof}
Because $x_\star \in \cX_\star$ be a solution to \eqref{eqn:VIP}, it satisfies $\inprod{F(x)}{x-x_\star} \ge 0$ for all $x \in \cX$.
Now using \cref{lemma:first-key-bound} with $u=x_\star$, we have
\begin{align*}
    \sqnorm{x_{k+1} - x_\star} & \le \sqnorm{x_k - x_\star} - \sqnorm{x_{k+1} - \halfitr} - \sqnorm{\halfitr - x_k} + 2\gamma_k \inprod{F(x_k; \xi_{i_k}) - F(\halfitr; \xi_{i_k})}{x_{k+1} - \halfitr} \\
    & \quad - 2\gamma_k \inprod{F(\halfitr; \xi_{i_k})}{\halfitr - x_\star} \\
    & \le \sqnorm{x_k - x_\star} - \sqnorm{x_{k+1} - \halfitr} - \sqnorm{\halfitr - x_k} + \gamma_k^2 \sqnorm{F(x_k; \xi_{i_k}) - F(\halfitr; \xi_{i_k})} + \sqnorm{x_{k+1} - \halfitr} \\
    & \quad -\underbrace{2\gamma_k \inprod{F(\halfitr)}{\halfitr - x_\star}}_{=e_k \ge 0} - 2\gamma_k \inprod{F(\halfitr; \xi_{i_k}) - F(\halfitr)}{\halfitr - x_\star} \\
    & \le \sqnorm{x_k - x_\star} - \sqnorm{\halfitr - x_k} + \gamma_k^2 \sqnorm{F(x_k; \xi_{i_k}) - F(\halfitr; \xi_{i_k})} - 2\gamma_k \inprod{F(\halfitr; \xi_{i_k}) - F(x_k; \xi_{i_k})}{\halfitr - x_\star} \\
    & \quad - 2\gamma_k \inprod{F(x_k; \xi_{i_k}) - F(x_k)}{\halfitr - x_\star} - 2\gamma_k \inprod{F(x_k) - F(\halfitr)}{\halfitr - x_\star} - e_k \\
    & \le \sqnorm{x_k - x_\star} - (1 - \gamma_k^2 L^2) \sqnorm{\halfitr - x_k} + \frac{1}{4L^2} \sqnorm{F(\halfitr; \xi_{i_k}) - F(x_k; \xi_{i_k})} + 4\gamma_k^2 L^2 \sqnorm{\halfitr - x_\star} \\
    & \quad - 2\gamma_k \inprod{F(x_k; \xi_{i_k}) - F(x_k)}{\halfitr - x_k} - 2\gamma_k \inprod{F(x_k; \xi_{i_k}) - F(x_k)}{x_k - x_\star} \\
    & \quad + \frac{1}{4L^2} \sqnorm{F(x_k) - F(\halfitr)} + 4\gamma_k^2 L^2 \sqnorm{\halfitr - x_\star} - e_k \\
    & \le \sqnorm{x_k - x_\star} - (1 - \gamma_k^2 L^2) \sqnorm{\halfitr - x_k} + \frac{1}{4} \sqnorm{x_k - \halfitr} + 4\gamma_k^2 L^2 (2\sqnorm{x_k - x_\star} + 2\sqnorm{x_k - \halfitr}) \\
    & \quad + 4\gamma_k^2 \sqnorm{F(x_k; \xi_{i_k}) - F(x_k)} + \frac{1}{4} \sqnorm{x_k - \halfitr} - 2\gamma_k \inprod{F(x_k; \xi_{i_k}) - F(x_k)}{x_k - x_\star} \\
    & \quad + \frac{1}{4} \sqnorm{x_k - \halfitr} + 4\gamma_k^2 L^2 (2\sqnorm{x_k - x_\star} + 2\sqnorm{x_k - \halfitr}) - e_k \\
    & \le (1 + 16\gamma_k^2 L^2) \sqnorm{x_k - x_\star} - \left[\left( \frac{1}{4} - 17\gamma_k^2 L^2 \right)\sqnorm{\halfitr - x_k} + e_k \right] + 4\gamma_k^2 \sqnorm{F(x_k; \xi_{i_k}) - F(x_k)} \\
    & \quad - 2\gamma_k \inprod{F(x_k; \xi_{i_k}) - F(x_k)}{x_k - x_\star} .
\end{align*} 
The conclusion follows by taking the expectation conditioned on $\cF_k$ and applying \cref{proposition:uniform-lipschitzness-implies-variance-bound} to bound the third term.

\end{proof}

\begin{proof}[\textbf{Proof of \cref{theorem:independent-sample-SEG-unbounded-domain-high-probability-bound}}]

Let $\nu = A + 2L^2$.
Applying \cref{lemma:almost-sure-sequence-lemma}(ii) to the bound of \cref{lemma:I-SEG-dist-to-sol-growth-control} with $b_k = 0 , c_k = 6\gamma_k^2\sigma_\star^2$ and $d_k = 6\gamma_k^2(A+2L^2)$,
we obtain
\[
    \expec{\sqnorm{x_k-x_\star}}{}
    \le
    \left(\sqnorm{x_0-x_\star}+6\sigma_\star^2\Gamma\right)e^{6\nu\Gamma}
    =
    \frac{\delta D_\delta^2}{4}
\]
for all $k\ge 0$, and an event $\cE_\delta$ with $\prob{\cE_\delta}\ge 1-\delta$
such that
\[
    \sup_{k\ge 0}\sqnorm{x_k-x_\star}
    \le
    \frac{
    \left(\sqnorm{x_0-x_\star}+6\sigma_\star^2\Gamma\right)e^{6\nu\Gamma}
    }{\delta}
    =
    \frac{D_\delta^2}{4}.
\]
It remains to prove the gap function bound conditioned on $\cE_\delta$. 
Let $\cX_{D_\delta/2}:=\cX\cap B \left( x_\star, \frac{D_\delta}{2} \right)$, and denote 
$G_k = F(x_k;\xi_{i_k}), \widehat{G}_k := F(\halfitr;\halfxi)$ and $\Delta_k = F(\halfitr) - \widehat G_k$.
For any $u\in\cX_{D_\delta/2}$, using \cref{lemma:first-key-bound}, we obtain
\begin{align}
    \inprod{\widehat{G}_k}{\halfitr-u} & \le 
    \frac{1}{2\gamma_k} \left( \sqnorm{u-x_k} - \sqnorm{u-x_{k+1}} \right)
    - \frac{1}{2\gamma_k}\sqnorm{\halfitr-x_k}
    - \frac{1}{2\gamma_k} \sqnorm{x_{k+1} - \halfitr} 
    + \inprod{G_k - \widehat{G}_k}{x_{k+1} - \halfitr} \nonumber \\
    & \le 
    \frac{1}{2\gamma_k} \left( \sqnorm{u-x_k} - \sqnorm{u-x_{k+1}} \right)
    - \frac{1}{2\gamma_k}\sqnorm{\halfitr-x_k}
    - \frac{1}{2\gamma_k} \sqnorm{x_{k+1} - \halfitr} \nonumber \\
    & \quad
    + \frac{\gamma_k}{2}\sqnorm{G_k-\widehat{G}_k} + \frac{1}{2\gamma_k} \sqnorm{x_{k+1} - \halfitr} \nonumber \\
    & \le 
    \frac{1}{2\gamma_k} \left( \sqnorm{u-x_k} - \sqnorm{u-x_{k+1}} \right)
    - \frac{1}{2\gamma_k}\sqnorm{\halfitr-x_k}
    + \frac{\gamma_k}{2}\sqnorm{G_k-\widehat{G}_k} .
\label{eqn:I-SEG-high-probability-theorem-initial-bound}
\end{align}
Now we add $T_k = \inprod{\Delta_k}{\halfitr - u}$ to the both sides of \eqref{eqn:I-SEG-high-probability-theorem-initial-bound} and upper bound it as
\begin{align*}
    T_k & = \inprod{\Delta_k}{x_k - u} + \inprod{\Delta_k}{\halfitr - x_k} \\
    & \le \inprod{\Delta_k}{x_k - u} + \frac{\gamma_k}{2} \sqnorm{\Delta_k} + \frac{1}{2\gamma_k} \sqnorm{\halfitr - x_k}
\end{align*}
to obtain
\begin{align*}
    \inprod{F(\halfitr) }{\halfitr-u} \le \frac{1}{2\gamma_k} \left( \sqnorm{u-x_k} - \sqnorm{u-x_{k+1}} \right) 
    + \frac{\gamma_k}{2}\sqnorm{G_k-\widehat{G}_k}
    + \inprod{\Delta_k}{x_k - u}
    + \frac{\gamma_k}{2} \sqnorm{\Delta_k} .
\end{align*}
Lower-bounding the left hand side with $\inprod{F(u)}{\halfitr - u}$ by monotonicity of $F$, summing from $k=0$ to $N$,
applying the same telescoping as in \cref{proposition:S-SEG-gap-compact} 
and maximizing over $u\in\cX_{D_\delta/2}$, we obtain on $\cE_\delta$
\begin{align}
    \max_{u\in\cX_{D_\delta/2}} \sum_{k=0}^N \inprod{F(u)}{\halfitr - u}
    = (N+1)\err_{x_\star,D_\delta/2}(\overline{\hat{x}}_N)
    &\le
    \frac{D_\delta^2}{2\gamma_{N+1}} + R_N + S_N ,
    \label{eqn:I-SEG-gap-summed}
\end{align}
where
\begin{align*}
    R_N & := \frac12\sum_{k=0}^N \gamma_k \left(
        \sqnorm{G_k-\widehat{G}_k} + \sqnorm{\Delta_k}
    \right), \\
    S_N & := \sum_{k=0}^N\inprod{\Delta_k}{x_k-x_\star}
    + \max_{u\in\cX_{D_\delta/2}} \sum_{k=0}^N\inprod{\Delta_k}{x_\star-u} .
\end{align*}
To bound the conditional expectations of $R_N$ and $S_N$, define
\begin{align*}
    \tau := \inf\left\{
        k\ge 0
        \,\middle|\,
        \norm{x_k-x_\star}>\frac{D_\delta}{2}
    \right\} , \qquad
    \widetilde{R}_N = \frac12\sum_{k=0}^N \gamma_k \ones_{\{\tau > k\}} \left(
        \sqnorm{G_k-\widehat{G}_k} + \sqnorm{\Delta_k}
    \right) .
\end{align*}
Here $\tau$ is a stopping time such that $\mathbf{1}_{\{\tau > k\}}$ is $\cF_k$-measurable, $\sqnorm{x_k - x_\star} \le \frac{D_\delta^2}{4}$ on the event $\{\tau > k\}$, and $\cE_\delta \subset \{\tau > k\}$ for any $k\ge 0$.
Therefore, we can bound the conditional expectation of $R_N$ as 
\begin{align}
\label{eqn:I-SEG-high-probability-RN-bound}
\begin{aligned}
    \condexp{R_N}{\cE_\delta} = \condexp{\widetilde{R}_N}{\cE_\delta} \le \frac{1}{1-\delta} \expec{\widetilde{R}_N}{} 
    & \le \frac{1}{1-\delta} \sum_{k=0}^N \gamma_k \expec{\ones_{\{\tau > k\}} \left( \sqnorm{G_k-\widehat{G}_k} + \sqnorm{\Delta_k} \right) }{} \\
    & \le \frac{1}{1-\delta} \sum_{k=0}^N \gamma_k \expec{\ones_{\{\tau > k\}} \left( 2\sqnorm{G_k} + 2\sqnorm{\widehat{G}_k} + \sqnorm{\Delta_k} \right) }{} .
\end{aligned}
\end{align}
Now, using \eqref{eqn:I-SEG-second-moment-growth} and the fact that $\sqnorm{x_k - x_\star} \le \frac{D_\delta^2}{4}$ on the event $\{\tau > k\}$, we have
\begin{align*}
    \expec{\ones_{\{\tau > k\}} \sqnorm{G_k} }{} \le \expec{ \ones_{\{\tau > k\}} \condexp{\sqnorm{G_k}}{\cF_k} }{} 
    \le \expec{ \ones_{\{\tau > k\}} \left( \nu\sqnorm{x_k - x_\star} + \sigma_\star^2 \right) }{} \le \nu D_\delta^2 + \sigma_\star^2 .
\end{align*}
Proceeding in a similar way using \eqref{eqn:I-SEG-second-moment-growth} together with \eqref{eqn:I-SEG-high-probability-lemma-xkhat-bound}, proved in \cref{lemma:I-SEG-dist-to-sol-growth-control}, gives
\begin{align*}
    \expec{\ones_{\{\tau > k\}} \sqnorm{\widehat G_k} }{} 
    & \le \expec{ \ones_{\{\tau > k\}} \condexp{\nu\sqnorm{\halfitr - x_\star} + \sigma_\star^2}{\cF_k} }{} \\
    & \le \expec{ \ones_{\{\tau > k\}} \left( \nu \left( (1+\gamma_k^2 \nu) \sqnorm{x_k - x_\star} + \gamma_k^2 \sigma_\star^2 \right) + \sigma_\star^2 \right) }{} \le 2\left( \nu D_\delta^2 + \sigma_\star^2 \right) 
\end{align*}
where we use $\gamma_k^2 \nu \le 1$ in the last inequality.
Finally, using \eqref{eqn:I-SEG-high-probability-lemma-xkhat-bound} again,
\begin{align}
    \expec{\ones_{\{\tau > k\}} \sqnorm{\Delta_k} }{} \nonumber
    & \le \expec{\ones_{\{\tau > k\}} \condexp{\sqnorm{\Delta_k}}{\hat{\cF}_k} }{} \nonumber \\
    & \le \expec{\ones_{\{\tau > k\}} \left( A\sqnorm{\halfitr - x_\star} + \sigma^2 \right) }{} \nonumber \\
    & \le \expec{\ones_{\{\tau > k\}} \left( A \; \condexp{\sqnorm{\halfitr - x_\star} }{\cF_k} + \sigma^2 \right) }{} \nonumber \\
    & \le \expec{\ones_{\{\tau > k\}} \left( A \left( (1+\gamma_k^2 \nu) \sqnorm{x_k - x_\star} + \gamma_k^2 \sigma_\star^2 \right) + \sigma^2 \right) }{} \nonumber \\
    & \le 2\left( \nu D_\delta^2 + \sigma_\star^2 \right) 
    \label{eqn:I-SEG-Delta-norm-bound}
\end{align}
where in the last line, we use $A \le \nu$ and $\gamma_k^2 \nu \le 1$.
Plugging these bounds back into \eqref{eqn:I-SEG-high-probability-RN-bound}, we obtain 
\begin{align}
\label{eqn:I-SEG-high-probability-RN-bound-final}
    \condexp{R_N}{\cE_\delta} \le \frac{8\left( \nu D_\delta^2 + \sigma_\star^2 \right)}{1-\delta} \sum_{k=0}^N \gamma_k .
\end{align}
Finally, to bound $S_N$, define $\widetilde{\Delta}_k =\Delta_k\ones_{\{\tau>k\}}$.
Because $\ones_{\{\tau>k\}}$ is $\hat{\cF}_k$-measurable and
$\condexp{\Delta_k}{\hat{\cF}_k}=0$, we also have $\condexp{\widetilde{\Delta}_k}{\hat{\cF}_k}=0$.
Moreover, on $\cE_\delta$, $\Delta_k=\widetilde{\Delta}_k$ for every $k \ge 0$, so
\[
    \ones_{\cE_\delta}S_N
    \le
    \left(
        \sum_{k=0}^N
        \inprod{\widetilde{\Delta}_k}{x_k-x_\star}
    \right)_+
    +
    D_\delta
    \norm{
        \sum_{k=0}^N\widetilde{\Delta}_k
    }.
\]
Therefore, by Cauchy-Schwarz inequality and the martingale difference property of $\widetilde \Delta_k$,
\begin{align}
    \condexp{S_N}{\cE_\delta} & = \frac{1}{\prob{\cE_\delta}} \expec{\mathbf{1}_{\cE_\delta} S_N }{} \nonumber \\
    & \le \frac{1}{1-\delta} \left( \expec{\left( \sum_{k=0}^N \inprod{\widetilde \Delta_k}{x_k - x_\star} \right)^2}{} ^{1/2} + \frac{D_\delta}{2} \left( \expec{\sqnorm{\sum_{k=0}^N \widetilde \Delta_k}}{} \right)^{1/2} \right) \nonumber \\
    & \le \frac{1}{1-\delta} \left( \left(\sum_{k=0}^N \expec{ \inprod{\widetilde \Delta_k}{x_k - x_\star}^2 }{} \right)^{1/2} + \frac{D_\delta}{2} \left( \sum_{k=0}^N \expec{\sqnorm{ \widetilde \Delta_k}}{} \right)^{1/2} \right) .
    \label{eqn:I-SEG-SN-pre-bound}
\end{align}
We have the upper bound \eqref{eqn:I-SEG-Delta-norm-bound} on $\expec{\sqnorm{\widetilde \Delta_k}}{}$, and 
because $\sqnorm{x_k-x_\star}\le \frac{D_\delta^2}{4}$ on $\{\tau>k\}$, this also implies
\begin{align}
    \expec{
        \inprod{\widetilde{\Delta}_k}{x_k-x_\star}^2
    }{} \le
    \frac{D_\delta^2}{4}
    \expec{\sqnorm{\widetilde{\Delta}_k}}{} 
    \label{eqn:I-SEG-stopped-scalar-second-moment}
\end{align}
Plugging \eqref{eqn:I-SEG-Delta-norm-bound} and \eqref{eqn:I-SEG-stopped-scalar-second-moment} into
\eqref{eqn:I-SEG-SN-pre-bound}, we obtain
\begin{align}
    \condexp{S_N}{\cE_\delta}
    &\le
    \frac{2D_\delta}{1-\delta}
    \sqrt{
        (N+1)
        \left(\nu D_\delta^2 + \sigma_\star^2\right)
    }.
    \label{eqn:I-SEG-high-probability-SN-bound}
\end{align}
Taking conditional expectation in \eqref{eqn:I-SEG-gap-summed} and using
\eqref{eqn:I-SEG-high-probability-RN-bound-final} and \eqref{eqn:I-SEG-high-probability-SN-bound} gives the desired bound
\begin{align*}
    \condexp{
        \err_{x_\star,D_\delta/2}(\overline{\hat{x}}_N)
    }{\cE_\delta}
    &\le
    \frac{D_\delta^2}{2(N+1)\gamma_{N+1}} +
    \frac{1}{1-\delta} \left[
        8\left( \nu D_\delta^2 + \sigma_\star^2 \right) \sum_{k=0}^N \gamma_k + 2D_\delta \sqrt{\frac{\nu D_\delta^2 + \sigma_\star^2}{N+1}} 
    \right] .
\end{align*}

\end{proof}

\begin{proof}[\textbf{Proof of \cref{theorem:same-sample-SEG-unbounded-domain-high-probability-bound}}]
Applying \cref{lemma:almost-sure-sequence-lemma}(ii) to the bound of \cref{lemma:S-SEG-dist-to-sol-growth-control}, we have 
\begin{align*}
    \expec{\sqnorm{x_k - x_\star}}{} \le \left( \sqnorm{x_0 - x_\star} + 4\sigma^2 \Gamma \right) e^{24L^2 \Gamma} = \frac{\delta D_\delta^2}{4}
\end{align*}
for all $k\ge 0$, and with probability at least $1-\delta$, 
\begin{align*}
    \sup_{k\ge 0} \sqnorm{x_k - x_\star} < \frac{\left( \sqnorm{x_0 - x_\star} + 4\sigma^2 \Gamma \right) e^{24L^2 \Gamma}}{\delta} = \frac{D_\delta^2}{4} .
\end{align*}
On this event $\cE_\delta$, we proceed using the same argument as in the proof of \cref{proposition:S-SEG-gap-compact}; 
for any $u \in \cX_{D_\delta/2} := \cX \cap B\left(x_\star, \frac{D_\delta}{2}\right)$, we have 
\begin{align}
\begin{aligned}
    \inprod{F(\halfitr)}{\halfitr - u} & \le \left( \frac{1}{2\gamma_k}\sqnorm{u - x_k} - \frac{1}{2\gamma_{k+1}}\sqnorm{u - x_{k+1}} \right) + \frac{1}{2}\left( \frac{1}{\gamma_{k+1}} - \frac{1}{\gamma_k} \right) D_\delta^2 \\
    & \quad - \frac{1}{2\gamma_k} (1-\gamma_k^2 L^2) \sqnorm{\halfitr - x_k} + \inprod{F(\halfitr) - F(\halfitr;\xi_{i_k})}{\halfitr - u} 
\end{aligned}
\label{eqn:high-probability-bound-before-summation}
\end{align}
because on $\cE_\delta$, $\norm{x_k - u} \le \norm{x_k - x_\star} + \norm{x_\star - u} \le D_\delta$ for all $k$.
Then we bound the last inner product as
\begin{align*}
    & \inprod{F(\halfitr) - F(\halfitr;\xi_{i_k})}{\halfitr - u} \\
    & = \inprod{F(\halfitr) - F(x_k)}{\halfitr - u} + \inprod{F(x_k) - F(x_k; \xi_{i_k})}{\halfitr - u} + \inprod{F(x_k; \xi_{i_k}) - F(\halfitr; \xi_{i_k})}{\halfitr - u} \\
    & \le \frac{1}{16\gamma_k L^2} \sqnorm{F(\halfitr) - F(x_k)} + 4\gamma_k L^2 \sqnorm{\halfitr - u} + \inprod{F(x_k) - F(x_k; \xi_{i_k})}{\halfitr - u} \\
    & \quad + \frac{1}{16\gamma_k L^2} \sqnorm{F(x_k; \xi_{i_k}) - F(\halfitr; \xi_{i_k})} + 4\gamma_k L^2 \sqnorm{\halfitr - u} \\
    & \le \frac{1}{8 \gamma_k} \sqnorm{\halfitr - x_k} + 8\gamma_k L^2 \sqnorm{\halfitr - u} + \inprod{F(x_k) - F(x_k; \xi_{i_k})}{\halfitr - x_k} + \inprod{F(x_k) - F(x_k; \xi_{i_k})}{x_k - u} \\
    & \le \frac{1}{8 \gamma_k} \sqnorm{\halfitr - x_k} + 8\gamma_k L^2 \left( 2\sqnorm{x_k - u} + 2\sqnorm{\halfitr - x_k} \right) + 2\gamma_k \sqnorm{F(x_k) - F(x_k;\xi_{i_k})} + \frac{1}{8\gamma_k} \sqnorm{\halfitr - x_k} \\
    & \quad + \inprod{F(x_k) - F(x_k; \xi_{i_k})}{x_k - u} \\
    & = \left(\frac{1}{4\gamma_k} + 16 \gamma_k L^2 \right) \sqnorm{\halfitr - x_k} + 16\gamma_k L^2 D_\delta^2 + 2\gamma_k \sqnorm{F(x_k) - F(x_k;\xi_{i_k})} + \inprod{F(x_k) - F(x_k; \xi_{i_k})}{x_k - u} .
\end{align*}
Plugging this back into \eqref{eqn:high-probability-bound-before-summation}, using monotonicity of $F$ and telescoping, we obtain
\begin{align}
\label{eqn:S-SEG-high-probability-total-bound}
\begin{aligned}
    & \max_{u\in \cX_{D_\delta/2}} \sum_{k=0}^N \inprod{F(u)}{\halfitr - u} \\
    & \le \frac{D_\delta^2}{2\gamma_{N+1}} - \sum_{k=0}^N \frac{1}{2\gamma_k} \left( \frac{1}{2} - 33\gamma_k^2 L^2 \right) \sqnorm{\halfitr - x_k} + 
    \underbrace{ \sum_{k=0}^N \gamma_k \left( 16L^2 D_\delta^2 + 2\sqnorm{F(x_k) - F(x_k; \xi_{i_k})} \right) }_{=R_N} \\
    & \quad + \underbrace{\sum_{k=0}^N \inprod{F(x_k) - F(x_k; \xi_{i_k})}{x_k - x_\star} + \max_{u\in \cX_{D_\delta/2}} \sum_{k=0}^N \inprod{F(x_k) - F(x_k; \xi_{i_k})}{x_\star - u}}_{= S_N} .
\end{aligned}
\end{align}
Note that $\frac{1}{2} - 33\gamma_k^2 L^2 > 0$ by our choice of $\gamma_k$, so the second term is nonpositive.
Next, define 
\begin{align*}
    \tau = \inf \left\{ k\ge 0 \,\middle|\, \norm{x_k - x_\star} > \frac{D_\delta}{2} \right\} , \qquad \widetilde{R}_N = \sum_{k=0}^N \gamma_k \ones_{\{\tau > k\}} \left( 16L^2 D_\delta^2 + 2\sqnorm{F(x_k) - F(x_k; \xi_{i_k})} \right) .
\end{align*}
Here $\tau$ is a stopping time such that $\mathbf{1}_{\{\tau > k\}}$ is $\cF_k$-measurable, $\sqnorm{x_k - x_\star} \le \frac{D_\delta^2}{4}$ on the event $\{\tau > k\}$, and $\cE_\delta \subset \{\tau > k\}$ for any $k\ge 0$.
Therefore, we can bound the conditional expectation of $R_N$ as 
\begin{align}
\label{eqn:S-SEG-high-probability-RN-bound}
\begin{aligned}
    \condexp{R_N}{\cE_\delta} = \condexp{\widetilde{R}_N}{\cE_\delta} \le \frac{1}{1-\delta} \expec{\widetilde{R}_N}{} 
    & \le \frac{1}{1-\delta} \sum_{k=0}^N \gamma_k \left( 16L^2 D_\delta^2 + 2\left( 2L^2\expec{\mathbf{1}_{\{\tau > k\}} \sqnorm{x_k - x_\star}}{} + \sigma^2 \right) \right) \\
    & \le \frac{17L^2 D_\delta^2 + 2\sigma^2}{1-\delta} \sum_{k=0}^N \gamma_k .
\end{aligned}
\end{align}
Finally, the conditional expectation of $S_N$ can be bounded similarly as in the proof of \cref{theorem:independent-sample-SEG-unbounded-domain-high-probability-bound}:
letting
\begin{align*}
    \Delta_k = F(x_k) - F(x_k; \xi_{i_k}) , \qquad \widetilde \Delta_k = \Delta_k \mathbf{1}_{\{\tau > k\}} ,
\end{align*}
we have $\condexp{\widetilde \Delta_k}{\cF_k} = 0$ and $\Delta_k = \widetilde{\Delta}_k$ on $\cE_\delta$, which implies
\[
    \mathbf{1}_{\cE_\delta} S_N \le 
    \left( \sum_{k=0}^N \inprod{ \widetilde \Delta_k}{x_k - x_\star} \right)_{+} + \frac{D_\delta}{2} \norm{\sum_{k=0}^N \widetilde \Delta_k} .
\]
Using Cauchy-Schwarz and the martingale difference property of $\widetilde \Delta_k$, 
we obtain 
\begin{align*}
    \condexp{S_N}{\cE_\delta} \le \frac{1}{1-\delta} \left( \left(\sum_{k=0}^N \expec{ \inprod{\widetilde \Delta_k}{x_k - x_\star}^2 }{} \right)^{1/2} + \frac{D_\delta}{2} \left( \sum_{k=0}^N \expec{\sqnorm{ \widetilde \Delta_k}}{} \right)^{1/2} \right) .
\end{align*}
Now we note that by \cref{proposition:uniform-lipschitzness-implies-variance-bound} and $\expec{\sqnorm{x_k - x_\star}}{} \le \frac{\delta D_\delta^2}{4}$,
\begin{align*}
    \expec{\sqnorm{\widetilde \Delta_k}}{} & = \expec{ \ones_{\{\tau > k\}} \condexp{ \sqnorm{F(x_k) - F(x_k; \xi_{i_k})}}{\cF_k}}{}  \le \expec{ 2L^2\sqnorm{x_k - x_\star} + \sigma^2}{} \le \frac{L^2 \delta D_\delta^2}{2} + 2\sigma^2
\end{align*}
and furthermore, because $\widetilde{\Delta}_k$ is supported on $\{\tau > k\}$ where $\sqnorm{x_k - x_\star} \le \frac{D_\delta^2}{4}$,
\begin{align*}
    \expec{ \inprod{\widetilde \Delta_k}{x_k - x_\star}^2 }{} \le \expec{\mathbf{1}_{\{\tau > k\}} \sqnorm{\widetilde \Delta_k} \sqnorm{x_k - x_\star}}{} \le \frac{D_\delta^2}{4} \expec{\sqnorm{\widetilde \Delta_k}}{} \le \frac{D_\delta^2}{4} \left( \frac{L^2 \delta D_\delta^2}{2} + 2\sigma^2 \right) .
\end{align*}
Plugging the resulting bound on $\condexp{S_N}{\cE_\delta}$ and \eqref{eqn:S-SEG-high-probability-RN-bound} into \eqref{eqn:S-SEG-high-probability-total-bound},
we obtain 
\begin{align*}
    \condexp{\max_{u\in \cX_{D_\delta/2}} \sum_{k=0}^N \inprod{F(u)}{\halfitr - u}}{\cE_\delta} \le \frac{D_\delta^2}{2\gamma_{N+1}} + \frac{17L^2 D_\delta^2 + 2\sigma^2}{1-\delta} \sum_{k=0}^N \gamma_k + \frac{D_\delta}{1-\delta} \sqrt{(N+1) \left( \frac{L^2 \delta D_\delta^2}{2} + 2\sigma^2 \right)} .
\end{align*}
The conclusion follows by dividing throughout by $N+1$.
\end{proof}

\subsubsection{Beyond the high-probability theorems: Impossibility via counterexamples}
\label{section:impossibility-via-counterexamples}

Theorems~\ref{theorem:independent-sample-SEG-unbounded-domain-high-probability-bound} and \ref{theorem:same-sample-SEG-unbounded-domain-high-probability-bound} from the previous section provide high-probability guarantees on the gap function (\cref{corollary:conditional-expectation-to-high-probability}) given the step-size satisfying $\extstep = \gamma_k$ and $\sum_{k=0}^\infty \gamma_k^2 < \infty$.
In this section, we provide the results complementary to them. 
First we show that for both I-SEG and S-SEG, the square summability assumption on the step-size cannot be removed without replacement, by presenting a stochastic VIP for which even the high-probability gap function convergence cannot be guaranteed with non-square-summable step-sizes.
Then we show that even when $\gamma_k = \extstep$ and $\sum_{k=0}^\infty \gamma_k^2 < \infty$, there exists a counterexample where each SEG trajectory escapes any bounded disk with positive probability.
Therefore, no pre-specified bounded region makes the gap function a meaningful measure of convergence with probability $1$, showing that the high-probability approach of Theorems~\ref{theorem:independent-sample-SEG-unbounded-domain-high-probability-bound} and \ref{theorem:same-sample-SEG-unbounded-domain-high-probability-bound} is suitably tight.

\begin{theorem}
\label{theorem:same-sample-unbounded-divergence-nonzero-mean}
There exists a stochastic VIP~\eqref{eqn:VIP} on the unconstrained domain $\cX=\reals^2$ with unique solution $x_\star = 0$, satisfying Assumptions~\ref{assumption:monotonicity}, \ref{assumption:uniform-lipschitzness} and \ref{assumption:unbiasedness}, 
but S-SEG with $\extstep = \gamma_k = \frac{1}{\sqrt{k+1}}$ satisfies $\norm{x_k} \to \infty$ and $\limsup_{k\to\infty} \norm{\overline{\hat{x}}_k} = \infty$ almost surely for every $x_0 \neq x_\star$, where $\overline{\hat{x}}_k = \frac{1}{k+1}\sum_{j=0}^k \hat{x}_j$. 
Consequently, for every $D>0$,
\[
    \err_{x_\star, D}(x_k) \to \infty,
    \qquad
    \limsup_{k\to\infty} \err_{x_\star, D}(\overline{\hat{x}}_k) = \infty
\]
almost surely.
\end{theorem}

\begin{theorem}
\label{theorem:independent-sample-unbounded-divergence-generalized-variance}
There exists a stochastic VIP~\eqref{eqn:VIP} on $\cX=\reals^2$ with unique solution $x_\star = 0$, satisfying Assumptions~\ref{assumption:monotonicity}, \ref{assumption:lipschitzness}, \ref{assumption:unbiasedness} and \ref{assumption:generalized-variance-bound}, 
but I-SEG with $\extstep=\gamma_k=\frac{1}{\sqrt{k+1}}$ satisfies 
$\norm{x_k} \to \infty$ and $\limsup_{k\to\infty} \norm{\overline{\hat{x}}_k} = \infty$ almost surely for every $x_0 \neq x_\star$, where $\overline{\hat{x}}_k = \frac{1}{k+1}\sum_{j=0}^k \hat{x}_j$. 
Consequently, for every $D>0$,
\[
    \err_{x_\star,D}(x_k)\to\infty,
    \qquad
    \limsup_{k\to\infty}\err_{x_\star,D}(\overline{\hat{x}}_k)=\infty
\]
almost surely.
\end{theorem}

\paragraph{Remark.} While the counterexamples for the above theorems use the step-sizes $\gamma_k = \frac{1}{\sqrt{k+1}}$ for simplicity, they can be replaced by $\frac{1}{\sqrt{k+r}}$ with any $r>0$ without any essential change to the proof, so that $\gamma_k$ satisfies the requirements in convergence theorems presented in earlier sections.
\medskip

In the proofs of Theorems~\ref{theorem:same-sample-unbounded-divergence-nonzero-mean} and \ref{theorem:independent-sample-unbounded-divergence-generalized-variance} below and several later results, we will consider the specific linear counterexample sharing the same structure. 
We present a handy lemma for analyzing the iterates of SEG when run on general problems of such form.

\begin{lemma}
\label{lemma:stochastic-vip-counterexample}
Let $J = \begin{bmatrix} 0 & -1 \\ 1 & 0 \end{bmatrix}$, and let $F(x;\xi_1) = (\rho I + (1+\nu)J)x$ and $F(x;\xi_2) = (-\rho I + (1-\nu)J)x$ for $x \in \cX = \reals^2$, where $\rho, \nu \in \reals$ are fixed.
Then the stochastic VIP with $F(x) = \frac{1}{2}\left( F(x; \xi_1) + F(x; \xi_2) \right)$ satisfies \cref{assumption:monotonicity} and \cref{assumption:uniform-lipschitzness} with $L = \sqrt{\rho^2 + (1+|\nu|)^2}$, has the mean operator $F(x)= Jx$ given that $i\in \{1,2\}$ is sampled uniformly, and has a unique solution $x_\star = 0$.
By identifying $\reals^2$ with $\mathbb{C}$ via $x=(x_1,x_2) \longleftrightarrow z(x) = x_1 + \mathbf{i} x_2$ where $\mathbf{i} = \sqrt{-1}$, the S-SEG iterates $x_k$ can be identified with $z_k = z(x_k) \in \mathbb C$ generated by the following rule:
\begin{itemize}
    \item For S-SEG, $\hat z_k = \left(1 - \extstep(s_k\rho + \imu(1+s_k\nu))\right)z_k$ and $z_{k+1} = z_k \cdot b_{s_k}(\gamma_k, \extstep)$ where $s_k \in \{+1,-1\}$ is sampled uniformly and 
    \[
        b_s (\gamma, \hat{\gamma}) = 1 + \gamma \bigl((\rho^2-(1+s\nu)^2)\hat\gamma - s\rho\bigr)
        + \imu\gamma \left( (1+s\nu)(2s\rho\hat\gamma - 1) \right).
    \]
    \item For I-SEG, $\hat z_k = \left(1 - \extstep(u_k\rho + \imu(1+u_k\nu))\right)z_k$ and $z_{k+1} = z_k \cdot c_{u_k, s_k}(\gamma_k, \extstep)$ where $u_k, s_k \in \{+1,-1\}$ are sampled uniformly and independently, and 
    \[
        c_{u,s} (\gamma, \hat{\gamma}) = 1 + \gamma \bigl( (\rho^2 su - (1+s\nu)(1+u\nu)) \hat\gamma - s\rho \bigr) + \imu\gamma \left( \rho(s+u+2su\nu)\hat\gamma - (1+s\nu) \right) .
    \]
\end{itemize}
\end{lemma}

\begin{proof}
It is clear that $F(x) = Jx$ is monotone (\cref{assumption:monotonicity}) with unique zero $x_\star = 0$.
Additionally, both $F(\cdot; \xi_1)$ and $F(\cdot; \xi_2)$ are Lipschitz continuous with
\[
    L = \max\left\{ \norm{\rho I + (1+\nu)J}, \norm{-\rho I + (1-\nu)J} \right\}
    = \sqrt{\rho^2+(1+|\nu|)^2} 
\]
so \cref{assumption:uniform-lipschitzness} holds.
Next, observe that under the identification of $x = (x_1, x_2)$ by $z(x) = x_1 + \imu x_2$, the operator $J$ acts as multiplication by $\mathbf{i} = \sqrt{-1}$ because $\imu z(x) = \imu (x_1 + \imu x_2) = -x_2 + \imu x_1 = z(-x_2, x_1) = z(Jx)$.
Now for S-SEG, let $s_k \in \{+1,-1\}$ denote the sampled sign at iteration $k$, with $s_k = +1$ corresponding to $\xi_1$ and $s_k = -1$ to $\xi_2$.
Then the update rule, written in terms of the complex number, is given by
\begin{align*}
    \hat{z}_k & = z_k - \extstep (s_k\rho + \imu(1+s_k\nu)) z_k = (1 - \extstep (s_k\rho + \imu(1+s_k\nu)) ) z_k \\
    z_{k+1} & = z_k - \gamma_k (s_k\rho + \imu(1+s_k\nu)) \hat{z}_k = \left( 1 - \gamma_k (s_k\rho + \imu(1+s_k\nu))(1 - \extstep (s_k\rho + \imu(1+s_k\nu)) ) \right) z_k
\end{align*}
and rearranging the last factor and using $s_k^2 = 1$ yields the claimed form of $b_{s_k}(\gamma_k, \extstep)$.
Similarly, for I-SEG, let $u_k, s_k \in \{+1,-1\}$ respectively denote the sampled sign for the extrapolation and update steps.
Then the recursion is given by 
\begin{align*}
    \hat{z}_k & = z_k - \extstep (u_k\rho + \imu(1+u_k\nu)) z_k = (1 - \extstep (u_k\rho + \imu(1+u_k\nu)) ) z_k \\
    z_{k+1} & = z_k - \gamma_k (s_k\rho + \imu(1+s_k\nu)) \hat{z}_k = \left( 1 - \gamma_k (s_k\rho + \imu(1+s_k\nu))(1 - \extstep (u_k\rho + \imu(1+u_k\nu)) ) \right) z_k 
\end{align*}
and rearranging the last factor, we obtain the form of $c_{u_k,s_k}(\gamma_k, \extstep)$.

\end{proof}

\begin{proof}[\textbf{Proof of \cref{theorem:same-sample-unbounded-divergence-nonzero-mean}.}]
Consider the stochastic VIP of \cref{lemma:stochastic-vip-counterexample} with $\rho=2$, $\nu=0$ and $\extstep = \gamma_k$, and let $z_k, \hat{z}_k$ be complex numbers corresponding to the S-SEG iterates $x_k, \hat{x}_k$.
We can write $z_{k+1} = b_{s_k}(\gamma_k) z_k$ where we write
\[
    b_s(\gamma) := b_s(\gamma, \gamma) = 1 - \gamma(2s+\imu) + \gamma^2 (2s + \imu)^2 = 
    (1 - 2s\gamma + 3\gamma^2) - \imu (\gamma - 4s\gamma^2)  .
\]
By direct computation, we obtain
\[
    |b_s(\gamma)|^2 = 1-4s\gamma+11\gamma^2-20s\gamma^3 + 25\gamma^4 .
\]
Hence, by plugging $t = |b_s (\gamma, \gamma)|^2 - 1 = -4s\gamma + 11\gamma^2 - 20s\gamma^3 + 25\gamma^4$ into the Taylor expansion $\frac12\log(1+t) = \frac{1}{2} \left(t - \frac12 t^2 + \frac13 t^3 - \cdots \right)$, we obtain
\begin{align*}
    \log|b_s(\gamma)| = \frac{1}{2} \log (1+t) = \frac{1}{2} t - \frac{1}{4} t^2 + \cO\left( t^3 \right) = -2s\gamma + \frac{3}{2} \gamma^2 + \cO\left( \gamma^3 \right)
\end{align*}
provided that $|t| < 1$.
Hence, if $\gamma < \frac{1}{8}$, there exists a constant $C > 0$ such that
\begin{align}
\label{eqn:same-sample-log-modulus-expansion}
    \left| \log |b_s(\gamma)| + 2s\gamma -\frac32\gamma^2 \right| \le C\gamma^3 , \quad \left| \log |b_s(\gamma)| \right| \le C\gamma .
\end{align}
Additionally, 
\begin{align*}
    \tan \arg b_s(\gamma) = -\frac{\gamma - 4s\gamma^2}{1 - 2s\gamma + 3\gamma^2} = -\gamma (1 - 4s\gamma) \left[ 1 + (2s\gamma - 3\gamma^2) + (2s\gamma - 3\gamma^2)^2 + \cdots \right] = -\gamma + 2s\gamma^2 + \cO\left(\gamma^3 \right)
\end{align*}
where we use the Taylor expansion $\frac{1}{1-t} = 1 + t + t^2 + \cdots$ with $t=2s\gamma - 2\gamma^2$ which holds when $|t|<1$, and in particular, if $\gamma < \frac{1}{8}$.
Because $\arctan t = t - \frac{t^3}{3} + \frac{t^5}{5} - \cdots$, this implies $\arg b_s(\gamma) = -\gamma + 2s\gamma^2 + \cO\left(\gamma^3 \right)$, so by possibly increasing the constant $C$ in \eqref{eqn:same-sample-log-modulus-expansion}, we have
\begin{align}
\label{eqn:same-sample-argument-expansion}
    \left| \arg b_s(\gamma) + \gamma - 2s\gamma^2 \right| \le C\gamma^3 , \quad
    \left| \arg b_s(\gamma) \right| \le C\gamma 
\end{align}
for $\gamma < \frac{1}{8}$.
Since $\gamma_k = \frac{1}{\sqrt{k+1}} \to 0$, there exists $K\in \mathbb{N}$ such that conditions \eqref{eqn:same-sample-log-modulus-expansion} and \eqref{eqn:same-sample-argument-expansion} will hold with $\gamma=\gamma_k$ for all $k\ge K$.
Then, using those bounds and $s_k^2 = 1$, we obtain
\begin{align}
    \log |z_{N+1}| & = \log |z_K| + \sum_{k=K}^N \log |b_{s_k}(\gamma_k)| \nonumber \\
    & \ge \log |z_K| - 2\sum_{k=K}^N \frac{s_k}{\sqrt{k+1}} + \frac32\sum_{k=K}^N \frac{1}{k+1} - \sum_{k=K}^N \frac{C}{(k+1)^{3/2}} \nonumber \\
    & \ge \log |z_K| - 2\sum_{k=K}^N \frac{s_k}{\sqrt{k+1}} + \frac32\sum_{k=K}^N \frac{1}{k+1} - 3C
    \label{eqn:zN-norm-bound}
\end{align}
where the last inequality follows from $\sum_{k=K}^N (k+1)^{-3/2} < \sum_{k=1}^\infty k^{-3/2} \le 1 + \int_1^\infty x^{-3/2} \, dx = 3$.
Now the $Y_k = \frac{s_k}{\sqrt{k+1}}$ are mutually independent and zero-mean random variables, and
\begin{align*}
    \mathrm{Var}\left( \frac{Y_k}{\log(k+2)} \right) = \frac{\mathrm{Var}(Y_k)}{\log^2 (k+2)} = \frac{1}{(k+1) \log^2 (k+2)}
\end{align*}
is summable, so by Kolmogorov's two-series theorem, $\sum_{k=0}^\infty \frac{Y_k}{\log(k+2)}$ converges almost surely.
On that event, by Kronecker's lemma, we have $\frac{1}{\log(N+2)} \sum_{k=0}^N Y_k \to 0$.
In other words, $\sum_{k=0}^N Y_k = \sum_{k=0}^N \frac{s_k}{\sqrt{k+1}} = o(\log N)$ almost surely.
Together with \eqref{eqn:zN-norm-bound}, this implies
\begin{align*}
    \log |z_{N+1}| \ge \log |z_K| - 3C - o(\log N)  + \frac{3}{2} \sum_{k=K}^N \frac{1}{k+1} = \frac{3}{2} \log N - o(\log N) = \log N \left(\frac{3}{2} - o(1) \right)
\end{align*}
and thus $\norm{x_N} = |z_N| \ge N^{\frac{3}{2} - o(1)} \to \infty$ almost surely.

Next, we turn to the ergodic average of the extrapolated iterates; let
\[
    \widehat{w}_N = \sum_{k=0}^N \hat{z}_k = (N+1)\overline{\hat{z}}_N,
    \qquad
    \overline{\hat{z}}_N = \frac{1}{N+1}\sum_{k=0}^N \hat{z}_k .
\]
We will show that the limit superior of $|\overline{\hat{z}}_N| = \norm{\overline{\hat{x}}_N}$ is infinity, where
$\overline{\hat{x}}_N = \frac{1}{N+1}\sum_{k=0}^N \hat{x}_k$.
Given $N \in \mathbb{N}$, let $m_N = \lfloor N^{1/4}\rfloor$.
For $N$ sufficiently large so that $N-m_N \ge K$, for any $N-m_N \le k \le N$, using the bound $\left| \log |b_s (\gamma)| \right| \le C\gamma$ from \eqref{eqn:same-sample-log-modulus-expansion} we have
\begin{align*}
\left|\log \frac{|z_k|}{|z_N|} \right|
&\le \sum_{j=k}^{N-1} \left| \log |b_{s_j}(\gamma_j)| \right|
\le \sum_{j=k}^{N-1} \frac{C}{\sqrt{j+1}}
\le \frac{C m_N}{\sqrt{N-m_N}}
= \cO(N^{-1/4}) .
\end{align*}
Similarly, using \eqref{eqn:same-sample-argument-expansion},
\begin{align*}
\left| \arg z_k - \arg z_N \right|
&\le \sum_{j=k}^{N-1} \left| \arg b_{s_j}(\gamma_j) \right|
\le \sum_{j=k}^{N-1} \frac{C}{\sqrt{j+1}}
\le \frac{C m_N}{\sqrt{N-m_N}}
= \cO(N^{-1/4}) .
\end{align*}
Therefore, for all sufficiently large $N$, if $N-m_N \le k \le N$, then
\begin{align}
\label{eqn:zk-zN-magnitude-and-argument-bound-hat-average}
    \frac12 |z_N| \le |z_k| \le 2|z_N|,
    \qquad
    \left| \arg z_k - \arg z_N \right| \le \frac{\pi}{6}.
\end{align}
Moreover, since
\[
    \hat{z}_k = \left(1-\gamma_k(2s_k+\imu)\right)z_k 
\]
and $\gamma_k \to 0$, by increasing $K$ if necessary, we may also assume that, for all $k\ge K$,
\[
    \left|1-\gamma_k(2s_k+\imu)\right| \ge \frac12,
    \qquad
    \left|\arg\left(1-\gamma_k(2s_k+\imu)\right)\right| \le \frac{\pi}{6}.
\]
Hence, for all sufficiently large $N$ and all $N-m_N \le k \le N$,
\[
    |\hat{z}_k| \ge \frac12 |z_k| \ge \frac14 |z_N|,
    \qquad
    \left|\arg \hat{z}_k - \arg z_N\right| \le \frac{\pi}{3}.
\]
Now consider the partial sum
\[
    \widehat{d}_N := \sum_{k=N-m_N}^N \hat{z}_k
    = \widehat{w}_N - \widehat{w}_{N-m_N-1}.
\]
Let $u_N = z_N/|z_N|$. Then, for all sufficiently large $N$ and all $k=N-m_N,\dots,N$,
\[
    \Re\left(u_N^* \hat{z}_k\right)
    =
    |\hat{z}_k|\cos(\arg \hat{z}_k-\arg z_N)
    \ge \frac12 |\hat{z}_k|
    \ge \frac18 |z_N|,
\]
where $u_N^*$ denotes the complex conjugate. Therefore,
\[
    |\widehat{w}_N-\widehat{w}_{N-m_N-1}|
    =
    |\widehat{d}_N|
    \ge
    \Re\left(u_N^* \widehat{d}_N\right)
    =
    \sum_{k=N-m_N}^N \Re\left(u_N^* \hat{z}_k\right)
    \ge
    \frac{m_N |z_N|}{8}.
\]
This implies
\[
    \max\{|\widehat{w}_N|,|\widehat{w}_{N-m_N-1}|\}
    \ge
    \frac{|\widehat{d}_N|}{2}
    \ge
    \frac{m_N |z_N|}{16}.
\]
Dividing by the appropriate denominators gives
\begin{align*}
    \max\{|\overline{\hat{z}}_N|,|\overline{\hat{z}}_{N-m_N-1}|\}
    &=
    \max\left\{
        \frac{|\widehat{w}_N|}{N+1},
        \frac{|\widehat{w}_{N-m_N-1}|}{N-m_N}
    \right\} \\
    &\ge
    \frac{1}{N+1}
    \max\{|\widehat{w}_N|,|\widehat{w}_{N-m_N-1}|\}
    \ge
    \frac{m_N |z_N|}{16(N+1)}.
\end{align*}
Because $m_N = \Theta(N^{1/4})$ and $|z_N| = N^{3/2-o(1)}$ almost surely, the right-hand side diverges to infinity almost surely. Therefore,
\[
    \limsup_{N\to\infty} |\overline{\hat{z}}_N| = \infty,
\]
and equivalently,
\[
    \limsup_{N\to\infty} \norm{\overline{\hat{x}}_N} = \infty
    \qquad \text{a.s.}
\]
Because $\err_{x_\star, D}(x) = D\norm{x}$ for any $x$ for this problem, we conclude that
\[
    \err_{x_\star, D}(x_k) \to \infty,
    \qquad
    \limsup_{k\to\infty} \err_{x_\star, D}(\overline{\hat{x}}_k) = \infty .
\]

\end{proof}

\paragraph{I-SEG dynamics on the same problem converges.}
When I-SEG is instead run on the same stochastic problem, $x_k$ converges to $x_\star = 0$ almost surely, unlike S-SEG. 
See Figure~\ref{fig:sseg-counterexample-plane-trajectory} for visual comparison.
To check this analytically, applying \cref{lemma:stochastic-vip-counterexample} and plugging in $\gamma_k = \extstep$, we have $z_{k+1} = c_{u_k, s_k} (\gamma_k) z_k$ with
\begin{align}
\label{eqn:I-SEG-factor-on-S-SEG-counterexample}
    c_{u,s}(\gamma) := c_{u,s}(\gamma, \gamma) = 1 - 2s\gamma + (4us - 1)\gamma^2 - \imu \left( \gamma - 2(u+s) \gamma^2 \right)
\end{align}
and using $u^2 = s^2 = 1$, we obtain
\begin{align}
    & |c_{u,s}(\gamma)|^2 = 1 - 4s\gamma + (3 + 8us) \gamma^2 - 20u \gamma^3 + 25\gamma^4 \notag \\
    & \implies \log |c_{u,s}(\gamma)| = -2s\gamma + \left(-\frac52 + 4u s\right)\gamma^2 + \cO(\gamma^3) .
\label{eqn:I-SEG-log-factor-on-S-SEG-counterexample}
\end{align}
Here the crucial difference is that the expectation of the second term is $-\frac{5}{2}\gamma^2$ has a negative coefficient, while for S-SEG, the coefficient of $\gamma^2$ is a positive constant.
Therefore, taking the logarithm we obtain 
\[
    \log |z_{N+1}| = \log |z_K| -2\sum_{k=K}^N \frac{s_k}{\sqrt{k+1}} - \frac52\sum_{k=K}^N \frac{1}{k+1} + 4\sum_{k=K}^N \frac{u_k s_k}{k+1} + \cO(1) .
\]
As before, we have $\sum_{k=K}^N \frac{s_k}{\sqrt{k+1}} = o(\log N)$ almost surely, $-\frac52\sum_{k=K}^N \frac{1}{k+1} = -\frac52 \log N + \cO(1)$, and $\sum_{k=K}^N \frac{u_k s_k}{k+1}$ converges almost surely by Kolmogorov's two-series as it is a sum of independent, zero-mean and square summable random variables. 
This shows that $z_N \to 0$ almost surely, which implies that both $x_k$ and $\overline{\hat{x}}_k$ converges to $0$ almost surely.

It is worth noting that, on the other hand, 
$\expec{\sqnorm{x_{k+1}} \,\middle|\, \cF_k}{} = \left(1+3\gamma_k^2+25\gamma_k^4\right) \sqnorm{x_k}$ for I-SEG run on this example.
That is, even if the expectation of $\sqnorm{x_k}$ strictly increases over iterations, we still have $x_k \to 0$ almost surely.
There is no contradiction here---expectation of random variables may increase, and even blow up, even if the random variable itself converges to $0$ almost surely: take, for example, the sample space $\Omega = [0,1]$ and random variables $X_n(\omega) = n^2 \mathbf{1}_{(0,1/n]}(\omega)$.
Hence, showing that $\sqnorm{x_k - x_\star}$ increases in expectation, as done in several prior work \citep{ChavdarovaGidelFleuretLacoste-Julien2019_reducing, HsiehIutzelerMalickMertikopoulos2020_explore, chaeStochasticExtragradientFlipflop2024}, is by itself not sufficient to guarantee that the algorithm's dynamics is truly divergent.
In the proof of \cref{theorem:independent-sample-unbounded-divergence-generalized-variance} below, we demonstrate a counterexample where I-SEG diverges almost surely.
\medskip

\begin{proof}[\textbf{Proof of \cref{theorem:independent-sample-unbounded-divergence-generalized-variance}}]
Consider the stochastic VIP from \cref{lemma:stochastic-vip-counterexample} with $\rho=0$ and $\nu=2$.
Then, because $\expec{\sqnorm{F(x;\xi)-F(x)}}{}=4\sqnorm{x}$, \eqref{eqn:generalized-variance-bound} holds with $A=4$ and $\sigma=0$.
By \cref{lemma:stochastic-vip-counterexample}, for $s_k, u_k \in \{+1, -1\}$, we have
\[
    \hat{z}_k = \left( 1-\imu \gamma_k(1+2u_k) \right) z_k, \qquad
    z_{k+1} = c_{u_k,s_k} \left( \gamma_k \right) z_k,
\]
where
\[
    c_{u,s}(\gamma):= c_{u,s}(\gamma,\gamma) = 1 - \gamma^2 (1+2s)(1+2u) - \imu \gamma (1+2s) .
\]
Then, using $u^2 = s^2 = 1$, we obtain $|c_{u,s}(\gamma)|^2 = 1 + \left( 3 - 4u - 8su \right) \gamma^2 + \cO\left(\gamma^4 \right)$.
Thus, for $\gamma$ small enough,
\[
    \log |c_{u,s}(\gamma)|
    =\left(\frac32-2u-4su\right)\gamma^2+\cO(\gamma^4) 
\]
and
\[
    |\arg c_{u,s}(\gamma)| = \left| \arctan \frac{-\gamma(1+2s)}{1-\gamma^2 (1+2s) (1+2u)} \right| = \left| -\gamma(1+2s) + \cO\pr{\gamma^3} \right| \le C\gamma
\]
for some fixed constant $C>0$.
Now, we have
\begin{align*}
    \log|z_N| = \log|z_0| + \sum_{k=0}^{N-1} \log |c_{u_k,s_k}(\gamma_k)| 
    = \log|z_0| + \sum_{k=0}^{N-1} \left(\frac32-2u_k-4s_k u_k\right) \gamma_k^2 + \sum_{k=0}^{N-1} \cO\pr{\gamma_k^4} ,
\end{align*}
and the last summation almost surely converges to a finite constant because $\sum_{k=0}^\infty \gamma_k^4<\infty$.
On the other hand, $\left(\frac32-2u_k-4s_k u_k\right) \gamma_k^2$ has mean $\frac{3}{2}\gamma_k^2$ and summable $\cO(\gamma_k^4)$ variance, 
so by Kolmogorov's two-series theorem, we have
\[
    \log |z_N|= \frac32 \sum_{k=0}^N \gamma_k^2 + \cO(1) = \frac32 \log N+\cO(1) .
\]
Hence $\norm{x_N}=|z_N|\to\infty$ almost surely.

It remains to show that the ergodic average also diverges along a subsequence.
Let $m_N=\lfloor N^{1/4}\rfloor$. Then for $N$ sufficiently large and $N-m_N\le k\le N$, 
similar argument as in the proof of \cref{theorem:same-sample-unbounded-divergence-nonzero-mean}, using
\[
    \left|\log\frac{|z_k|}{|z_N|}\right|
    \le C\sum_{j=k}^{N-1}\gamma_j^2=\cO(N^{-3/4}),
    \qquad
    |\arg z_k-\arg z_N|
    \le C\sum_{j=k}^{N-1}\gamma_j=\cO(N^{-1/4}) ,
    \qquad
    \hat{z}_k = \left(1 + \cO(\gamma_k)\right) z_k
\]
yields
\[
    \max\{\norm{\overline{\hat{x}}_N},\norm{\overline{\hat{x}}_{N-m_N-1}}\}
    \ge \frac{m_N\norm{x_N}}{16(N+1)} .
\]
Since $m_N=\Theta(N^{1/4})$ and $\norm{x_N}=N^{3/2+o(1)}$ almost surely, the right-hand side diverges almost surely, showing that
$\limsup_{k\to\infty} \norm{\overline{\hat{x}}_k} = \infty$.
\end{proof}

Interestingly, although we do not formally argue it here, the S-SEG dynamics on the same problem considered in the above proof converges (Figure~\ref{fig:iseg-counterexample-plane-trajectory}).
This indicates that neither I-SEG nor S-SEG is not strictly preferred over the other, and they could behave differently depending on the geometry of the stochastic problem.

Finally, we show that even when $\gamma_k = \extstep$ is square summable, we cannot have any uniform gap function guarantee.
For the stochastic monotone VIP considered in the proof of \cref{theorem:same-sample-unbounded-divergence-nonzero-mean}, for any $M>0$, albeit with low probability scaling down as $M$ grows, we have $\limsup_{k\to\infty} \err_{x_\star,D}(\overline{\hat{x}}_k) > M$ with both types of SEG.
In this sense, the high probability bounds established in Theorems~\ref{theorem:independent-sample-SEG-unbounded-domain-high-probability-bound} and \ref{theorem:same-sample-SEG-unbounded-domain-high-probability-bound} are tight results for SEG using symmetric step-sizes in unbounded domains.

\begin{theorem}
\label{theorem:unbounded-gap-high-probability-tightness}
Consider the stochastic VIP from
Theorem~\ref{theorem:same-sample-unbounded-divergence-nonzero-mean}.
If SEG (either same-sample or independent-sample) 
is run on this problem with initial point $x_0 \ne 0$ and $\extstep = \gamma_k$ satisfying $0 < \gamma_k < \overline{\gamma}$ for all $k\ge 0$ and sufficiently small $\overline{\gamma}$ and $\sum_{k=0}^\infty \gamma_k^2 < \infty$, then there exists a positive random variable $R_\infty$ such that $\norm{x_k}\to R_\infty$ almost surely.
In particular, if $\extstep = \gamma_k = \frac{\eta}{k+1}$ where $\eta > 0$ is sufficiently small, then we also have
\[
    \norm{\overline{x}_k}\to \frac{R_\infty}{\sqrt{1+\eta^2}},
    \qquad
    \norm{\overline{\hat{x}}_k}\to \frac{R_\infty}{\sqrt{1+\eta^2}}
\]
almost surely, where $\overline{x}_k = \frac{1}{k+1}\sum_{j=0}^k x_j$ and $\overline{\hat{x}}_k = \frac{1}{k+1}\sum_{j=0}^k \hat{x}_j$.
Moreover, the law of $R_\infty$ is unbounded, and consequently, for every $D>0$,
\begin{align*}
    \mathrm{Prob} \left[\lim_{k\to\infty}\err_{x_\star,D}(\overline{x}_k) > M\right] > 0,
    \qquad
    \mathrm{Prob} \left[\lim_{k\to\infty}\err_{x_\star,D}(\overline{\hat{x}}_k) > M\right] > 0,
    \quad \forall M > 0 .
\end{align*}
\end{theorem}

\begin{proof}
We first provide the proof for S-SEG. The I-SEG case follows from the essentially same arguments with minor modifications, 
which we specify at the end of the proof.

Using \cref{lemma:stochastic-vip-counterexample} again, write $z_{k+1} = b_{s_k}(\gamma_k)z_k$ where $s_k\in\{+1,-1\}$ is chosen uniformly at each iteration.
Because $b_{s_k} (\gamma_k) = 1 + \cO(\gamma_k)$, for $\gamma_k$ sufficiently small, we have $\left| b_{s_k} (\gamma_k) - 1 \right| \le \frac{1}{2}$.
Then, using \eqref{eqn:same-sample-log-modulus-expansion} and \eqref{eqn:same-sample-argument-expansion}, we obtain
\begin{align}
    \log \left| b_{s_k}(\gamma_k) \right| & = -2s_k \gamma_k + \frac{3}{2} \gamma_k^2 + \cO\left(\gamma_k^3 \right)
    \label{eqn:harmonic-log-rho}
    \\
    \arg b_{s_k}(\gamma_k) & = -\gamma_k + 2s_k \gamma_k^2 + \cO\left( \gamma_k^3 \right) .
    \label{eqn:harmonic-arg}
\end{align}
In \eqref{eqn:harmonic-log-rho}, the series
$\sum_{k=0}^\infty -2s_k \gamma_k$ converges almost surely because each term has mean zero and their variances $\mathrm{Var}\left(-2s_k \gamma_k\right) = 4\gamma_k^2$ is summable.
On that event, $\log |z_{N+1}| = \log |z_0| +\sum_{k=0}^N \log \left| b_{s_k}(\gamma_k) \right|$ converges to a real limit because $\sum_{k=0}^\infty \gamma_k^2$ and $\sum_{k=0}^\infty \gamma_k^3$ are finite, and thus $\norm{x_k} = |z_k| \to R_\infty$ almost surely for some finite, positive random variable $R_\infty > 0$.

Next, fix our step-size selection to be $\extstep = \gamma_k = \frac{\eta}{k+1}$ with $\eta > 0$ small enough so that $|b_{s_0}(\eta) - 1| \le \frac{1}{2}$ for both $s_0 = \pm 1$.
Let $H_N = \sum_{k=0}^{N-1}\frac{1}{k+1}$ be the partial sum of the harmonic series.
Then
\begin{align*}
    \arg z_N + \eta H_N = \arg z_0 + \sum_{k=0}^{N-1} \arg b_{s_k}(\gamma_k) + \sum_{k=0}^{N-1} \frac{\eta}{k+1} = \arg z_0 + \sum_{k=0}^{N-1} \frac{2\eta^2 s_k}{(k+1)^2} + \cO \left( \frac{1}{(k+1)^3} \right)
\end{align*}
by \eqref{eqn:harmonic-arg}, and the series in the right hand side converge absolutely, say, to $\theta_\infty$.
Therefore, we obtain
Since $|z_N|\to R_\infty$, it follows that
\[
    z_N e^{\imu \eta H_N} = |z_N| e^{\imu \left( \arg z_N + \eta H_N\right)} \to R_\infty e^{\imu \theta_\infty}
\]
almost surely.
It is well-known that $\lim_{N\to\infty} \left( H_N - \log N \right) = \lim_{N\to\infty} \left( H_N - \log (N+1) \right) = \gamma_\text{Euler} \approx 0.577$ (the Euler constant), so we have
\begin{align*}
    z_N e^{\imu \eta H_N} = z_N e^{\imu \eta (\log (N+1) + \gamma_\text{Euler}) + o(1)} = z_N e^{\imu \eta \gamma_\text{Euler}} (N+1)^{\imu \eta} (1 + o(1)) \to R_\infty e^{\imu \theta_\infty} 
\end{align*}
so letting $A_\infty = R_\infty e^{\imu \left(\theta_\infty - \eta \gamma_\text{Euler}\right)}$, we can write
\[
    z_N = A_\infty (N+1)^{-\imu \eta} + |z_N| \, o(1) = A_\infty (N+1)^{-\imu \eta} + o(1) = A_\infty (N+1)^{-\imu \eta} + r_N
\]
almost surely, where $r_N = z_N - A_\infty (N+1)^{-\imu \eta} = |z_N| \, o(1) = o(1)$ because $|z_N| \to R_\infty$.
Then the sequence $r_k$ is Ces\`aro summable with $\frac{1}{N+1}\sum_{k=0}^N r_k \to 0$ as $N\to\infty$.
Hence
\[
    \overline{z}_N = \frac{1}{N+1} \sum_{k=0}^N z_k = \frac{A_\infty}{N+1} \sum_{k=0}^N (k+1)^{-\imu\eta} + \frac{1}{N+1} \sum_{k=0}^N r_k = \frac{A_\infty}{N+1} \sum_{k=0}^N (k+1)^{-\imu\eta} + o(1).
\]
Now observe that 
\[
    \frac{1}{N+1}\sum_{k=0}^N (k+1)^{-\imu\eta} = (N+1)^{-\imu\eta} \frac{1}{N+1}\sum_{k=1}^{N+1} \left(\frac{k}{N+1}\right)^{-\imu\eta} = (N+1)^{-\imu \eta} \left[ \int_0^1 x^{-\imu\eta} \, dx + o(1) \right] = \frac{(N+1)^{-\imu\eta}}{1-\imu\eta} + o(1) ,
\]
where the last equality uses the fact that $\int_0^1 x^{-\imu\eta} \, dx = \frac{1}{1-\imu \eta}$ and $\left|(N+1)^{-\imu \eta}\right| = 1$.
The second equality follows from the convergence of the Riemann sum representing the improper integral $\int_0^1 x^{-\imu \eta} \, dx$, which exists because $\left|x^{-\imu\eta}\right| = 1$ for any $x>0$.
Combining the above, we see that 
\[
    \overline{z}_N = \frac{A_\infty (N+1)^{-\imu\eta}}{1-\imu\eta} + o(1) \implies \lim_{N\to\infty} \norm{\overline{x}_N} = \lim_{N\to\infty} |\overline{z}_N| = \frac{|A_\infty|}{|1-\imu\eta|} = \frac{R_\infty}{\sqrt{1 + \eta^2}} 
\]
almost surely.
The same limit holds for the ergodic average of the extrapolated iterates.
Indeed, since
\[
    \hat{z}_k = \left(1-\gamma_k(2s_k+\imu)\right)z_k
    = z_k - \frac{\eta}{k+1}(2s_k+\imu)z_k ,
\]
we have
\begin{align*}
    \left|\overline{\hat{z}}_N-\overline{z}_N\right|
    &=
    \left|
    \frac{1}{N+1}\sum_{k=0}^N(\hat{z}_k-z_k)
    \right| \le
    \frac{\eta\sqrt{5}}{N+1}\sum_{k=0}^N \frac{|z_k|}{k+1}.
\end{align*}
Since $|z_k|\to R_\infty$ almost surely, the sequence $(|z_k|)_{k\ge 0}$ is almost surely bounded. Hence, almost surely,
\[
    \left|\overline{\hat{z}}_N-\overline{z}_N\right|
    \le
    \frac{\eta\sqrt{5}\sup_{k\ge 0}|z_k|}{N+1}
    \sum_{k=0}^N \frac{1}{k+1}
    \to 0 
\]
as $N \to \infty$.
Therefore,
\[
    \lim_{N\to\infty}\norm{\overline{\hat{x}}_N}
    =
    \lim_{N\to\infty}|\overline{\hat{z}}_N|
    =
    \lim_{N\to\infty}|\overline{z}_N|
    =
    \frac{R_\infty}{\sqrt{1+\eta^2}}
\]
almost surely.

Finally, we show that the law of $R_\infty$ is unbounded.
For $m=1,2,\dots$, consider 
\[
    \log T_m := \sum_{k=m}^{\infty} \log \left| b_{s_k}(\gamma_k) \right| = -2\eta \sum_{k=m}^{\infty} \frac{s_k}{k+1} + o(1)
\]
where the last expression follows from \eqref{eqn:harmonic-log-rho}, which shows that the residual term is $\cO\left( \sum_{k=m}^\infty \frac{1}{(k+1)^2} + \frac{1}{(k+1)^3} \right)$, which converges to zero as $m\to\infty$.
But because $\sum_{k=0}^\infty \frac{s_k}{k+1}$ is almost surely convergent, we have $\sum_{k=m}^\infty \frac{s_k}{k+1} \to 0$ almost surely as $m\to\infty$.
Therefore, we have $T_m \to 1$ almost surely, and because
\[
    \{T_m \to 1\} \subset \liminf_{m\to\infty} \left\{ T_m \ge \frac{1}{2} \right\} = \bigcup_{m=1}^\infty \bigcap_{j\ge m} \left\{ T_j \ge \frac{1}{2} \right\}
\]
this in particular shows that $\mathrm{Prob}[T_m\ge 1/2] \to 1$ as $m\to\infty$.
Next, consider the event $\cE_m = \{s_0 = \cdots = s_{m-1} = -1\}$. 
On $\cE_m$, we have
\[
    \log P_m := \sum_{k=0}^{m-1} \log \left| b_{s_k}(\gamma_k) \right| = 2\eta \sum_{k=0}^{m-1} \frac{1}{k+1} + \cO(1) = 2\eta H_m + \cO(1) 
\]
and the last quantity grows arbitrarily large with $m$.
Hence, for any given $M>0$, we can choose $m$ that is large enough so that 
\begin{align*}
    \mathrm{Prob}[T_m\ge 1/2] \ge \frac12 \quad \text{ and } \quad \frac{|z_0| P_m}{2\sqrt{1+\eta^2}} > M
     \, \text{ on } \, \cE_m .
\end{align*}
Additionally, because $\{T_m \ge 1/2\}$ is independent of $\cE_m$, we have
\begin{align*}
    \mathrm{Prob} \left[ \cE_m\cap\{T_m\ge 1/2\}\right] = \mathrm{Prob} \left[ \cE_m \right] \cdot \mathrm{Prob}[T_m\ge 1/2] \ge \frac{1}{2^{m+1}}
\end{align*}
and on that event,
\begin{align*}
    R_\infty = |z_0| P_m T_m \ge \frac{|z_0| P_m}{2} > \sqrt{1+\eta^2} M \implies \lim_{N\to\infty} \norm{\overline{\hat{x}}_N} = \frac{R_\infty}{\sqrt{1+\eta^2}} > M .
\end{align*}
Note that this can be done for any $M$; hence, given any $M,D>0$, we can choose an event of positive probability on which $\lim_{N\to\infty} \norm{\overline{\hat{x}}_N} > \frac{M}{D}$, so that
\[
    \lim_{N\to\infty} \err_{x_\star,D}(\overline{\hat{x}}_N) = \lim_{N\to\infty} D\norm{\overline{\hat{x}}_N} > M .
\]

\paragraph{Case of I-SEG.}
\cref{lemma:stochastic-vip-counterexample} yields $z_{k+1}=c_{u_k,s_k}(\gamma_k)z_k$, where $c_{u,s}(\gamma)$ is given by \eqref{eqn:I-SEG-factor-on-S-SEG-counterexample}, and as in \eqref{eqn:I-SEG-log-factor-on-S-SEG-counterexample} we can derive 
$\log|c_{u,s}(\gamma)| = -2s\gamma+\left(-\frac52+4us\right)\gamma^2+\cO(\gamma^3)$. 
Furthermore, $\arg c_{u,s}(\gamma)=-\gamma+2u\gamma^2+\cO(\gamma^3)$, so the first-order terms in $\gamma$ are the same as in the S-SEG case,
so we can proceed in the exact same way.
Showing that $R_\infty$ has an unbounded law also uses the identical argument; we can take $\cE_m = \{s_0 = \cdots = s_{m-1} = -1\}$, and show that for any $M>0$, $R_\infty > M$ for sufficiently large $m$, uniformly over $u_0, \dots, u_{m-1}$.
\end{proof}

\subsubsection{Effect of DSEG step-sizes on the convergence of SEG}
\label{section:DSEG-gap-convergence}

The \ref{eqn:DSEG} step-size rule \citep{HsiehIutzelerMalickMertikopoulos2020_explore} in \cref{proposition:I-SEG-almost-sure-DSEG} was designed to ensure almost sure last-iterate convergence of 
I-SEG for unconstrained problems under generalized variance condition similar to~\eqref{eqn:generalized-variance-bound}.
Although \cite{HsiehIutzelerMalickMertikopoulos2020_explore} did not explicitly state it, as a consequence of their algorithm design and analysis, it is relatively straightforward to see that \ref{eqn:DSEG} step-size also induces the gap function convergence of I-SEG on averaged iterates (\cref{proposition:independent-sample-DSEG-restricted-gap}).
On the other hand, this is not the case for S-SEG, as we demonstrate in Section~\ref{subsection:DSEG-fails-to-fix-S-SEG}.

\begin{proposition}
\label{proposition:independent-sample-DSEG-restricted-gap}
Let $\cX=\reals^d$, and let $x_\star\in\cX_\star \cap B(y,D)$, so $F(x_\star)=0$. 
Suppose that Assumptions~\ref{assumption:monotonicity}, \ref{assumption:lipschitzness}, \ref{assumption:unbiasedness} and \ref{assumption:generalized-variance-bound} hold.
Then for every $y\in\reals^d$ and $D>0$,
I-SEG run with initial point $x_0 \in B(y,D)$ and $\gamma_k > 0$, $\extstep \in \left(0, \frac{1}{3L}\right)$ satisfying \eqref{eqn:DSEG} satisfies
\begin{align}
\label{eqn:independent-DSEG-gap-rate}
    \expec{\err_{B(y,D)}\left(\overline{\hat{x}}_N^{\gamma}\right)}{}
    \le \frac{C_{y,D}}{\Gamma_N} .
\end{align}
where $C_{y,D}>0$ is a constant depending on $y, D$ and the fixed problem data but not on $N$ and $x_0 \in B(y,D)$, while $\Gamma_N = \sum_{k=0}^N \gamma_k$ and $\overline{\hat{x}}_N^{\gamma} := \frac{1}{\Gamma_N}\sum_{k=0}^N \gamma_k \halfitr$.
In particular, if 
\[
    \gamma_k=\gamma_0(k+1)^{-a} , \qquad \extstep=\hat{\gamma}_0(k+1)^{-b}
\]
with $\frac{1}{2}<a<1, a+b\le 1$ and $a+2b>1$, then $\expec{\err_{B(y,D)}\left(\overline{\hat{x}}_N^{\gamma}\right)}{} = \order{N^{-(1-a)}}$.
\end{proposition}

\begin{proof}
For convenience, we write $G_k = F(x_k;\xi_{i_k}), \widehat{G}_k = F(\halfitr;\halfxi), \Delta_k:=F(\halfitr)-\widehat{G}_k$, and let $\eta_k = \gamma_k^2 + 2L\gamma_k\extstep^2$.
Note that because $\cX = \reals^d$, we have
\[
    x_{k+1}=x_k-\gamma_k \widehat{G}_k,
    \qquad
    \halfitr=x_k-\extstep G_k .
\]
We quickly establish some key estimates needed for the proof. First, we have
\begin{align}
\label{eqn:DSEG-gap-Gk-moment-estimate}
    \expec{\sqnorm{G_k}\,\middle|\,\cF_k}{}
    =\sqnorm{F(x_k)}+
      \expec{\sqnorm{G_k-F(x_k)}\,\middle|\,\cF_k}{}
    \le (L^2+A)\sqnorm{x_k-x_\star}+\sigma^2 
    \le C\left(1+\sqnorm{x_k-x_\star}\right) 
\end{align}
where we take $C \ge \max\left\{L^2 + A, \sigma^2 \right\}$.
Next, using $\condexp{G_k}{\cF_k} = F(x_k)$, we have
\begin{align}
\label{eqn:DSEG-gap-halfitr-distance-estimate}
    \expec{\sqnorm{\halfitr-x_\star}\,\middle|\,\cF_k}{} 
    = \sqnorm{x_k - x_\star} - 2\extstep \inprod{F(x_k)}{x_k - x_\star} + \extstep^2 \condexp{\sqnorm{G_k}}{\cF_k}
    \le C \left( 1 + \sqnorm{x_k - x_\star} \right)
\end{align}
by enlarging $C$ if necessary, because of \eqref{eqn:DSEG-gap-Gk-moment-estimate} and the fact that $\extstep$ is upper-bounded by $\frac{1}{3L}$.
This implies, again by enlarging $C$ if necessary,
\begin{align}
\label{eqn:DSEG-gap-Gk-hat-moment-estimate}
    \expec{\sqnorm{\widehat{G}_k}\,\middle|\,\cF_k}{} \le (L^2 + A) \condexp{\sqnorm{\halfitr - x_\star}}{\cF_k} + \sigma^2 \le C \left( 1 + \sqnorm{x_k - x_\star} \right) .
\end{align}
Now, we have
\begin{align}
    \notag
    \sqnorm{x_{k+1}-x_\star} &=
    \sqnorm{x_k-x_\star} -2\gamma_k\inprod{\widehat{G}_k}{x_k-x_\star} +\gamma_k^2\sqnorm{\widehat{G}_k} \\
    &=    
    \sqnorm{x_k-x_\star}
      -2\gamma_k\inprod{\widehat{G}_k}{\halfitr-x_\star}
      -2\gamma_k\extstep\inprod{\widehat{G}_k}{G_k}
      +\gamma_k^2\sqnorm{\widehat{G}_k} 
    \label{eqn:DSEG-gap-xkp1-distance-estimate}
\end{align}
where we use $x_k - \halfitr = \extstep G_k$.
Conditioning first on $\hat{\cF}_k$ and taking conditional expectation, the term
$\inprod{\widehat{G}_k}{\halfitr-x_\star}$ becomes $\inprod{F(\halfitr)}{\halfitr-x_\star}\ge 0$ so we can drop it.
Furthermore, by $L$-Lipschitzness of $F$,
\begin{align}
\label{eqn:DESG-gap-cross-term-bound}
    \expec{\inprod{\widehat{G}_k}{G_k}\,\middle|\,\cF_k}{}
    &=\expec{\inprod{F(\halfitr)}{G_k}\,\middle|\,\cF_k}{} 
    =\sqnorm{F(x_k)}
      +\expec{\inprod{F(\halfitr)-F(x_k)}{G_k}\,\middle|\,\cF_k}{} 
    \ge \sqnorm{F(x_k)}
      -L\extstep\expec{\sqnorm{G_k}\,\middle|\,\cF_k}{} .
\end{align}
Plugging this into \eqref{eqn:DSEG-gap-xkp1-distance-estimate} and also using \eqref{eqn:DSEG-gap-Gk-moment-estimate} and \eqref{eqn:DSEG-gap-Gk-hat-moment-estimate}, we obtain 
\begin{align*}
    \expec{\sqnorm{x_{k+1}-x_\star}\,\middle|\,\cF_k}{}
    & \le
    \sqnorm{x_k - x_\star} + 2L \gamma_k \extstep^2 C \left(1 + \sqnorm{x_k - x_\star}\right) + \gamma_k^2 C \pr{1 + \sqnorm{x_k - x_\star}} \\
    & =
    \left(1+C \eta_k\right)\sqnorm{x_k-x_\star} + C \eta_k .
\end{align*}
The \ref{eqn:DSEG} step-size condition implies that $\sum_{k=0}^\infty \eta_k < \infty$, so by \cref{lemma:almost-sure-sequence-lemma}(ii) we have
\begin{align}
\label{eqn:DSEG-distance-expec-bound}
    \sup_{k\ge 0}\expec{\sqnorm{x_k-x_\star}}{} \le \left( \sqnorm{x_0 - x_\star} + C\sum_{k=0}^\infty \eta_k \right) e^{C\sum_{k=0}^\infty \eta_k} := M < \infty .
\end{align}
Now fix $u\in B(y,D)$ and rearrange the following identity similar to \eqref{eqn:DSEG-gap-xkp1-distance-estimate}:
\begin{align*}
    \sqnorm{x_{k+1} - u} = \sqnorm{x_k-u}
      -2\gamma_k\inprod{\widehat{G}_k}{\halfitr-u}
      -2\gamma_k\extstep\inprod{\widehat{G}_k}{G_k}
      +\gamma_k^2\sqnorm{\widehat{G}_k} ,
\end{align*}
to obtain
\begin{align*}
    \gamma_k\inprod{F(u)}{\halfitr-u} & \le \gamma_k\inprod{F(\halfitr)}{\halfitr-u} \\
    & \le \frac12\left(\sqnorm{u-x_k}-\sqnorm{u-x_{k+1}}\right)
    +\frac{\gamma_k^2}{2}\sqnorm{\widehat{G}_k}
    -\gamma_k\extstep\inprod{\widehat{G}_k}{G_k}
    +\gamma_k\inprod{\Delta_k}{\halfitr-x_\star}
    +\gamma_k\inprod{\Delta_k}{x_\star-u} .
\end{align*}
Taking the sum of the above from $k=0$ to $N$ and using $\norm{x_\star - u} \le \norm{x_\star - y} + \norm{y - u} \le \norm{x_\star - y} + D := D_\star$, we have
\begin{align*}
    \Gamma_N \err_{B(y,D)}(\overline{\hat{x}}_N^\gamma) & = \max_{u \in B(y,D)} \sum_{k=0}^N \gamma_k\inprod{F(u)}{\halfitr-u} \\
    & \le \frac{1}{2} \sqnorm{u - x_0} + \sum_{k=0}^N
    \left(\frac{\gamma_k^2}{2}\sqnorm{\widehat{G}_k}
    -\gamma_k\extstep\inprod{\widehat{G}_k}{G_k}\right)
    +\sum_{k=0}^N\gamma_k\inprod{\Delta_k}{\halfitr-x_\star}
    +D_\star\norm{\sum_{k=0}^N\gamma_k\Delta_k} .
\end{align*} 
Here we have $\frac{1}{2}\sqnorm{u - x_0} \le 2D^2$ because $u,x_0 \in B(y,D)$, the third term $\sum_{k=0}^N\gamma_k\inprod{\Delta_k}{\halfitr-x_\star}$ vanishes under total expectation, 
and the second term can be bounded using \eqref{eqn:DSEG-gap-Gk-hat-moment-estimate}, \eqref{eqn:DESG-gap-cross-term-bound} and \eqref{eqn:DSEG-distance-expec-bound} as
\begin{align}
\label{eqn:DSEG-gap-residual-bound-corrected}
    \expec{\frac{\gamma_k^2}{2}\sqnorm{\widehat{G}_k}
    -\gamma_k\extstep\inprod{\widehat{G}_k}{G_k}}{}
    =
    \expec{\condexp{\frac{\gamma_k^2}{2}\sqnorm{\widehat{G}_k}
    -\gamma_k\extstep\inprod{\widehat{G}_k}{G_k}}{\cF_k}}{}
    \le
    \expec{C \eta_k \left(1+\sqnorm{x_k-x_\star}\right)}{}
    \le C\eta_k (1+M) .
\end{align}
The final term can be bounded using Cauchy-Schwarz, Jensen inequality, $\condexp{\Delta_k}{\cF_k} = 0$, \eqref{eqn:DSEG-gap-halfitr-distance-estimate} and \eqref{eqn:DSEG-distance-expec-bound}:
\begin{align*}
    \expec{\norm{\sum_{k=0}^N\gamma_k\Delta_k}}{} \le
    \left(\sum_{k=0}^N\gamma_k^2\expec{\sqnorm{\Delta_k}}{}\right)^{1/2}
    \le \pr{AC(1+M) + \sigma^2}^{1/2} \left(\sum_{k=0}^\infty\gamma_k^2\right)^{1/2} .
\end{align*}
Putting altogether, we have
\begin{align*}
    \Gamma_N\expec{\err_{B(y,D)}(\overline{\hat{x}}_N^\gamma)}{}
    &\le
    2D^2 + C(1+M)\sum_{k=0}^N \eta_k
    + D_\star \pr{AC(1+M) + \sigma^2}^{1/2} \left(\sum_{k=0}^\infty\gamma_k^2\right)^{1/2} .
\end{align*}
Letting the right-hand side $C_{y,D}$, which is finite and independent of $N$, we obtain the desired result.
\end{proof}

\section{Almost sure convergence of the last iterate}
\label{section:almost-sure-convergence}

Classical works on the almost-sure iterate convergence of I-SEG with square-summable step-sizes $\gamma_k = \extstep$ relied on an additional condition that the operator is quasi-strictly monotone.
Indeed, without the additional assumption of quasi-strict monotonicity (\cref{assumption:quasi-strict-monotone}), I-SEG fails to achieve almost sure convergence $x_k \to x_\star$ with step-size $\extstep = \gamma_k$, even when the step-size is square-summable.
However, the \emph{double step-size extragradient} \eqref{eqn:DSEG} by \citep{HsiehIutzelerMalickMertikopoulos2020_explore}, 
where $\extstep$ and $\gamma_k$ decay at different rates, allows I-SEG to attain $x_k \to x_\star$ almost surely for general monotone $F$ that is not necessarily quasi-strict monotone,
in the unconstrained domain $\cX = \reals^d$ under the weaker variance assumption \eqref{eqn:generalized-variance-bound}.

In this section, we investigate the analogous results for S-SEG, which are missing in the literature. 
We first show that S-SEG also exhibits the iterate convergence $x_k \to x_\star$ almost surely for quasi-strictly monotone problems with square-summable step-sizes.
Our counterexample from \cref{theorem:unbounded-gap-high-probability-tightness} in the previous section shows that this cannot be extended to general monotone operators.
Furthermore, we show that even the \ref{eqn:DSEG} step-sizes fail to rescue the divergence of S-SEG, unlike it does for I-SEG;
both $x_k$ and the weighted average of the iterates can diverge to infinity,
and therefore, neither almost sure convergence nor gap function convergence can be expected for S-SEG.

\subsection{Convergence under quasi-strict monotonicity}

In the following, we show that S-SEG converges for quasi-strictly monotone problems with square-summable $\gamma_k$ and $\extstep$.

\begin{theorem}
\label{theorem:S-SEG-almost-sure-strict-monotone}
Under Assumptions~\ref{assumption:monotonicity}, \ref{assumption:uniform-lipschitzness}, \ref{assumption:unbiasedness} and \ref{assumption:quasi-strict-monotone}, S-SEG with deterministic $\extstep, \gamma_k > 0$ satisfying $\sum_{k=0}^\infty \gamma_k = \infty$, $\sum_{k=0}^\infty \extstep^2 < \infty$ and $\sum_{k=0}^\infty \gamma_k^2 < \infty$ converges almost surely to a solution to \eqref{eqn:VIP}.
\end{theorem}

\begin{proof}

Let $z\in \cX_\star$ be a solution to \eqref{eqn:VIP}.
Because $\sum_{k=0}^\infty \extstep^2 < \infty$ and $\sum_{k=0}^\infty \gamma_k^2 < \infty$, the sequences $\extstep, \gamma_k$ converge to 0, and in particular upper bounded; choose $\overline{\gamma}$ such that $\gamma_k, \extstep \le \overline{\gamma}$ for all $k\ge 0$.
Now observe that because $x_{k+1} = \proj_{\cX} (x_k - \gamma_k F(\halfitr; \xi_{i_k})$ and $z \in \cX$, by nonexpansiveness of $\proj_{\cX}$,
\begin{align*}
    \sqnorm{x_{k+1} - z} & \le \sqnorm{x_k - \gamma_k F(\halfitr; \xi_{i_k}) - z} \\
    & = \sqnorm{x_k - z} - 2\gamma_k \inprod{F(\halfitr; \xi_{i_k})}{x_k - z} + \gamma_k^2 \sqnorm{F(\halfitr; \xi_{i_k})} \\
    & = \sqnorm{x_k - z} - 2\gamma_k \inprod{F(x_k; \xi_{i_k})}{x_k - z} - 2\gamma_k \inprod{F(\halfitr; \xi_{i_k}) - F(x_k; \xi_{i_k})}{x_k - z} + \gamma_k^2 \sqnorm{F(\halfitr; \xi_{i_k})} .
\end{align*}
The second term becomes $-2\gamma_k \inprod{F(x_k)}{x_k - z}$ under conditional expectation on $\cF_k$.
Next,
\begin{align*}
    - 2\gamma_k \inprod{F(\halfitr; \xi_{i_k}) - F(x_k; \xi_{i_k})}{x_k - z} & \le 2\gamma_k \norm{F(\halfitr; \xi_{i_k}) - F(x_k; \xi_{i_k})} \norm{x_k - z} \\
    & \le 2\gamma_k L \norm{\halfitr - x_k} \norm{x_k - z} \\
    & \le 2\gamma_k L \norm{\extstep F(x_k; \xi_{i_k})} \norm{x_k - z} \\
    & \le \gamma_k \extstep L \left( \sqnorm{F(x_k; \xi_{i_k})} + \sqnorm{x_k - z} \right) \\
    & \le \gamma_k \extstep L \left( (1 + 2L^2) \sqnorm{x_k - z} + \sigma^2 \right) 
\end{align*}
where we again use nonexpansiveness of $\proj_\cX$ for the third inequality, and the last line follows from \cref{proposition:uniform-lipschitzness-implies-variance-bound} with $\sigma^2 = 2\max_{i=1,\dots,n} \sqnorm{F(z; \xi_k)}$.
We also have
\begin{align*}
    \sqnorm{F(\halfitr; \xi_{i_k})} & \le \left( \norm{F(x_k; \xi_{i_k})} + L\norm{\halfitr - x_k} \right)^2 \le \left( \norm{F(x_k; \xi_{i_k})} + L\extstep \norm{F(x_k; \xi_{i_k})} \right)^2 \\
    & \le (1 + \overline{\gamma} L)^2 \sqnorm{F(x_k; \xi_{i_k})} \le 
    (1 + \overline{\gamma} L)^2 \left( 2L^2 \sqnorm{x_k - z} + \sigma^2 \right) .
\end{align*}
Putting altogether, we obtain
\begin{align*}
    \condexp{\sqnorm{x_{k+1} - z}}{\cF_k}
    & \le \left(1 + \gamma_k \extstep L(1+2L^2) + 2\gamma_k^2 L^2 (1+\overline{\gamma}L)^2 \right) \sqnorm{x_k - z} - 2\gamma_k \inprod{F(x_k)}{x_k - z} \\
    & \quad + \left( \gamma_k \extstep L + \gamma_k^2 (1 + \overline{\gamma}L)^2 \right) \sigma^2 .
\end{align*}
Because $\sum_{k=0}^\infty \gamma_k \extstep \le \left( \sum_{k=0}^\infty \gamma_k^2 \right)^{1/2} \left( \sum_{k=0}^\infty \extstep^2 \right)^{1/2} < \infty$, 
we can apply \cref{lemma:almost-sure-sequence-lemma} with
\begin{gather*}
    a_k = \sqnorm{x_k - z}, \qquad b_k = 2\gamma_k \inprod{F(x_k)}{x_k - z} \ge 0, \qquad c_k = \left( \gamma_k \extstep L + \gamma_k^2 (1 + \overline{\gamma}L)^2 \right) \sigma^2, \\
    d_k = \gamma_k \extstep L(1+2L^2) + 2\gamma_k^2 L^2 (1+\overline{\gamma}L)^2 
\end{gather*}
to see that $\sqnorm{x_k - z}$ converges and $\sum_{k=0}^\infty \gamma_k \inprod{F(x_k)}{x_k - z} < \infty$ almost surely.

Take any countable, dense subset $\cD \subseteq \cX_\star$.
Consider the event $\cE$ of probability $1$ on which $\norm{x_k - z}$ and the series $\sum_{k=0}^\infty \gamma_k \inprod{F(x_k)}{x_k - z}$ converge for every $z \in \cD$.
Fix any $z_0 \in \cD$. On $\cE$, $\norm{x_k - z_0}$ converges, and in particular, $x_k$ stays bounded in some closed ball $B(z_0, R)$ for some random variable $R>0$.
Now we show that on $\cE$, $x_k$ has a limit point $x_\star \in \cX_\star$.
Suppose to the contrary---since $\cX_\star \cap B(z_0, R+1)$ is compact, this implies there exists $\delta \in (0, 1)$ and $N \in \mathbb{N}$ such that 
\[
    \mathrm{dist} (x_k, \cX_\star \cap B(z_0, R+1)) \ge \delta 
\]
for all $k \ge N$.
As $x_k \in B(z_0, R)$, the closure $\cC$ of the set $\{x_k \,|\, k \ge N\}$ is compact, and 
$\mathrm{dist}(\cC, \cX_\star) \ge \delta$.
Then, by Assumption~\ref{assumption:quasi-strict-monotone}, the continuous function $x \mapsto \inprod{F(x)}{x-z_0}$ is positive on $\cC$ and thus  attains a minimum $m > 0$ on $\cC$.
This is a contradiction, since then $\inprod{F(x_k)}{x_k - z_0} \ge m$ for all $k \ge N$, which implies $\sum_{k=N}^\infty 2\gamma_k \inprod{F(x_k)}{x_k - z_0} \ge 2m \sum_{k=N}^\infty \gamma_k = \infty$. This shows that $x_k$ must have a limit point $x_\infty \in \cX_\star$.

Finally, it remains to show that in fact, $x_k \to x_\infty$. 
Choose a sequence $z_j \in \cD$ such that $z_j \to x_\infty$.
On $\cE$, the limits $e_j = \lim_{k\to\infty} \|x_k - z_j\|$ exist for all $j$.
Hence, for any $\epsilon > 0$, take sufficiently large $j$ for which $\norm{z_j - x_\infty} < \epsilon$.
Then, taking $k\to \infty$ in $\left| \norm{x_k - x_\infty} - \norm{x_k - z_j} \right| \le \norm{z_j - x_\infty}$, we obtain
\begin{align*}
    e_j - \epsilon \le \liminf_{k\to\infty} \norm{x_k - x_\infty} \le \limsup_{k\to\infty} \norm{x_k - x_\infty} \le e_j + \epsilon .
\end{align*}
Since $\epsilon$ was arbitrary, this shows that $\lim_{k\to\infty} \norm{x_k - x_\infty}$ exists, and because $x_\infty$ is a limit point of $x_k$, that limit must be $0$.
This shows that $x_k \to x_\infty$ on $\cE$, completing the proof.

\end{proof}

\subsection{S-SEG may diverge even with DSEG step-sizes}
\label{subsection:DSEG-fails-to-fix-S-SEG}

We now turn to the question of what happens without quasi-strict monotonicity.
\cref{theorem:unbounded-gap-high-probability-tightness} shows that we cannot expect the iterate convergence $x_k \to x_\star \in \cX_\star$ with symmetric step-sizes.
However, because \ref{eqn:DSEG} step-sizes resolves this issue by breaking the symmetry, one may naturally hypothesize that S-SEG may exhibit almost-sure convergence with DSEG step-sizes.
Below, we show that this is in fact not the case, by presenting a counterexample.

\begin{theorem}
\label{theorem:same-sample-DSEG-gamma-average-divergence}
For every fixed $0<\eta<\frac{1}{10}$,
there exists a stochastic VIP on $\cX=\reals^2$ satisfying 
Assumptions~\ref{assumption:monotonicity}, \ref{assumption:uniform-lipschitzness}
and \ref{assumption:unbiasedness} with unique solution $x_\star=0$, and a sufficiently large $k_0 > 1$ such that for every $x_0 \ne 0$, S-SEG with step-sizes
$\extstep = \frac{1}{\log(k+k_0)}, \gamma_k=\frac{\eta}{k+k_0}$, 
which satisfy \eqref{eqn:DSEG}, exhibits $\norm{x_k}\to\infty$ almost surely.
Furthermore, denoting $\Gamma_N=\sum_{k=0}^N\gamma_k$ and $\overline{\hat{x}}_N^\gamma=\frac{1}{\Gamma_N}\sum_{k=0}^N \gamma_k \hat{x}_k$,
we have $\norm{\overline{\hat{x}}_N^\gamma}\to\infty$ almost surely as $N\to\infty$.
\end{theorem}

\begin{proof}
Fix $0<\eta<1/10$, and choose an integer $\rho\ge2$ such that
$q:=\eta(\rho^2-1)>1$. Let $L=\sqrt{\rho^2+1}$.
Choose an integer $k_0>1$ sufficiently large that
\[
    \frac{1}{\log k_0}<\frac{1}{3L}
    , \qquad
    \sup_{k\ge0,\;s\in\{-1,+1\}}
    |b_s(\gamma_k)-1|\le\frac12.
\]
The step-sizes satisfy the \eqref{eqn:DSEG} conditions
\[
    \sum_k\gamma_k\hat\gamma_k=\infty,
    \qquad
    \sum_k\gamma_k^2<\infty,
    \qquad
    \sum_k\gamma_k\hat\gamma_k^2<\infty,
\]
because
\[
    \gamma_k\hat\gamma_k
    =
    \frac{\eta}{(k+k_0)\log(k+k_0)},
    \quad
    \gamma_k^2
    =
    \frac{\eta^2}{(k+k_0)^2},
    \quad
    \gamma_k\hat\gamma_k^2
    =
    \frac{\eta}{(k+k_0)\log^2(k+k_0)}.
\]
Now consider the stochastic VIP from \cref{lemma:stochastic-vip-counterexample} with $\nu=0$ and $\rho$.
By \cref{lemma:stochastic-vip-counterexample}, this problem satisfies the required assumptions with $F(x) = Jx$ and has unique solution $x_\star=0$.
In addition, with $s_k\in\{-1,+1\}$ sampled uniformly, we have
\[
    \hat z_k = \left(1-\extstep(s_k\rho+\imu)\right)z_k,
    \qquad
    z_{k+1}=b_{s_k}(\gamma_k)z_k,
\]
where $b_{s_k}(\gamma_k) = 1+\gamma_k\bigl((\rho^2-1)\extstep-s_k\rho\bigr)+\imu\gamma_k\bigl(2s_k\rho \extstep-1\bigr)$.
Because $|b_{s_k}(\gamma_k) - 1| \le \frac{1}{2}$ for all $k\ge 0$ and any selection of $s_k$, the principal complex logarithm of $b_{s_k}(\gamma_k)$ is well-defined.
Now using the series expansion $\operatorname{Log} (1+u) = u - \frac{u^2}{2} + \frac{u^3}{3} - \cdots$ for $|u|<1$, we have
\begin{align}
\label{eqn:S-SEG-DSEG-log-expansion}
    \operatorname{Log} b_{s_k}(\gamma_k)
    =
    - s_k\rho\gamma_k
    - \imu\gamma_k
    + (\rho^2-1)\gamma_k \extstep
    + 2\imu s_k\rho\gamma_k \extstep
    + \cO(\gamma_k^2) .
\end{align}
Because we have $\sum_{k=0}^\infty \gamma_k^2<\infty$ and $\sum_{k=0}^\infty \gamma_k^2\extstep^2<\infty$, by Kolmogorov's two-series theorem, the random series 
$\sum_{k=0}^\infty s_k\gamma_k$ and $\sum_{k=0}^\infty s_k\gamma_k \extstep$ converge almost surely.
Moreover, we have 
\[
    \sum_{k=0}^{N-1}\gamma_k \extstep
    =
    \eta\log\log N + c_1 + o(1),
    \qquad
    \sum_{k=0}^{N-1}\gamma_k
    =
    \eta\log N + c_2 + o(1)
\]
for some deterministic constants $c_1, c_2$ depending on $k_0$.
The remainder $\cO(\gamma_k^2)$ in
\eqref{eqn:S-SEG-DSEG-log-expansion} is uniform in $s_k$ and summable.
Therefore, summing \eqref{eqn:S-SEG-DSEG-log-expansion} from $k=0$ to $N-1$, 
we see that there exists an almost surely finite nonzero complex random variable $C$ such that
\[
    z_N
    =
    C(\log N)^q N^{-\imu\eta}(1+o(1))
\]
where $q=\eta(\rho^2-1) > 1$. In particular, this implies $\norm{x_N}=|z_N|\to\infty$ almost surely.

To analyze the averaged iterates, we need an estimate of the quantity $T_N := \sum_{k=0}^N \gamma_k \extstep z_k$.
Rearranging the recursion
\begin{align*}
    z_{k+1}-z_k
    =
    -\imu\gamma_k z_k
    +
    \gamma_k\left((\rho^2-1)\extstep-s_k\rho+2\imu s_k\rho \extstep\right)z_k 
\end{align*}
we obtain
\begin{align}
\label{eqn:S-SEG-DSEG-gamma-k-zk-expression}
    \gamma_k z_k = \imu (z_{k+1} - z_k)
    -
    \imu\gamma_k\left((\rho^2-1)\extstep-s_k\rho+2\imu s_k\rho \extstep\right)z_k .
\end{align}
Now multiplying the above equation by $\extstep$ and summing gives
\begin{align}
\label{eqn:S-SEG-DSEG-TN-expansion}
    T_N
    &=
    \imu\sum_{k=0}^N \extstep(z_{k+1}-z_k)
    -\imu(\rho^2-1)\sum_{k=0}^N \gamma_k \extstep^2 z_k
    +\imu\rho\sum_{k=0}^N s_k\gamma_k \extstep z_k
    +2\rho\sum_{k=0}^N s_k\gamma_k \extstep^2 z_k .
\end{align}
Observe that we have $|\hat\gamma_{k+1}-\extstep| = \cO \pr{\frac{1}{k(\log k)^2}}$ and we already obtained $|z_k| = \cO((\log k)^q)$ almost surely.
Hence, by summation by parts we have
\begin{align*}
    \left| \sum_{k=0}^N \extstep(z_{k+1}-z_k) \right| & = \left| -\hat\gamma_0 z_0 + \hat\gamma_N z_{N+1} + \sum_{k=0}^{N-1} \left(\hat\gamma_k - \hat\gamma_{k+1}\right) z_{k+1}  \right| \\
    & = \cO(1) + |\hat{\gamma}_N| |z_{N+1}| + \sum_{k=0}^{N-1} \left| \hat{\gamma}_k - \hat{\gamma}_{k+1} \right| \cdot |z_{k+1}| \\
    & = \cO\pr{(\log N)^{q-1}} + \cO\pr{ \sum_{k=0}^{N-1} \frac{(\log k)^q}{k (\log k)^2}} = \cO\pr{(\log N)^{q-1}} 
\end{align*}
almost surely.
We also have
\[
    \left| \sum_{k=0}^N \gamma_k \extstep^2 z_k \right|
    \le
    \sum_{k=0}^N \gamma_k \extstep^2 |z_k|
    =
    \cO\left(\sum_{k=2}^N \frac{(\log k)^q}{k (\log k)^2}\right)
    =
    \cO\pr{(\log N)^{q-1}} 
\]
almost surely.
Finally, the two random series $\sum_{k=0}^\infty s_k\gamma_k \extstep z_k$ and $\sum_{k=0}^\infty s_k\gamma_k \extstep^2 z_k$ in \eqref{eqn:S-SEG-DSEG-TN-expansion} converge almost surely because $\condexp{s_k}{\cF_k} = 0$ makes them martingale and
\[
    \sum_{k=0}^\infty \gamma_k^2\extstep^2|z_k|^2 = \cO\pr{\sum_k \frac{(\log k)^{2q}}{k^2 (\log k)^2}} < \infty,
    \qquad
    \sum_{k=0}^\infty \gamma_k^2\extstep^4|z_k|^2 = \cO\pr{\sum_k \frac{(\log k)^{2q}}{k^2 (\log k)^4}} < \infty 
\]
almost surely. 
This shows that $T_N=\cO\pr{(\log N)^{q-1}} = o\pr{(\log N)^q}$.
Now summing the identity \eqref{eqn:S-SEG-DSEG-gamma-k-zk-expression} and telescoping we obtain
\begin{align*}
    \sum_{k=0}^N \gamma_k z_k
    &=
    \imu(z_{N+1}-z_0)
    -\imu(\rho^2-1) T_N 
    +\imu\rho\sum_{k=0}^N s_k\gamma_k z_k
    +2\rho\sum_{k=0}^N s_k\gamma_k \extstep z_k .
\end{align*}
The last two series $\sum_{k=0}^\infty s_k\gamma_k z_k$ and $\sum_{k=0}^\infty s_k\gamma_k \extstep z_k$ 
converge almost surely because
\[
    \sum_{k=0}^\infty \gamma_k^2|z_k|^2 = \cO\pr{\sum_k \frac{(\log k)^{2q}}{k^2}} < \infty,
    \qquad
    \sum_{k=0}^\infty \gamma_k^2\extstep^2|z_k|^2 = \cO\pr{\sum_k \frac{(\log k)^{2q}}{k^2 (\log k)^2}} <\infty .
\]
Because we have $z_N = C(\log N)^q N^{-\imu \eta} (1 + o(1))$ and $T_N=o((\log N)^q)$, this yields
\[
    \sum_{k=0}^N \gamma_k z_k
    =
    \imu z_{N+1} + o((\log N)^q)
    =
    \imu C(\log N)^q N^{-\imu\eta} (1+o(1)).
\]
Finally, using $\hat z_k = \left(1-\extstep(s_k\rho+\imu)\right)z_k$ we have 
\[
    \sum_{k=0}^N \gamma_k \hat z_k
    =
    \sum_{k=0}^N \gamma_k z_k
    -
    \rho\sum_{k=0}^N s_k\gamma_k \extstep z_k
    -
    \imu\sum_{k=0}^N \gamma_k \extstep z_k ,
\]
and we have already checked that the second series converges almost surely, while the last series is $\imu T_N=o(L_N^q)$. 
Therefore, obtain
\[
    \sum_{k=0}^N \gamma_k \hat z_k
    =
    \imu C(\log N)^q N^{-\imu\eta}(1+o(1)).
\]
Because $\Gamma_N=\sum_{k=0}^N\gamma_k=\eta\log N+\cO(1)$, we conclude that
\[
    \overline{\hat z}_N^\gamma
    :=
    \frac{1}{\Gamma_N}\sum_{k=0}^N \gamma_k\hat z_k
    =
    \frac{\imu C}{\eta}(\log N)^{q-1}N^{-\imu\eta}(1+o(1))
\]
which implies $\norm{\overline{\hat x}_N^\gamma} = |\overline{\hat z}_N^\gamma| \to \infty$ almost surely, since $q > 1$.
\end{proof}

\section{Conclusion}

We present a convergence theory for SEG on stochastic monotone VIPs, providing results that have been missing from the literature despite the algorithm's simplicity and popularity.
Our analysis highlights the differences between I-SEG and S-SEG, both in the assumptions required for convergence and their behavior across distinct problems and step-size choices.
We further characterize the dynamics of SEG when the domain and operator variance are unbounded.

While our results address several gaps in the literature on SEG, they also raise interesting questions for future work.
A similar analysis may be possible for the single-oracle variant of EG, now commonly called optimistic gradient (OG) \cite{rakhlinOptimizationLearningGames2013}, which dates back to the work of \citet{popovModificationArrowHurwiczMethod1980}.
Prior work on stochastic OG (SOG) assumed uniformly bounded variance \cite{cuiAnalysisReflectedGradient2016,hsiehConvergenceSinglecallStochastic2019,GidelBerardVignoudVincentLacoste-Julien2019_variational,bohmSolvingNonconvexnonconcaveMinmax2023} or large batch sizes depending on the desired accuracy level \cite{bohmSolvingNonconvexnonconcaveMinmax2023,ChoudhuryGorbunovLoizou2023_singlecall},
and it would be informative to determine whether SOG behaves similarly to SEG on stochastic monotone VIPs with unbounded variance over general closed convex domains.
Another open direction is to identify a simple mechanism that induces convergence under the generalized variance condition~\eqref{eqn:generalized-variance-bound}.
Recent methods use regularization or variance reduction to achieve convergence in similar settings \cite{neuDealingUnboundedGradients2024,alacaogluWeakerVarianceAssumptions2025,alacaogluSolvingStochasticVariational2026},
and it remains to determine whether these ideas can be distilled into a minimal update rule comparable to SEG.
Finally, although we have shown that the \ref{eqn:DSEG} step-size rule of \citet{HsiehIutzelerMalickMertikopoulos2020_explore} resolves the non-convergence issue for I-SEG but not for S-SEG, it remains an open question whether the non-convergence of S-SEG is unavoidable or can be remedied by another carefully chosen step-size rule.

\section*{Acknowledgments}
The authors thank Sayantan Choudhury and Ahmet Alacaoglu for useful discussions during the early stages of this project.
TaeHo Yoon’s contribution to this work was supported by NSF CCF 2504626.
Nicolas Loizou’s contribution to this work was supported by NSF CCF 2504626 and NSF CAREER 2542902.

\bibliographystyle{plainnat}
\bibliography{ref}

\end{document}